%% file: main.tex
\documentclass{amsart} 
\usepackage{graphicx}  
\usepackage{hyperref}
\usepackage{amsmath, amsthm, amssymb, amscd}
\usepackage{mathrsfs}
\usepackage{mathtools}
\usepackage{enumerate}
\usepackage{hyperref}
\usepackage{verbatim}
\usepackage{subcaption}
\usepackage{cite}
\usepackage{color}
\usepackage{esint}
\usepackage{tikz-cd}
\usetikzlibrary{shapes.geometric}
\usepackage[shortlabels]{enumitem}
 
\usepackage{bm}
\usepackage{bbm, dsfont}
\definecolor{mycolor}{rgb}{0.0, 0.75, 1.0}

\makeatletter
\newcommand{\labitem}[2]{%
\def\@itemlabel{\textbf{#1}}
\item
\def\@currentlabel{#1}\label{#2}}
\makeatother

\title[Construction of energy forms on Laakso-type spaces]{Construction of self-similar energy forms and singularity of Sobolev spaces on Laakso-type fractal spaces}
\author{Riku Anttila}
\address[Riku Anttila]{Department of Mathematics and Statistics, University of Jyväskylä, P.O. Box 35, FI-40014 Jyväskylä, Finland}
\email{riku.t.anttila@jyu.fi}

\author{Sylvester Eriksson-Bique}
\address[Sylvester Eriksson-Bique]{Department of Mathematics and Statistics, University of Jyväskylä, P.O. Box 35, FI-40014 Jyväskylä, Finland}
\email{sylvester.d.eriksson-bique@jyu.fi}

\author{Ryosuke Shimizu}
\address[Ryosuke Shimizu]{Waseda Research Institute for Science and Engineering, Waseda University, 3-4-1 Okubo, Shinjuku-ku, Tokyo 169-8555, Japan.}
\curraddr{Graduate School of Informatics, Kyoto University, Yoshida-honmachi, Sakyo-ku, Kyoto 606-8501, Japan.}
\email{r.shimizu@acs.i.kyoto-u.ac.jp}

\thanks{}
\subjclass[2020]{28A80, 31E05, 31C45, 46E36}
\keywords{Laakso-type fractals, (edge)-iterated graph systems, energy forms, energy measures, Sobolev spaces, singularity} 
\date{\today}

\usepackage[shortlabels]{enumitem}
 
\newtheorem{theorem}[equation]{Theorem}
\newtheorem{lemma}[equation]{Lemma}

\newtheorem{proposition}[equation]{Proposition}
\newtheorem{corollary}[equation]{Corollary}
\newtheorem{assumption}[equation]{Assumption}

\newtheorem{question}[equation]{Question}
\numberwithin{equation}{section}

\theoremstyle{definition}
\newtheorem{definition}[equation]{Definition}

\theoremstyle{remark}
\newtheorem{remark}[equation]{Remark}
\newtheorem{example}[equation]{Example}
    \newcommand*{\N}{\mathbb{N}}
    \newcommand*{\R}{\mathbb{R}}
        \DeclarePairedDelimiter\Span{\langle}{\rangle}
        \DeclareMathOperator{\id}{id}
        \DeclareMathOperator{\diam}{diam}
        \DeclareMathOperator{\Fib}{Fib}
        \DeclareMathOperator{\dist}{dist}
        \newcommand{\pwalk}{d_{w,p}}
        
        \DeclareMathOperator{\divr}{div}
        \DeclareMathOperator{\len}{len}
        \DeclareMathOperator{\cCap}{\rm Cap}
        \DeclareMathOperator{\cRes}{\rm Res}
        \DeclareMathOperator{\cE}{\mathcal{E}}
        \DeclareMathOperator{\fn}{\mathfrak{n}}
        \DeclareMathOperator{\fe}{\mathfrak{e}}

        \DeclareMathOperator{\cM}{\mathcal{M}}
        
        \DeclareMathOperator{\supp}{supp}
        \DeclareMathOperator{\Lip}{Lip}
        
        \DeclareMathOperator{\sgn}{sgn}

        \DeclareMathOperator{\restr}{\upharpoonright}
        \DeclarePairedDelimiter\abs{\lvert}{\rvert}
        \DeclarePairedDelimiter\norm{\lVert}{\rVert}
        \def\vint_#1{\mathchoice%
          {\mathop{\kern 0.2em\vrule width 0.6em height 0.69678ex depth -0.58065ex
                  \kern -0.8em \intop}\nolimits_{\kern -0.4em#1}}%
          {\mathop{\kern 0.1em\vrule width 0.5em height 0.69678ex depth -0.60387ex
                  \kern -0.6em \intop}\nolimits_{#1}}%
          {\mathop{\kern 0.1em\vrule width 0.5em height 0.69678ex
              depth -0.60387ex
                  \kern -0.6em \intop}\nolimits_{#1}}%
          {\mathop{\kern 0.1em\vrule width 0.5em height 0.69678ex depth -0.60387ex
                  \kern -0.6em \intop}\nolimits_{#1}}}
\def\vintslides_#1{\mathchoice%
          {\mathop{\kern 0.1em\vrule width 0.5em height 0.697ex depth -0.581ex
                  \kern -0.6em \intop}\nolimits_{\kern -0.4em#1}}%
          {\mathop{\kern 0.1em\vrule width 0.3em height 0.697ex depth -0.604ex
                  \kern -0.4em \intop}\nolimits_{#1}}%
          {\mathop{\kern 0.1em\vrule width 0.3em height 0.697ex depth -0.604ex
                  \kern -0.4em \intop}\nolimits_{#1}}%
          {\mathop{\kern 0.1em\vrule width 0.3em height 0.697ex depth -0.604ex
                  \kern -0.4em \intop}\nolimits_{#1}}}

\newcommand{\kint}{\vint}

\newcommand{\aveint}[2]{\mathchoice%
          {\mathop{\kern 0.2em\vrule width 0.6em height 0.69678ex depth -0.58065ex
                  \kern -0.8em \intop}\nolimits_{\kern -0.45em#1}^{#2}}%
          {\mathop{\kern 0.1em\vrule width 0.5em height 0.69678ex depth -0.60387ex
                  \kern -0.6em \intop}\nolimits_{#1}^{#2}}%
          {\mathop{\kern 0.1em\vrule width 0.5em height 0.69678ex depth -0.60387ex
                  \kern -0.6em \intop}\nolimits_{#1}^{#2}}%
          {\mathop{\kern 0.1em\vrule width 0.5em height 0.69678ex depth -0.60387ex
                  \kern -0.6em \intop}\nolimits_{#1}^{#2}}}
\usepackage{systeme}
\makeatletter
\let\c@equation\c@figure
\makeatother

\begin{document}

\begin{abstract}
We construct self-similar $p$-energy forms $\mathscr{E}_p$ on a rich class of \emph{Laakso-type fractal spaces} and study the properties of the associated Sobolev spaces $\mathscr{F}_p$.
The main result of the paper is the discovery of a new analytic phenomenon, which we refer to as \emph{singularity of Sobolev spaces}. 
This means that the associated Sobolev spaces $\mathscr{F}_{p_1}$ and $\mathscr{F}_{p_2}$ for distinct $p_1,p_2 \in (1,\infty)$ intersect only at constant functions.
We show that the Laakso diamond space of Lang--Plaut is one such example, and explain why this does not contradict the inverse limit construction of Cheeger--Kleiner which proves Laakso Diamond to support Poincar\'e inequality of Heinonen--Koskela.
\end{abstract}

\maketitle

\section{Introduction}
In this paper, we introduce a unified and self-contained framework of analysis on a certain class of \emph{Laakso-type fractal spaces} that was recently introduced by the first two authors \cite{anttila2024constructions}.
Let us begin the paper by describing the constructions and some results regarding one particular example, the \emph{Laakso diamond space}, which is a well-known example in embeddings and analysis on metric spaces communities \cite{AndreaSchioppa2015,LangPlaut,LeeNaor,David_Schul_2017,CK_PI}.
\begin{figure}[h!]
    \centering\includegraphics[width=260pt]{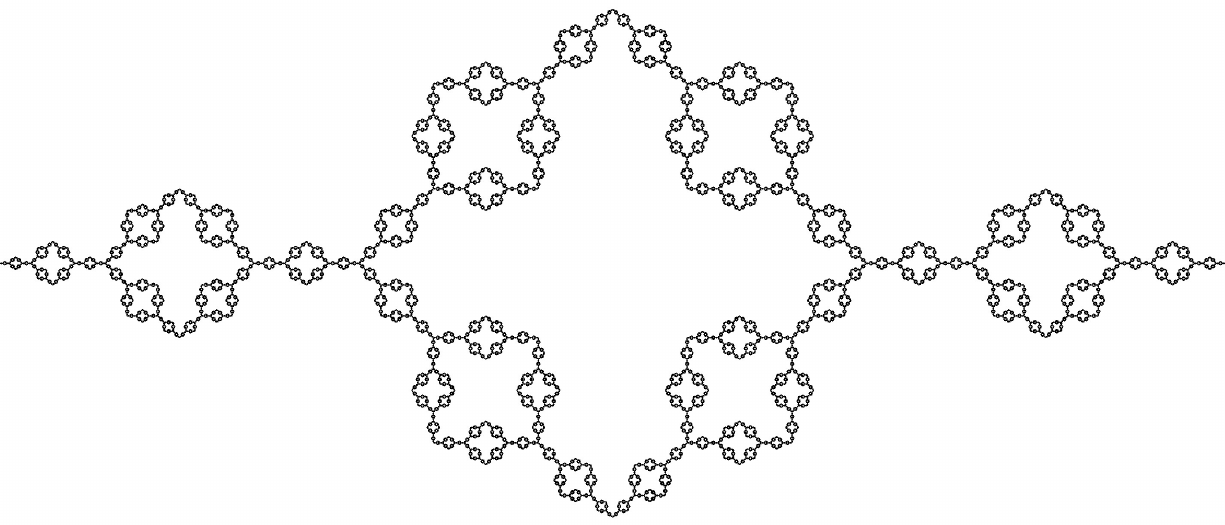} 
    \caption{The Laakso diamond space snowflake-embedded onto a self-similar set of $\R^2$.
    }
    \label{Fig: Laakso diamond limit}
\end{figure}

\begin{example}\label{intro: Construction example}
The Laakso diamond space can be constructed as a limit of graphs, and the idea is understood best from Figures \ref{Fig: Laakso Diamond replacement} and \ref{fig: Laakso G_3}.
\begin{figure}[!ht]
\centering\includegraphics[width=290pt]{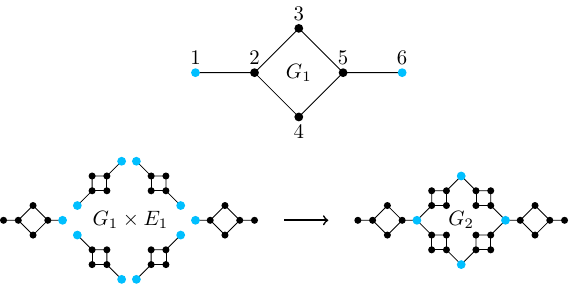}
\caption{The first replacement in the construction of Laakso diamond space.}
\label{Fig: Laakso Diamond replacement}
\end{figure}
The resulting object is a metric measure space $(X,d,\mu)$ which is geodesic and $Q$-Ahlfors regular for $Q := \log(6)/\log(4)$. Recall that $Q$-Ahlfors regularity means $\mu(B(x,r)) \approx r^Q$ for all $x \in X$ and $r \in (0,2\diam(X))$.
The construction goes as follows. Let $G_1$ be the diamond graph as in Figure \ref{Fig: Laakso Diamond replacement}.
We construct the graph $G_2$ by first replacing each edge in $G_1$ by $G_1$ and get the disconnected graph $G_1 \times E_1$. 
Then we identify certain vertices according to a predetermined rule which are indicated by \textcolor{mycolor}{blue colorings} in Figure \ref{Fig: Laakso Diamond replacement}. This produces the graph $G_2$. Next we construct the graph $G_3$ in Figure \ref{fig: Laakso G_3} by first replacing each edge in $G_2$ by $G_1$, creating a disconnected graph $G_1 \times E_2$, and then identifying the vertices in an analogous manner.
By following this pattern, we have an infinite sequence of graphs $\{G_n\}_{n \in \N}$. Finally, we give each graph $G_n$ a rescaled path-metric $4^{-n}d_{G_n}$ and the limiting object $(X,d,\mu)$ is obtained by taking a Gromov-Hausdorff limit of $G_n$ as $n \to \infty$. The scaling constants $4^{-n}$ are natural because they preserve the diameters equal to $1$.
\begin{figure}[!ht]
\centering\includegraphics[width=250pt]{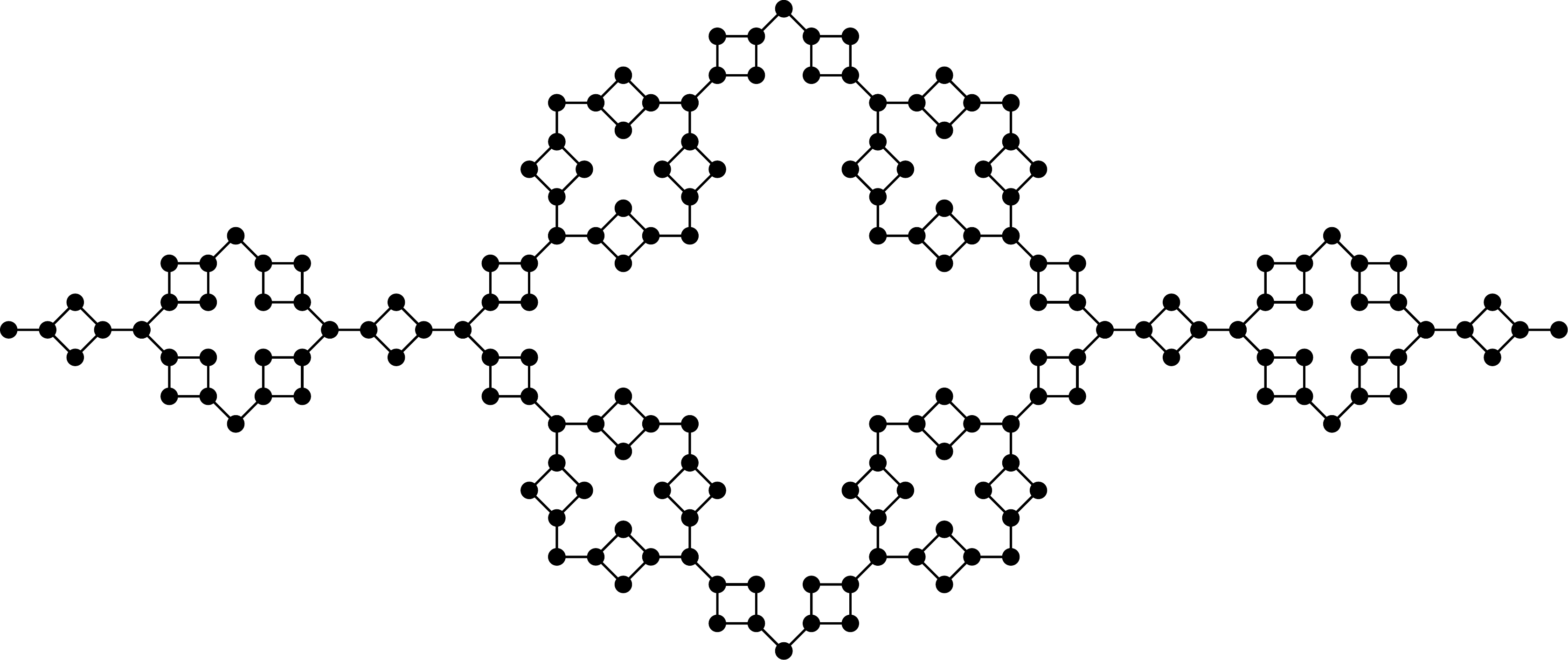}
\caption{The second replacement in the construction of the Laakso diamond space.}
\label{fig: Laakso G_3}\end{figure}
\end{example}

Our general construction of Laakso-type fractals generalizes the construction above, and the details can be found in Section \ref{sec:ssLaakso}.
The methods are strongly inspired by Laakso's work \cite{Laakso,La00} and the inverse limit construction of Cheeger--Kleiner \cite{CK_PI}, which produce metric measure spaces satisfying strong analytic properties, most notably a Poincar\'e inequality of Heinonen--Koskela \cite{HK}; see the end of Introduction for further related discussion.
However, a major motivation for this work, as well as for the previous study \cite{anttila2024constructions}, lies in the fact that our construction also produces metric spaces with weaker analytic properties compared those in \cite{Laakso,La00,CK_PI}.
This observation allows for more flexible constructions, and in \cite{anttila2024constructions}, led to a negative resolution of Kleiner's conjecture \cite[Conjecture 7.5]{KleinerICM}. For further details, see \cite[Theorem 1.3]{anttila2024constructions} and Example \ref{example: Counterexample}.

Laakso's work has influenced many other later studies as well. For instance, see Barlow \cite{Barlow04}, Barlow--Evans \cite{barlow2004markov} which was later studied in \cite{steinhurst2010diffusions,alonso2018explicit,ruiz2021heat}, a recent work of Murugan \cite{murugan2024diffusionsrandomwalksprescribed} and its extension \cite{YangLaakso2025}, studies on differentiability on Laakso spaces \cite{UnivDiffLaakso,capolli2025maximal} and the other references above.
See also \cite{LaaksoOverview} for a nice overview on Laakso spaces.

\subsection{Construction of energy forms}
The main objective of this work is to construct natural  \emph{$p$-energies} $\mathscr{E}_p : L^p(X,\mu) \to [0,\infty]$ for $p \in (1,\infty)$ on Laakso-type fractal spaces. The role of $\mathscr{E}_p$ is a counterpart of the Dirichlet $p$-energy $f \mapsto \norm{\nabla f}^p_{L^p(\R^n)} $. We are interested in the associated \emph{Sobolev spaces}
\[
    \mathscr{F}_p := \{ f \in L^p(X,\mu) : \mathscr{E}_p(f) < \infty \}.
\]
The construction of $\mathscr{E}_p$ is quite easy to describe on the Laakso diamond, and we have sketched it below.
The same approach can be used to reconstruct the usual Sobolev space on the unit interval; see \cite[Introduction]{kigami}.
See also \cite{murugan2023first,shimizu,herman2004p,cao2022p,AnalOnFractals,KusuokaZhou,Baudoin_Chen_2023,scaleirregular,kigamiota} for related techniques on other fractal spaces.

Let $(X,d,\mu)$ be the Laakso diamond space and $\{G_n\}_{n \in \N}$ be the sequence of graphs as in Example \ref{intro: Construction example}. The vertex sets and edge sets of $G_n$ are denoted $V_n$ and $E_n$, respectively.
First, observe from the figures above that we have natural inclusions $V_n \subseteq X$ for all $n \in \N$. Given $f : X \to \R$ we can define the natural discrete energies
\[
    \mathcal{E}_p^{G_n}(f) := \sum_{\{ v,w \} \in E_n} \abs{f(v) - f(w)}^p.
\]
The desired $p$-energy is then obtained by taking the normalized limit
\[
    \mathscr{E}_p(f) := \lim_{n \to \infty} \cM_p^{-n}\mathcal{E}_p^{G_n}(f) \in [0,\infty].
\]
The constant $\cM_p \in (0,\infty)$ is given explicitly as the value of the discrete optimization problem on $G_1$,
\begin{equation}\label{intro:Capacity}
    \cM_p := \min_{U} \sum_{ \{v,w\} \in E_1 } \abs{U(v) - U(w)}^p.
\end{equation}
where the minimum is taken over all functions $U : V_1 \to \R$ satisfying the boundary conditions $U(v_+) = 1$ and $U(v_-) = 0$. Here $v_+$ and $v_-$ are the left most and right most vertices in $G_1$, respectively. A minimizer $U$ exists and is unique by the convexity properties of the $p$-energy.
The solution is computed in Figure \ref{intro: potential/flow}.
By this choice of $\cM_p$, the sequence of normalized energies $\cM_p^{-n}\mathcal{E}_p^{G_n}(f)$ is non-decreasing in $n \in \N$ for all $f : X\to \R$. Hence $\mathscr{E}_p(f) \in [0,\infty]$ is always well-defined.
\begin{figure}[!ht]
    \begin{tikzpicture}[scale=0.55]
    \input{plots/capacity_laakso}
    \end{tikzpicture} 
    \caption{The values in the figure indicate the unique function $U$ that solves the optimization problem \eqref{intro:Capacity}.}
    \label{intro: potential/flow}
\end{figure}

Details on our general constructions can be found in Sections \ref{sec: DPT}-\ref{Sec: Sobolev spaces}.
Many of our methods, especially the averaging procedure used to prove the Poincar\'e inequality in Section \ref{sec: DPT}, are strongly influenced by Cheeger--Kleiner \cite[Section 5]{CK_PI}.

\subsection{Singularity of Sobolev spaces}
The constructed $p$-energies $\mathscr{E}_p$ enjoy a comprehensive analytic structure. We have for instance a Poincar\'e type inequality, natural capacity bounds, contraction properties, Clarkson's inequality, density of continuous functions and so forth. See Sections \ref{Sec: Sobolev spaces}-\ref{sec:otherSob} for proofs and precise statements. Our energies are also \emph{self-similar}, meaning they solve an equation which respects the canonical dynamics of the underlying fractal space; see Theorem \ref{Thm: Existence of p-Energy}-\eqref{item: p-energySS}.
The study of self-similar energies is a core concept in analysis on fractals; we refer to the survey for references \cite{KajinoSGsurvey}.

We also show that these energies are comparable to a \emph{Korevaar--Schoen} type geometric energies, which are quite common and important in the modern literature of analysis on general metric spaces. See for instance \cite{kajino2024korevaarschoenpenergyformsassociated,Baudoin_2024,shimizu2024characterizationssobolevfunctionsbesovtype,HKST}.
On the Laakso diamond, the Korevaar--Schoen formula is written 
\begin{align}\label{intro:Korevaar-Schoen}
    C^{-1}\mathscr{E}_p(f) & \leq \liminf_{r \to 0^+} \int_X \kint_{B(x,r)} \frac{\abs{f(x) - f(y)}^p}{r^{d_{w,p}}}\, \mu(dy)\mu(dx)\\
    & \leq \hspace{12pt}\sup_{r > 0} \int_X \kint_{B(x,r)} \frac{\abs{f(x) - f(y)}^p}{r^{d_{w,p}}}\,\mu(dy)\mu(dx) \leq C\mathscr{E}_p(f) \nonumber
\end{align}
where $d_{w,p} := \log\left(6\cM_p^{-1}\right)/\log(4) > p$ and $\mu$ is the natural Ahlfors regular $Q$-Hausdorff measure for $Q := \log(6)/\log(4)$ with respect to the geodesic metric; see Section \ref{sec:otherSob} for details.

We have introduced all the basic ideas of the work. Now, we state the main result of the work which is the discovery of the following rather surprising behavior of the Laakso diamond space.
The methods heavily rely on the description of $\mathscr{E}_p$, and the Korevaar-Schoen formula does not seem to be explicit enough to this end. The proof is sketched below and see Theorem \ref{thm:SingularitySobolev} for a more general result.

\begin{theorem}[Singularity of Sobolev spaces]\label{intro:Thm}
    Let $(X,d,\mu)$ be the Laakso diamond space as in Example \ref{intro: Construction example}.
    For a given $p \in (1,\infty)$, let
    $\mathscr{E}_p : L^p(X,\mu) \to [0,\infty]$ be the self-similar $p$-energy form described above and denote the associated Sobolev space
    \[
        \mathscr{F}_p := \{ f \in L^p(X,\mu) : \mathscr{E}_p(f) < \infty \}.
    \]
    Then for all distinct $p_1,p_2 \in (1,\infty)$,
    \[
        \mathscr{F}_{p_1} \cap \mathscr{F}_{p_2} = \{ f \in L^1(X,\mu) : f \text{ is constant $\mu$-almost everywhere} \}.
    \]
\end{theorem}

Namely, the Sobolev spaces $\mathscr{F}_p$ for distinct $p \in (1,\infty)$ on the Laakso diamond are essentially disjoint.
We refer to this analytic phenomenon as \emph{singularity of Sobolev spaces}, emphasizing their disjointness. A partial reason for this terminology is to reflect another phenomenon, the \emph{singularity of energy measures}, which was recently discovered by Kajino and third author \cite{KajinoSingularityPCF,KajinoSGsurvey}; see Proposition \ref{intro:SingularityEM} for analogous result on the Laakso diamond and Proposition \ref{prop:SingularityEM} in the general setting.

We also verify another curious property of the Laakso diamond space regarding its Lipshchitz functions. See Theorem \ref{thm:NoLips} for the general version.

\begin{theorem}\label{intro:NoLips}
    For all $p \in (1,\infty)$ the Sobolev space $\mathscr{F}_p$ on the Laakso diamond space $(X,d,\mu)$ contains no non-constant Lipschitz functions. In other words,
    \[
        \mathscr{F}_p \cap \Lip(X,d) = \{ f \in L^1(X,\mu) : f \text{ is constant $\mu$-almost everywhere} \}.
    \]
\end{theorem}

Now that the theorems are stated, let us elaborate why we find them interesting. We have outlined our thoughts below into three bullet points.

First, they give clear qualitative distinctions between the more classical forms of analysis studied by the analysis on metric spaces community; see the standard references \cite{He,HKST,bjorn2011nonlinear}. This is because much of the related literature heavily employs Lipschitz functions, and most Sobolev type spaces studied therein contain the compactly supported Lipschitz functions.

Second, what we found particularly surprising is that these behaviors are realized on the Laakso diamond space. On one hand, Laakso diamond is known to carry one the foundational inequalities of analysis on metric spaces, the Poincar\'e inequality of Heinonen--Koskela \cite{HK}. This follows from Cheeger--Kleiner \cite{CK_PI}. But on the other hand, the validity of (Heinonen--Koskela) Poincar\'e inequality indicates that both singularity of Sobolev spaces and Theorem \ref{intro:NoLips} should be false.
This seeming contradiction is further discussed and clarified in the end of Introduction.

Third, while our constructions are quite standard in the literature of analysis on fractals (see \cite{kigami,AnalOnFractals} for further references), neither singularity of Sobolev spaces or Theorem \ref{intro:NoLips} have been positively verified before in the related works. Let us emphasize that both are open for the Sierpi\'nski gasket, which is arguably one of the simplest and most studied examples in Analysis on fractals.
Although, negative cases exist by certain tree-like fractals in \cite[Section 2.5]{Baudoin_Chen_2023} and \cite[Remark 3.10]{scaleirregular}.
These two questions, nevertheless, seem to have attracted some attention.
The validity of singularity of Sobolev spaces in a general context, especially on the Sierpi\'nski carpet, was posed by third author and Murugan \cite[Problem 10.7]{murugan2023first}.
In \cite[Page 129, Section 4]{herman2004p} the authors performed a numerical study analyzing whether the usual inclusions $\mathscr{F}_{p_2} \subseteq \mathscr{F}_{p_1}$ for $p_1 < p_2$ hold on the Sierpi\'nski gasket. The conclusion of the authors seems to be that the data suggests the negative.
In \cite[Section 6.3]{murugan2024heat} Murugan conjectured that the canonical Sobolev space $\mathscr{F}_2$ on the Sierpi\'nski gasket contains no non-constant (Euclidean) Lipschitz function, namely for the positive of Theorem \ref{intro:NoLips}.
See also a related conjecture of Barlow and Perkins from 1988 \cite[Section 9]{barlow1988brownian}, which currently is still open.

\subsection{Proof of Theorem \ref{intro:Thm} (Outline)}
Let $\{G_n\}_{n \in \N}$ be as in Example \ref{intro: Construction example}.
Our key idea is to study the relationship between $\mathscr{E}_p$ and the \emph{optimal flows} $\mathcal{J}_{p,n}$ from the left most vertex $v_+$ to the right most $v_-$ on all $G_n$.
Flows on graphs are basically analogues to divergence free vector fields, and the optimal flow is the one with minimal $q$-mass for $q:=p/(p-1)$; see Subsection \ref{subsec: (preli) DPT} for precise definitions.
\begin{figure}[!ht]
    \begin{tikzpicture}[scale=0.55]
    \input{plots/LaaksoDflow}
    \end{tikzpicture} 
    \caption{Optimal flow $\mathcal{J}_{p,1}$.}
    \label{intro: flow}
\end{figure}

First, we prove that $\mathscr{E}_p$ can be expressed as a Dirichlet type $p$-energy,
\[
    \mathscr{E}_p(f) = \int_X \abs{D_{p}\Span{f}}^p \, d\Gamma_p\Span{\mathscr{U}_p}
\]
where $\Gamma_p\Span{\mathscr{U}_p}$ is a certain Radon measure and $D_p\Span{f} \in L^p(X,\Gamma_p\Span{\mathscr{U}_p})$ is some analogue of the usual gradient $\nabla f$.

The proof studies the measures $d\mu_{p}\Span{f} := \abs{D_p\Span{f}}d\Gamma_p\Span{\mathscr{U}_p}$. They have another nice and explicit description. For instance, the total mass can be computed
\begin{equation}\label{intro:TV}
    \mu_{p}\Span{f}(X) = \lim_{n \to \infty} \sum_{\{ v,w \} \in E_n} \abs{f(v) - f(w)}\abs{\mathcal{J}_{p,n}(v,w)}.
\end{equation}
So the value $\mu_{p}\Span{f}(X)$ can be understood as a weighted total variation where the weights are given by the optimal flows. By suitably localizing the weighted sums, $\mu_{p_i}\Span{f}(A)$ can be similarly computed for any Borel set $A \subseteq X$.

Now, the most crucial detail in our argument is the following property of the Laakso diamond. The weights that determine $\mu_{p}\Span{f}$ are \emph{independent} of $p$, meaning the optimal flows satisfy $\mathcal{J}_{p_1,n}(v,w) = \mathcal{J}_{p_2,n}(v,w)$ for all pairs $\{ v,w \} \in E_n$.
This can be seen from the ``horizontal symmetry'' of the Laakso diamond which forces the optimal flows to behave the same way. The case $n = 1$  can be seen from Figure \ref{intro: flow}.

Thus, if we had $f \in \mathscr{F}_{p_1} \cap \mathscr{F}_{p_2}$, then obviously $\mu_{p_1}\Span{f}(X) = \mu_{p_2}\Span{f}(X)$ due to \eqref{intro:TV}. By localizing the computation to an arbitrary Borel subset, we have the equality of Radon measures
\begin{equation}\label{intro:eqnuu}
	\mu_{p_1}\Span{f} = \mu_{p_2}\Span{f}. 
\end{equation}

On the other hand, the measures $\Gamma_p\Span{\mathscr{U}_p}$ on the Laakso diamond are given as Bernoulli type measures whose weight are different for every $p$. It is well-known fact in probability theory that two distinct Bernoulli measures, for instance on the interval $[0,1]$, are mutually singular. We use this to prove that
\begin{equation}\label{intro:SingularityEM}
	\text{$\Gamma_p\Span{\mathscr{U}_{p_1}}$ and $\Gamma_p\Span{\mathscr{U}_{p_2}}$ for distinct $p_1,p_2 \in (1,\infty)$ are mutually singular.}
\end{equation}

The combination of \eqref{intro:eqnuu} and \eqref{intro:SingularityEM} shows that the measures $\mu_{p_1}\Span{f}$ and $\mu_{p_2}\Span{f}$ are simultaneously equal and mutually singular. This is possible only when they are zero measures. Finally, it is a basic property that $\mu_{p}\Span{f}$ is the zero measure if and only if $f$ is constant.

\subsection{Inverse limits of Cheeger--Kleiner}\label{subsec:CK}
The literature actually already contains a class of natural Sobolev spaces on the Laakso diamond space, and also on several other examples this work studies, thanks to the inverse limit construction of Cheeger and Kleiner \cite{CK_PI}. Their work produces metric measure spaces, including Laakso diamond, supporting Poincar\'e inequality in the sense of Heinonen--Koskela \cite{HK}.
Moreover, it is well-known that Poincar\'e inequality and doubling property of the measure implies that Korevaar--Schoen Sobolev space is equivalent to the Newton--Sobolev; see \cite[Chapter 10]{HKST}.
Because Lipschitz functions are contained in the Newton--Sobolev spaces, Theorems \ref{intro:Thm}-\ref{intro:NoLips} in the Cheeger--Kleiner setting of Laakso diamond are obviously false.

Now, it might seem that our theorems contradict Cheeger--Kleiner but let us elaborate why this is not the case.
In our construction, the underlying measure is always the natural Ahlfors regular Hausdorff measure. In Cheeger--Kleiner the measure is some general doubling measure satisfying certain axioms. But in the case of Laakso diamond, the Hausdorff measure fails to satisfy some of these axioms. 
In conclusion, the frameworks differ at the choice of the measure.

In our general construction, there are other subtle differences from Cheeger--Kleiner in the projective structures, but we will not pursue this discussion any further. An interested reader may for instance compare Examples \ref{example: Counterexample} and \ref{ex: Singularities not equivalent} with the axioms of \cite{CK_PI}.
On the other hand, for some other examples where the Hausdorff measure satisfies the axioms of \cite{CK_PI}, our construction actually agrees with Cheeger--Kleiner; see Remark \ref{rem:CK} for further discussion.
We find it quite striking that such seemingly harmless difference causes decisively different forms of analysis.

\subsection{Organization of the paper}
This paper is organized to generalize the construction of the $p$-energies $\mathscr{E}_p$, that was described above on the Laakso diamond space, to a larger class of fractal spaces. We use the framework of (edge)-iterated graph systems (IGS) that was introduced by the first two authors \cite{anttila2024constructions}.

In Section \ref{Sec: Preliminary} we recall some classical definitions and results from metric spaces and graphs.
In Section \ref{sec:ssLaakso}, we introduce the general IGS construction.

Then in Section \ref{sec: DPT}, we study discrete potential theory of IGSs. We establish two key ingredients, the strong monotonicity principle and a Poincar\'e type inequality.

The construction of self-similar $p$-energy forms, as well as the investigation of their general properties, are covered by Sections \ref{Sec: Discretizations and mollifiers}-\ref{Sec: Sobolev spaces}.
In Section \ref{sec:otherSob}, we compare our framework to two other typical constructions of Sobolev spaces, Korevaar--Schoen and Newton--Sobolev.

In Section \ref{sec:SingularitySobolev}, we prove and discuss the main result of the paper, the Singularity of Sobolev spaces.

\section*{Acknowledgments}
The first author is supported by the Finnish Ministry of Education and Culture’s Pilot for Doctoral Programmes (Pilot project Mathematics of Sensing, Imaging and Modelling). Part of the research work was conducted during his participation in the Trimster program \emph{Metric Analysis}, organized by Hausdorff Research Institute of Mathematics (HIM) in Bonn, and was financially supported by the Deutsche Forschungsgemeinschaft (DFG, German Research Foundation) under Germany's Excellence Strategy - EXC-2047/1 - 390685813. He thanks the institute and local organizers for their hospitality and other participants for the stimulating atmosphere.
The second author is supported by the Research Council of Finland via the project \emph{GeoQuantAM: Geometric and Quantitative Analysis on Metric spaces}, grant no. 354241.
The third author (JSPS Research Fellow PD) is supported in part by JSPS KAKENHI Grant Number JP23KJ2011. The work was started during the visit of the third author to University of Jyväskylä. The authors thank Mathav Murugan for many inspiring discussions, and Naotaka Kajino for his comments on an earlier version of the paper.

\section{Preliminary}\label{Sec: Preliminary}
We begin by recalling some classical terminology and results from metric geometry and graph theory.

\subsection{Metric spaces and measures}
Let $(X,d,\mu)$ be a metric measure space where $(X,d)$ is a complete metric space and $\mu$ is a Radon measure on $(X,d)$.
We denote the \emph{open balls}
\[
    B(x,r):=\{y\in X : d(x,y)<r\} \text{ for } x \in X \text{ and } r \in (0,\infty).
\]
The \emph{diameter} of a non-empty subset $A \subseteq X$ is $\diam(A,d) := \sup_{x,y \in A} d(x,y)$. If $B \subseteq X$ is another non-empty subset, then the \emph{distance} between $A$ and $B$ is $\dist(A,B,d) := \inf_{x \in A,y \in B} d(x,y)$.
If $A \subseteq X$ is any Borel set with $\mu(A) > 0$ and $f : X \to \R$ is any Borel measurable function that is integrable on $A$, we denote the \emph{integral average} of $f$ over $A$ by
\[
    f_A := \kint_{A} f \,d\mu := \frac{1}{\mu(A)}\int_A f\,d\mu.
\]
When $A \subseteq X$ is a non-empty Borel subset we denote $\mu\restr_{K}$ the Radon measure on $(K,d)$ obtained by restricting $\mu$ to the Borel subsets of $K$.

\begin{definition}\label{def: AhlforsRegular}
A Radon measure $\mu$ on a metric space $(X,d)$ is \emph{doubling} if there is a constant $D > 0$ so that for all $x \in X$ and $r > 0$,
\[
    \mu(B(x,2r)) \leq D \cdot  \mu(B(x,r)).
\]
For $Q \in (0,\infty)$, we say that $(X,d,\mu)$ is \emph{$Q$-Ahlfors regular} if there is $C \geq 1$ so that for all $x \in X$ and $r \in (0,2\diam(X))$,
\[
    C^{-1} \cdot r^Q \leq \mu(B(x,r)) \leq C \cdot r^Q.
\]
\end{definition}

\subsection{Terminologies of graphs}
A \emph{graph} is a pair $G := (V,E)$ where $V$ is a non-empty finite set of \emph{vertices} and $E$ is a finite multiset of \emph{edges} consisting of unordered pairs $\{ v,w\}$ for distinct vertices $v,w \in V$.
We do not allow loops (edges of the form $\{v,v\}$) in our graphs. The \emph{degree} of a vertex $v \in V$ is the number of edges in $E$ containing $v$,
\[
    \deg(v) := \abs{\{ e \in E : v \in e \}}
\]
Note that the multiplicities are counted, i.e., if there are distinct edges with same endpoints, then these are counted separately.
In this paper, graphs are always assumed to be finite, meaning that $\abs{V},\abs{E} < \infty$.

\begin{definition}
    Let $G := (V,E)$ be a graph and $k \in \N$.
    A sequence of vertices $\theta :=  [v_1,\dots,v_k]$ is a \emph{path} in $G$ if $\{v_i,v_{i+1}\} \in E$ for all $1 \leq i < k$. The length of a $\theta$ is $\len(\theta) := k-1$.
We say that $G$ is \emph{connected} if for every pair of vertices $v,w \in V$ there is a path $[v = v_1,\dots, v_k = w]$ in $G$.
For a connected graph we define the \emph{(shortest) path metric} $d_{G}(v,w) := \min_{\theta} \len(\theta)$, where the minimum is taken over all paths connecting $v$ and $w$.
\end{definition}

Note that if $G$ is a connected graph then $(V,d_G)$ is a metric space.

\begin{definition}\label{def: independent}
    Let $G := (V,E)$ be a graph and $A \subseteq V$ be any subset. We say that $A$ is  \emph{connected} if every pair of points in $A$ can be connected by path entirely contained in $A$. We also say that $A$ is \emph{independent} if for every pair of vertices in $v,w \in A$ we have $\{ v,w \} \notin E$.
\end{definition}

\begin{definition}
    Let $G := (V,E)$ be a graph and $A \subseteq V$ be any subset. The \emph{boundary} of $A$ is the subset $\partial{A} := \{ v \in V \setminus A : \{ v,w \} \in E \text{ for some } w \in A \}$. The \emph{closure} of $A$ is $\overline{A} := A \cup \partial{A}$.
\end{definition}

\subsection{Discrete potential theory}\label{subsec: (preli) DPT}
Next we recall some relevant definitions and classical results in discrete potential theory. 
See e.g. \cite{HolopainenpHarmonic} for further background.
Let $G := (V,E)$ be a graph and $p \in (1,\infty)$. For any $U : V \to \R$ we define its \emph{gradient} as the function $\nabla U : V \times V \to \R$ given by $\nabla U(v,w) := U(w) - U(v)$. If $e = \{v,w\} \in E$ we denote $\abs{\nabla U(e)} := \abs{\nabla U(v,w)}$.
The \emph{graph $p$-energy} $\mathcal{E}_p : \R^V \to \R_{\geq 0}$ associated to $G$ is given by
\[
    \mathcal{E}_p(U) := \sum_{e \in E} \abs{\nabla U(e)}^p.
\]
If $A,B \subseteq V$ are non-empty disjoint subsets of $V$, we define the \emph{$p$-capacity} between $A$ and $B$ as
\[
    \cCap_p(A,B,G) := \inf\{ \mathcal{E}_p(U) : U|_A \equiv 1 \text{ and } U|_B \equiv 0 \}.
\]
Functions satisfying $U|_A \equiv 1$ and $U|_B \equiv 0$ are sometimes referred as \emph{potential functions} to the $p$-capacity problem $\cCap_p(A,B,G)$. Potential functions with minimal energy are referred as \emph{optimal potential functions}.
If $A \subseteq V$, we say that $U : V \to \R$ is \emph{$p$-harmonic} in $A$ if
\begin{equation*}\label{eq:p-harmonic}
        \sum_{\{ v,w \} \in E} \sgn( \nabla U(v,w) ) \abs{\nabla U(v,w)}^{p-1} = 0 \, \text{ for all } v \in A,
    \end{equation*}
where $\sgn(c)$ denotes the sign of $c \in \R$.
\begin{lemma}\label{lemma: Potential function exists}
    Let $G$ be a connected graph, $A,B \subseteq V$ be non-empty disjoint subsets and $p \in (1,\infty)$. There is a unique function $U : V \to \R$ so that $U|_A \equiv 1, U|_B \equiv 0$ and $\mathcal{E}_p(U) = \cCap_p(A,B,G) > 0$. Moreover, this function is $p$-harmonic in $V \setminus (A \cup B)$.
    The converse is also true. If $U$ is $p$-harmonic in $V \setminus (A \cup B)$ satisfying $U|_A \equiv 1, U|_B \equiv 0$ then $\mathcal{E}_p(U) = \cCap_p(A,B,G)$.
\end{lemma}

\begin{lemma}[Strong maximum principle]\label{lemma: Strong max principle}
    Let $G = (V,E)$ be a graph and $A \subsetneq V$ be a non-empty connected subset. If $U : V \to \R$ is $p$-harmonic in $A$ then
    \[
        \max_{x \in \overline{A}} U(x) = \max_{x \in \partial{A}} U(x) \text{ and } \min_{x \in \overline{A}} U(x) = \min_{x \in \partial{A}} U(x).
    \]
    Moreover, if $U(y) = \max_{x \in \overline{A}} U(x)$ for some $y \in A$ then $U$ is constant in $\overline{A}$.
\end{lemma}

Then we review the dual problem of capacity.
A function $\mathcal{J} : V \times V \to \R$ is said to be \emph{anti-symmetric} on $G$ if $\mathcal{J}(v,w) = -\mathcal{J}(w,v)$ for all $v,w \in V$, and $\mathcal{J}(v,w) = 0$ unless $\{v,w\} \in E$. For simplicity, if $e = \{v,w\} \in E$, we denote $\abs{\mathcal{J}(e)} := \abs{\mathcal{J}(v,w)}$.
Further, given non-empty disjoint subsets $A,B \subseteq V$ we say that an anti-symmetric function $\mathcal{J}$ is a \emph{flow} from $A$ to $B$ if its \emph{divergence} at a vertex $v \in V$
\[
    \divr(\mathcal{J})(v) := \sum_{ \substack{ \{v,w\} \in E } } \mathcal{J}(v,w)
\]
is equal to $0$ for all $v \in V \setminus (A \cup B)$. The flow $\mathcal{J}$ from $A$ to $B$ is a \emph{unit flow} if
\[
    \sum_{v \in A} \divr(\mathcal{J})(v) = 1.
\]
The \emph{$p$-energy} of a unit flow $\mathcal{J}$ is defined as
\[
    \mathcal{E}_q(\mathcal{J}) := \sum_{e \in E} \abs{\mathcal{J}(e)}^q,
\]
where $q := \frac{p}{p-1} \in (1,\infty)$ is the \emph{dual exponent} of $p$. Given non-empty disjoint subsets $A,B \subseteq V$, we define the \emph{$p$-resistance} between $A$ and $B$ as
\[
    \cRes_p(A,B,G) := \inf\{ \mathcal{E}_q(\mathcal{J}) : \mathcal{J} \text{ is a unit flow from } A \text{ to } B \}.
\]
The unit flow from $A$ to $B$ with minimal $p$-energy is referred as the \emph{optimal unit flow}.
Whenever the exponent $p$ is clear from the context, $q \in (1,\infty)$ always denotes the dual exponent of $p \in (1,\infty)$.

\begin{lemma}[Divergence theorem]
    Let $U : V \to \R$ be any function and $\mathcal{J}: V \times V \to \R$ be any antisymmetric function on a graph $G = (V,E)$. Then
    \begin{equation}\label{eq:divthm}
        \sum_{v \in V} \divr(\mathcal{J})(v) \cdot U(v)
        =  -\sum_{ \{v,w\} \in E } \mathcal{J}(v,w)\cdot \nabla U(v,w).
    \end{equation}
\end{lemma} 

Next we recall the duality of potentials and flows. See e.g. \cite[Theorem 5.1]{NakamuraYamasakiDuality}.

\begin{lemma}[Duality]\label{lemma:duality}
    Let $G$ be a connected graph, $A,B \subseteq V$ be non-empty disjoint subsets and $p \in (1,\infty)$. Then there is a unique unit flow $\mathcal{J}$ from $A$ to $B$ satisfying $\mathcal{E}_q(\mathcal{J}) = \cRes_p(A,B,G)$. Moreover, if $U$ is the unique solution to the $p$-capacity problem $\cCap_p(A,B,G)$ as in Lemma \ref{lemma: Potential function exists}, then
    \begin{equation}\label{eq:duality (potentials and flows)}
        \abs{\mathcal{J}(e)} = \mathcal{E}_p(U)^{-1} \cdot \abs{\nabla U(e)}^{p-1}.
    \end{equation}
    It also holds that
    \begin{equation}\label{eq:duality}
        \mathcal{E}_p(U)^{\frac{1}{p}} \cdot \mathcal{E}_q(\mathcal{J})^{\frac{1}{q}} = 1.
    \end{equation}
\end{lemma}

\section{Construction of Laakso-type fractals}\label{sec:ssLaakso}
This section covers the geometric framework of Laakso-type fractal spaces that was introduced in \cite{anttila2024constructions}.

\subsection{Iterated graph systems}
Our construction of Laakso-type fractal spaces is based on the framework \emph{(edge)-iterated graph systems} introduced in \cite{anttila2024constructions}.
The terminology arises from the iteration procedure that can be viewed as an analogue of iterated function system of graphs.
Vertex variant of the construction was studied in \cite{ReplacementGraphs2024}.
The usage of the term ``iterated graph systems'' was inspired by the work of Neroli Z. \cite{neroli2024fractal} where the author studied a similar construction with the primary focus on graph theory. See also \cite{neroli2024fractal,li2025reducibleiteratedgraphsystems,generator_term,XI2017ScaleFree,XI2019}.

\begin{definition}\label{def: IGS}
An \emph{iterated graph system (IGS)} consists of the data $(V_1,E_1,I)$ and of the collection of functions $\phi_{v,e} : I \to V_1$ for all pairs $e \in E_1$ and $v \in V_1$ so that $v \in e$.
We assume that they satisfy the following properties.
\begin{enumerate}
\vskip.3cm
    \item $G_1 := (V_1,E_1)$ is a connected finite graph and $I$ is a non-empty finite set called the \emph{gluing set}.
    \vskip.3cm
    \item  \label{def: IGS (2)}
    For each $e \in E_1$ and its end point $v\in e$ the function $\phi_{v,e}:I \to V_1$ is an injection whose image $I_{v,e} := \phi_{v,e}(I)$ is an independent set of $V_1$ (recall Definition \ref{def: independent}).
    \vskip.3cm
    \item For each $e = \{v,w\} \in E_1$ it holds that $I_{v,e}\cap I_{w,e}=\emptyset$. 
\end{enumerate}
The graph $G_1$ will be referred as the \emph{generator} of the IGS, and the set $I$ together with the maps $\phi_{v,e}$ will be referred as the \emph{gluing rules}. The use of the term ``generator'' was inspired by \cite{generator_term}.
\end{definition}

% \begin{remark}A variant of IGS-framework, the vertex-iterated graph systems, was introduced by the first and second author in \cite{ReplacementGraphs2024}. Current work does not regard vertex-IGSs in any kind of form.\end{remark}

\begin{remark}
    In later sections, the notations $V_1,E_1,G_1,I,\phi_{v,e}$ always refer to the data associated to the IGS whenever the IGS is clear from the context. For the most part, we will not explicitly restate this association.
\end{remark}

\subsection{Replacement rule}\label{subsec:RR}
We discuss an iterative procedure associated to a given IGS consisting of the data $V_1,E_1,I,\{ \phi_{v,e} \}_{v \in e \in E_1}$ and $G_1 := (V_1,E_1)$ that produces an infinite sequence of graphs $(G_n)_{n \in \N}$. The reader is strongly encouraged to compare the procedure below to the construction of Laakso diamond described in Example \ref{intro: Construction example}.

In general, a graph $G=(V,E)$ is said to be \emph{labeled} by $G_1$ if for every $v \in e \in E$ there is an associated injective mapping $\phi_{v,e}:I\to V_1$, whose image is an independent set.
Given a labeled graph $G$, we can form a new graph $\widehat{G}=(\widehat{V},\widehat{E})$, whose vertices are $\widehat{V}=V_1\times E / \sim$, where we identify 
\[
(\phi_{v,e}(a),e)\sim (\phi_{v,f}(a),f)
\]
for every $e,f\in E$ which share an end point $v \in V$ and $a\in I$. Further, we define the set of edges
\[
\widehat{E}=\{[(v,e)],[(w,e)]: \{v,w\}\in E_1 \text{ and } e\in E\}.
\]
This amounts to replacing each edge in $G$ by a copy of $G_1$, which are glued along the images of the mappings $\phi_{v,e}$. We can also define a labeling for $\widehat{G}$ by 
\[
\phi_{[(v,e)], \{[(v,e)],[(w,e)]\}} = \phi_{v,\{v,w\}}
\]
so that $\widehat{G}$ is labeled by $G_1$. This procedure is called a \emph{replacement rule}.

By applying the replacement rule recursively to the generator $G_1$, we can construct an infinite sequence of graphs: Set $G_{n+1}:=\widehat{G}_n$ for $n \in \N$.
The edges and vertices of $G_{n+1}=(V_{n+1},E_{n+1})$ can be described as follows.
\begin{enumerate}
\vskip.3cm
    \item Let $V_{n+1}=V_1\times E_n / \sim$, where we identify vertices with the relationships $(\phi_{v,e}(a),e)\sim (\phi_{v,f}(a),f)$ for every $e,f\in E_n$ which share an end point $v$ and $a\in I$.
    \vskip.3cm
    \item $\{[v,e],[w,e]\}\in E_{n+1}$ if $\{v,w\}\in E_1$. 
    \vskip.3cm
    \item $\phi_{[v,e], \{[v,e],[w,e]\}} = \phi_{v,\{v,w\}}$.
    \vskip.3cm
\end{enumerate}
Notice that in the notation of equivalence classes $[(v,e)]$, we drop the parenthesis and write $[v,e]$. We call the graphs $G_n$ \emph{replacement graphs}.

The replacement graphs have a natural \emph{projective structure}. For $n \in \N$ we define $\pi_{n+1}:V_{n+1}\cup E_{n+1} \to V_n \cup E_n$ as follows. For each vertex $[v,e]\in V_{n+1}$ define $\pi_{n+1}([v,e])=e$ if $v\not\in I_{w,e}$ for any $w\in e$, and otherwise set $\pi_{n+1}([v,e])=w$ if $v\in I_{w,e}$. For an edge $\{[v,e],[w,e]\}\in E_{n+1}$ we define $\pi_{n+1}(\{[v,e],[w,e]\}) = e$. For $n,m \in \N$ so that $n > m$, we define $\pi_{n,m} := \pi_{m+1} \circ \dots \circ \pi_{n-1} \circ \pi_{n}$. We also define $\pi_{n,n} \coloneqq \id_{V_{n} \cup E_{n}}$. If $n>m$, $w\in V_{n}$ and $\pi_{n,m}(w)=v\in V_m$, then we call $w$ an \emph{ancestor} of $v$.

\subsection{Some examples}\label{Subsec: First examples}
Let us introduce some relevant and simple examples.
In the figures we only include the first two graphs $G_1$ and $G_2$. 
It also may be helpful for the reader to revisit the construction of the Laakso diamond in Example \ref{intro: Construction example} where the first three graphs $G_1,G_2$ and $G_3$ are drawn.

\begin{example}[Laakso space]\label{example: Laakso space} Laakso introduced a construction of metric spaces \cite{La00} (today known as \emph{Laakso spaces}) that satisfy arguably one the strongest pair of analytic conditions in modern analysis on metric spaces; Ahlfors regularity and $1$-Poincar\'e inequality.
\begin{figure}[!ht] 
\centering\includegraphics[width=240pt]{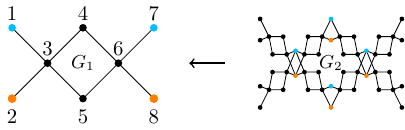}
\caption{Figure of the first replacement that produces the variant of a Laakso space discussed in Example \ref{example: Laakso space}. The colors of vertices indicate the gluing rules: \textcolor{mycolor}{Blue} vertex is connected to the other \textcolor{mycolor}{blue} vertex and \textcolor{orange}{orange} to the other \textcolor{orange}{orange}.} \label{fig: Laakso 1} \end{figure}
The following IGS produces a fractal space that can be regarded as a variant of a Laakso space. A further generalization is discussed in Proposition \ref{prop: dwp = p}, and see also Figure \ref{fig: genlaakso}.
The generator $G_1 = (V_1,E_1)$ is given in Figure \ref{fig: Laakso 1}. The gluing rules are given by $I := \{a,b\}$
\[
    (\phi_+(a),\phi_+(b)) := (1,2) \text{ and }  (\phi_-(a),\phi_-(b)) := (7,8).
\]
For $\{ v,w \} \in E_1$, we define $\phi_{v,\{v,w\}} = \phi_+$ if and only if $v < w$. When $v > w$, we define $\phi_{v,\{v,w\}} = \phi_-$.
\end{example}

\begin{example}\label{example: Counterexample}
A major motivation for the IGS-framework is that it can produce interesting metric spaces with weaker analytic properties compared to e.g. the Laakso space in Example \ref{example: Laakso space}.
\begin{figure}[!ht]
    \centering\includegraphics[width=240pt]{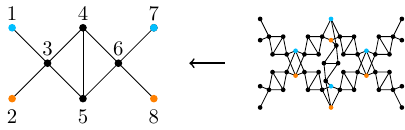}
    \caption{An IGS that produces a counterexample to Kleiner's conjecture.}
    \label{fig: Counterexample}
    \end{figure}
The construction in Figure \ref{fig: Counterexample} produces a metric space that negatively resolved Kleiner's conjecture \cite[Conjecture 7.5]{KleinerICM}; see \cite[Theorem 1.3]{anttila2024constructions}.
In Theorem \ref{thm:EMdichotomy} we see that the validity of Kleiner's conjecture is closely related to the structure of the certain Sobolev space; see Remark \ref{rem:Confdim} for a detailed explanation.
\end{example}

\begin{example}[Vicsek set]\label{example: Vicsek set}
    We discuss an IGS that produces a variant of the Vicsek set; see e.g. \cite{Baudoin_Chen_2023} for the standard construction.
    The generator $G_1$ is in the left of Figure \ref{fig: Vicsek set}.
    The gluing rules are defined as follows. First, we set $I := \{a\}, (\phi_-,(a),\phi_+(a)) := (5,1)$. For a given edge $\{v,w\} \in E_1$, we define
    $\phi_{v,\{v,w\}} = \phi_+$ if and only if $v < w$.
    Notice that the standard Vicsek set has 5 similarity maps and our version has 4.
    \begin{figure}[!ht] 
    \centering\includegraphics[width=310pt]{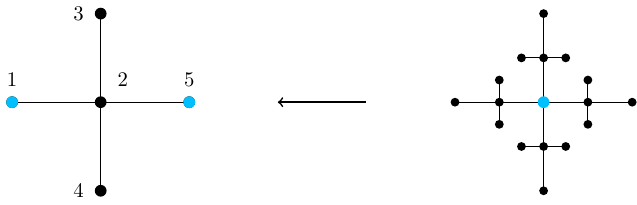}
    \caption{A figure of the IGS discussed in Example \ref{example: Vicsek set}.} \label{fig: Vicsek set} 
    \end{figure}
\end{example}

\begin{example}
    In Figure \ref{fig: multigraph} we provide an example of a multigraph IGS.
    \begin{figure}[!ht]
    \centering\includegraphics[width=310pt]{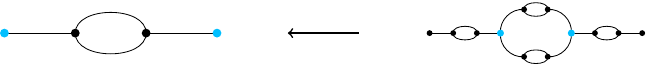}
    \caption{Figure of an IGS whose generator is a multigraph.}
    \label{fig: multigraph}
    \end{figure}
\end{example}

\begin{remark}\label{rmk:nonselfsim}
    The IGS framework exclusively produces self-similar metric spaces. By modifying the replacement rule to allow replacements with multiple different generators, we could easily construct metric spaces without exact self-similarity. We could also consider weighted graphs.
    Nevertheless, for such extensions, the treatment of the analogous notion of conductively uniform property (Definition \ref{def: Cond.unif}) would be somewhat delicate.
\end{remark}

\subsection{Geometric conditions}\label{subsec: Geometric assumptions}
In this subsection, we discuss some conditions on IGSs that ensure controlled and non-trivial geometry for the Laakso type fractals.

\begin{definition}\label{def: Simple}
    An IGS is \emph{simple} if there are two functions $\phi_+,\phi_- : I \to V_1$ so that for all $\{ v,w \} \in E_1$ we have 
    \[
        \left\{ \phi_{v,\{v,w\}},\phi_{w,\{v,w\}} \right\} = \{ \phi_+,\phi_- \}.
    \]
    When the IGS is simple, the gluing sets are denoted as
    \[
      I_{+} := \phi_+(I) \text{ and } I_{-} := \phi_-(I).
    \]
\end{definition}

For simple IGSs it is convenient to introduce the level-0 graph $G_0 := (V_0,E_0)$. We define 
\[
    V_0 := \{ v_+,v_- \} \text{ and } E_0 :=  \{e_0\} := \{ \{ v_+,v_- \} \}.
\]
Moreover, we give $G_0$ the gluing rules
\[
    \phi_{v_+,e_0} := \phi_+ \text{ and } \phi_{v_-,e_0} := \phi_-.
\]
Notice that the replacement rule on $G_0$ gives $G_1$. We then define the projection $\pi_1 : V_{1} \cup E_{1} \to V_0 \cup E_0$ so that
\[
    \pi_1(e) := e_{0} \text{ for all } e \in E_1 \text{ and } \pi_1(v) := 
    \begin{cases}
        v_+ & \text{ if } v \in I_+ \\
        v_- & \text{ if } v \in I_- \\
        e_{0} & \text{ otherwise. }
    \end{cases}
\]
In particular, it holds that
\[
    \pi_1^{-1}(v_+) = I_+ \text{ and } \pi_1^{-1}(v_-) = I_-.
\]
We also denote $\pi_{n,0} := \pi_{1} \circ  \pi_{n,1}$.
The set of all edges/vertices in all replacement graphs are denoted
\[
    E_{\#} := \bigsqcup_{n = 0}^\infty E_n \text{ and } V_{\#} := \bigsqcup_{n = 0}^\infty V_n.
\]
For all edges $e \in E_{\#}$ we shall denote the vertices $e^+,e^- \in e$ so that
\[
(\phi_{e^+,e},\phi_{e^-,e}) = (\phi_+,\phi_-).
\] 
Lastly, we define the higher-order gluing sets by
\[
    I_{+}^{(n)} := \pi_{n,0}^{-1}(v_+) \text{ and } I_{-}^{(n)} := \pi_{n,0}^{-1}(v_-).
\]
For $v \in e \in E_{\#}$ we define
\[
    I_{v,e}^{(n)} :=
    \begin{cases}
        I_{+}^{(n)} \text{ if } v = e^+ \\
        I_{-}^{(n)} \text{ if } v = e^-.
    \end{cases}
\]

\begin{definition}\label{def: doubling}
A simple IGS is said to be \emph{doubling} if for every $v \in I_{+} \cup I_-$ we have $\deg(v) = 1$.
For such $v \in I_{+} \cup I_-$ we denote the unique neighbor of $v$ by $\mathfrak{n}(v)$, and the corresponding edge by $\mathfrak{e}(v)$, i.e., $\mathfrak{e}(v) = \{ v, \mathfrak{n}(v)\}$.
\end{definition}

The doubling property ensures that the degrees do not grow uncontrollably. The following lemma is proven in \cite[Lemma 3.21]{anttila2024constructions}.

\begin{lemma}\label{lemma: Degree of ancestors}
    If the simple IGS is doubling then
    \[
        \sup_{v \in V_{\#}} \deg(v) = \max_{v \in V_1} \deg(v).
    \]
\end{lemma}

\begin{definition}
    A simple IGS is \emph{non-degenerate} if $\dist(I_+,I_{-},d_{G_1}) \geq 2$. Equivalently, there is no edge between the sets $I_+$ and $I_-$.
\end{definition}

Non-degeneracy removes uninteresting examples such as the case where $G_1$ contains only one edge. It also ensures a suitable growth rate of discrete distances. This is seen in the following lemma, which was essentially proven in \cite[Lemma 3.18]{anttila2024constructions} using the path decomposition in \cite[Proposition 3.16]{anttila2024constructions}. The proof is identical so we omit the details.

\begin{lemma}\label{lemma: Distance between ancestors}
    Let $L_* := \dist(I_+,I_-,d_{G_1} )$. If $n,m \in \N \cup \{0\}$ and $v,w \in V_m$ are distinct vertices, then
    \[
        \dist(\pi_{n+m,m}^{-1}(v),\pi_{n+m,m}^{-1}(w),d_{G_{n+m}}) \geq L_*^n.
    \]
\end{lemma}

\begin{lemma}\label{Lemma: Degree of boundary}
    If the IGS is simple and doubling then $\deg(v) = 1$ for all $n \in \N$ and $v \in I_{+}^{(n)} \cup I_{-}^{(n)}$. Let $\fn(v)$ denote the unique neighbor of such $v$. If the IGS is also non-degenerate, then
    \[
        \fn(v) \in V_n \setminus \left(I_{+}^{(n)} \cup I_{-}^{(n)}\right).
    \]
\end{lemma}

\begin{proof}
    The existence and the uniqueness of $\fn(v)$ is a direct inductive argument on $n \in \N$. The base case is the doubling property of the IGS and the inductive step follows from the gluing rules.
    
    Next, we argue the latter part, and assume $v \in I_-^{(n)}$ for simplicity.
    First, $\fn(v) \notin I_-^{(n)}$ follows from the fact that the gluing sets are independent sets (Definition \ref{def: IGS}-\eqref{def: IGS (2)}). According to Lemma \ref{lemma: Distance between ancestors}, $\fn(v) \in I_+^{(n)}$ would violate the fact that the IGS is non-degenerate.
\end{proof}

\subsection{Symbolic dynamics}\label{subsec: Symbolic dynamics}
The purpose of this subsection is to simplify the symbolic language of IGSs into a more geometric and intuitive form.

The following proposition gives a precise meaning to the self-similarity of our construction, and was essentially proven in \cite[Proposition 3.11]{anttila2024constructions}. There is a slight difference in \ref{SM2} because in this paper we allow multigraphs. The proof would be identical so we omit the details.

\begin{proposition}\label{prop: SM}
For every $n \in \mathbb{N} \cup \{ 0 \}$, $m \in \N $ and $e \in E_n$ there is a mapping $\sigma_{e,m} : V_m \to V_{n + m}$, the image of $\sigma_{e,m}$ denoted as $e \cdot G_m$ and the edges contained in this image as $e \cdot E_m$, with the following properties.
\begin{enumerate}[label={\textup{(\textcolor{blue}{SM\arabic*})}}, widest=a, leftmargin=*]
    \vskip.3cm
    \item \label{SM1} For every $e \in E_n$ the mapping $\sigma_{e,m}$ is injective and the collection of subsets $\{ e \cdot G_m \}_{e \in E_n}$ is a covering of $V_{n + m}$. For $v,w \in V_{m}$, $\{v,w\} \in E_m$ if and only if $\{ \sigma_{e,m}(v), \sigma_{e,m}(w) \} \in E_{n + m}$. Furthermore,  
    \begin{equation}\label{eq:phieq}
    \phi_{v,\{ v,w \}} = \phi_{\sigma_{e,m}(v),\{\sigma_{e,m}(v),\sigma_{e,m}(w)\}}.
    \end{equation}
    \vskip.3cm
    \item \label{SM2} For distinct edges $e,f \in E_n$ the subsets $e \cdot G_m$ and $f \cdot G_m$ intersect if and only if $e,f$ have a common vertex $v$.
    Moreover, their intersection is
    \[
        \bigcup_{v \in e \cap f} \sigma_{e,m}\left(I_{v,e}^{(m)}\right) = \pi_{n+m,n}^{-1}(e \cap f) = \bigcup_{v \in e \cap f} \sigma_{f,m}\left(I_{v,f}^{(m)}\right)
    \]
    \vskip.3cm
    \item \label{SM3} For every $e \in E_n$ we have $e \cdot E_m = \pi_{n+m,n}^{-1}(e)$.
    In particular, $\{ e \cdot E_m \}_{e \in E_n}$ is a partition of $E_{n + m}$.
\end{enumerate}
\end{proposition}

\begin{remark}\label{Remark: Describe similarity maps}
    The mappings $\sigma_{e,m} : V_m \to V_{n+m}$ are given inductively by the following conditions.
    \begin{enumerate}
        \item $\sigma_{e,0} : V_0 \to V_{n}$ is defined so that $\sigma_{e,0}(v_{\pm}) := e^{\pm}$.
        \item $\sigma_{e,m+1} : V_{m+1} \to V_{n + m+1}, \, [z,\{ v,w \}] \mapsto [z,\{\sigma_{e,m}(v),\sigma_{e,m}(w)  \}]$.
    \end{enumerate}
\end{remark}

\begin{corollary}\label{cor: Neigbours of ancestors}
     Suppose that the IGS is simple and doubling, and let $n,m \in \N$ and $v \in e \in E_n$. Then every $w \in \pi_{n+m,n}^{-1}(v)$ has a unique neighbor in $e \cdot G_{m}$.
     Let $\fn(w,e)$ denote the unique neighbor of such $w$ in $e \cdot G_{m}$. If the IGS is also non-degenerate then 
     \begin{equation}\label{eq: neighbour in interior}
         \fn(w,e) \in e \cdot G_m \setminus 
         (\pi_{n+m,n}^{-1}(e^+) \cup \pi_{n+m,n}^{-1}(e^-) )
     \end{equation}
\end{corollary}

\begin{proof}
    The claim follows from \ref{SM1} and Lemma \ref{Lemma: Degree of boundary}. Indeed, if $w = \sigma_{e,m}(z)$, then $\fn(w,e) = \sigma_{e,m}(\fn(z))$.
\end{proof}

When the IGS is simple and doubling, $v \in e \in E_n$ and $w \in \pi_{n+m,n}^{-1}(v)$, we define $\fe(w,e) := \{w, \fn(w,e)\}$.

\begin{definition}\label{def: n-Gluing maps}
For any $v \in V_n$ and $m \in \N$, we define the \emph{higher order gluing map} $\Phi_{v,m}:I^m \to \pi_{n+m,n}^{-1}(v)$ recursively as follows.
\begin{enumerate}
        \item $\Phi_{v,1}(a) := [\phi_{v,e}(a),e]$ for all $a \in I$.
        \item If $w := \Phi_{v,m}(a_1\dots a_{m}) \in \pi_{n+m,n}^{-1}(v)$ is given for all $a_1\dots a_m \in I^m$, we define
        \[
          \Phi_{v,m+1}(a_1\dots a_m a_{m+1}) := [\phi_{w,\fe(w,e)}(a_{m+1}),\fe(w,e)].  
        \]
\end{enumerate}
If $v = v_{\pm}$ we denote $\phi_{\pm,m} := \Phi_{v_{\pm},m}$.
\end{definition}

Note that the mappings $\Phi_{v,\bullet}$ do not depend on the choice of the edge $e$ containing $v$ due to the identifications at the end of a replacement. 

\begin{lemma}\label{lemma: Phi are bijective}
    For all $v \in V_n$ and $m \in \N$ the higher order gluing map $\Phi_{v,m}:I^m \to \pi_{n+m,n}^{-1}(v)$ is bijective.
\end{lemma}

\begin{proof}
    This follows from the injectivity of the first order gluing maps $\phi_{v,e}$ (Definition \ref{def: IGS}-\eqref{def: IGS (2)}).
\end{proof}

\begin{definition}\label{def: Symbolic product}
    We define the following symbolic operations.
    \begin{enumerate}
        \item If $e \in E_{\#}$ and $v \in V_n$ we define $e \cdot v := \sigma_{e,n}(v)$.
        \item If $e \in E_{\#}$ and $f \in E_{\#}$ we define $e \cdot f := \{ e \cdot f^-, e \cdot f^+ \}$.
        \item If $v \in V_{\#}$ and $a \in I^m$ we define $v \cdot a := \Phi_{v,m}(a)$.
        \item If $a = a_1\dots a_n \in I^n$ and $b_1\dots b_m \in I^m$ we define
    \[
        a \cdot b := a_1\dots a_n b_1\dots b_m \in I^{n+m}.
    \]
    \end{enumerate}
\end{definition}

\begin{lemma}\label{lemma: Symbolic associativity}
    Suppose that the IGS is simple and doubling.
    Let $e,f,g \in E_{\#}, v \in V_{\#}$ and $a \in I^n, b \in I^m$. Then the following associativity properties hold.
    \begin{enumerate}
        \item $(e \cdot f) \cdot v = e \cdot (f \cdot v)$ \label{item: Symbolic associativity (1)}
        \item $(e \cdot f) \cdot g = e \cdot (f \cdot g)$ \label{item: Symbolic associativity (2)}
        \item $(e \cdot v) \cdot a = e \cdot (v \cdot a)$ \label{item: Symbolic associativity (3)}
        \item $(v \cdot a) \cdot b = v \cdot (a \cdot b)$ \label{item: Symbolic associativity (4)}
    \end{enumerate}
\end{lemma}

\begin{proof}
    \eqref{item: Symbolic associativity (1)}-\eqref{item: Symbolic associativity (2)}: Follows from the explicit form of the mappings $\sigma_{\bullet,\bullet}$ in Remark \ref{Remark: Describe similarity maps}.

    \eqref{item: Symbolic associativity (3)}:
    Choose any edge $v \in e_v$, and first assume that $a \in I$. It follows from \eqref{eq:phieq} and the explicit form of of the mappings $\sigma_{\bullet,\bullet}$ that
    \begin{align*}
        (e \cdot v) \cdot a =
        \Phi_{e \cdot v}(a) = [\phi_{e \cdot v,e \cdot e_v}(a), e \cdot e_v] = e \cdot [\phi_{v,e_v}(a),e_v] = e \cdot (v \cdot a).
    \end{align*}
    This proves the case $n = 1$, and the general case $a \in I^n$ for $n \geq 1$ follows from a similar computation and induction.

    \eqref{item: Symbolic associativity (4)}:
    Suppose $e_v$ is again an edge containing $v$, and first assume $b \in I$. Then 
    \begin{align*}
        (v \cdot a) \cdot b = \Phi_{v \cdot a,1}(b) = [\phi_{v \cdot a, \fe(v \cdot a,e_v)}(b),\fe(v \cdot a,e_v)]
        = v \cdot (a \cdot b)
    \end{align*}
    and this proves the case $m = 1$. The general case would again follow from a direct inductive computation.
\end{proof}

\begin{remark}\label{remark: Identify as words}
    It is easy to see from \ref{SM3} that the set of edges $E_n$ can be identified as the set of sequences of length $n$,
    \[
        E_1^n := \{ e_1e_2\dots e_n : e_i \in E_1 \text{ for all } i = 1,2,\dots,n  \}.
    \]
    Using the product operation defined in Definition \ref{def: Symbolic product} and the associativity property proven in Lemma \ref{lemma: Symbolic associativity}-\eqref{item: Symbolic associativity (2)}, this identification can be naturally given by
    \[
        e_1e_2\dots e_n \mapsto e_1 \cdot e_2 \cdot \ldots \cdot e_n \in E_n.
    \]
    Motivated by this, we sometimes denote $e_1e_2\dots e_n := e_1 \cdot e_2 \cdot \ldots \cdot e_n \in E_n$. This notation particularly useful for defining measures. However, in order to avoid confusion, we always explicitly mention when this notation is used.
\end{remark}

\subsection{Limit space}\label{subsec: Limit space}
The finish line in the geometric framework of IGSs is the construction of the limit space.
The precise construction is quite different from the one described in Example \ref{intro: Construction example} for the Laakso diamond, but this is because the generality we aim to achieve seems to require an altered approach.

For this subsection, we fix a simple and non-degenerate IGS satisfying the doubling property.
Let us first define the \emph{symbol space} as the family of projective sequences
    \[
      \Sigma := \{ (e_i)_{i = 0}^{\infty} : e_i \in E_i \text{ and } \pi_{i+1}(e_{i+1}) = e_{i} \text{ for all } i \in \N \cup \{0\} \}.  
    \]
Here $\pi_i$ are as in Subsection \ref{subsec:RR}.
For $n \in \N \cup\{0\}$ and $e \in E_n$ we define the subsets $\Sigma_e \subseteq \Sigma$ by
    \[
        \Sigma_e := \left\{ (e_i)_{i = 0}^{\infty} \in \Sigma : e_n = e \right\},
    \]
    which are also obtained as the images of the mappings
    \[
        \sigma_e : \Sigma \to \Sigma, \, (e_i)_{i = 0}^{\infty} \mapsto (f_i)_{i = 0}^{\infty},
    \]
    where $f_{n + i} := e \cdot e_i$ for all $i \in \N \cup \{0\}$ and $f_{n-i} := \pi_{n,n-i}(e)$ for all $0 \leq i \leq n$.
    Note that $\sigma_{e_0} = \id_{\Sigma}$.
    We equip $\Sigma$ with the natural \emph{word metric}
    \[
        \delta_{\Sigma}((e_i)_{i = 0}^\infty,(f_i)_{i = 0}^\infty) :=
        \begin{cases}
            0 & \text{ if } e_i = f_i \text{ for all } i \in \N \cup \{0\} \\
            2^{-\min \left\{ k \ : \  (e_i)_{i = 0}^k \neq (f_i)_{i = 0}^k \right\}} & \text{ otherwise.}
        \end{cases}
    \]
    \begin{remark}\label{rem: Alternative symbol space}
        The symbol space alternatively can be defined as the space of sequences
    \[
        E_1^\N := \{ (e_i)_{i = 1}^\infty : e_i \in E_1 \text{ for all } i \in \N \}
    \]
    through the identification $T : E_1^{\N} \to \Sigma$,
    \[
        (e_i)_{i = 1}^{\infty} \mapsto (e_1e_2 \ldots e_i)_{i = 0}^{\infty} \in \Sigma.
    \]
    Here, $e_1e_2 \ldots e_0$ is understood as $e_0$.
    It follows from the discussion in Remark \ref{remark: Identify as words} that $T$ is a well-defined bijection.
    \end{remark}
    \begin{definition}\label{def: Measures on sigma}
    Let $\alpha : E_{\#} \to [0,\infty)$ be any function satisfying
    \begin{equation}\label{eq: Measures on sigma}
        \alpha(e) = \sum_{f \in E_1} \alpha(e \cdot f) \quad \text{ for all } e \in E_{\#}.
    \end{equation}
    Then we define $\mathfrak{m}_\alpha$ to be the Radon measure on $\Sigma$ satisfying
    \begin{equation}\label{def: mu_alpha}
         \mathfrak{m}_{\alpha}(\Sigma_e) = \alpha(e) \text{ for all } e \in E_{\#}.
    \end{equation}
    To see that $\mathfrak{m}_\alpha$ exists, we can construct it by first defining it on the alternative symbol space $E_1^{\N}$. Kolmogorov's extension theorem and the conditions \eqref{eq: Measures on sigma} ensure that there is a unique Radon measure $\widehat{\mathfrak{m}}_\alpha$ on $E_1^{\N}$ satisfying
    \[
        \widehat{\mathfrak{m}}_\alpha(\{ (e_i)_{i = 1}^\infty \in E_1^{\N} : e_1e_2\ldots e_n = e \}) = \alpha(e) \text{ for all } e \in E_n.
    \]
    By using the identification $T : E_1^{\N} \to \Sigma$ discussed in Remark \ref{rem: Alternative symbol space}, the pushforward measure $\mathfrak{m}_\alpha := T_*(\widehat{\mathfrak{m}}_\alpha)$ satisfies \eqref{def: mu_alpha}.
    If it holds that $\alpha(e_0) = 1$ and $\alpha(e_1e_2\dots e_n) = \alpha(e_1)\alpha(e_2)\dots\alpha(e_n)$ for all $n \in \N$ and $e_1,\dots,e_n \in E_1$, we say that the measure $\mathfrak{m}_\alpha$ is a \emph{Bernoulli measure} of $\Sigma$. Since $\mathfrak{m}_\alpha(\Sigma) = \alpha(e_0) = 1$, $\mathfrak{m}_\alpha$ is a probability measure.
    The values $\{ \alpha(e) \}_{e \in E_1}$ are called the \emph{weights} of $\mathfrak{m}_\alpha$.
    Observe that if we identify $\Sigma$ with $E_1^{\N}$, then $\mathfrak{m}_\alpha$ would be a Bernoulli measure in the usual sense.
    Lastly, we define the \emph{uniform Bernoulli measure} to be the Bernoulli measure $\mathfrak{m}_{\text{unif}}$ with weights $\alpha(e) = \abs{E_1}^{-1}$ for all $e \in E_1$. 
    \end{definition}

We are now prepared to introduce the precise definition of the limit space. Certain aspects, in particular the definition of the metric, are quite technical. Propositions \ref{prop: LS well-defined} and \ref{prop: Geom of LS} are collections of the geometric properties of the limit space. 
%Therefore in the later sections, we will reference the geometric properties of the limit space through Propositions \ref{prop: LS well-defined} and \ref{prop: Geom of LS}.

\begin{definition}\label{def: Limit space}
    The \emph{limit space} of an IGS consists of the metric measure space $(X,d,\mu)$ and mappings $\chi : \Sigma \to X$ and $F_e : X \to X$ for all $e \in E_{\#}$, which are constructed as follows. The underlying set $X$ is obtained as the quotient set
\begin{equation}\label{eq: Identifications in limit space}
    X := \Sigma/\sim, \text{ where } (e_i)_{i = 0}^{\infty} \sim (f_i)_{i = 0}^{\infty} \iff e_i \cap f_i \neq \emptyset \text{ for all } i \in \N \cup \{0\}.
\end{equation}
The \emph{coding map} is the canonical projection $\chi : \Sigma \to X$, and we define the measure $\mu$ on $X$ as the push-forward of the uniform Bernoulli measure $\mu := \chi_*(\mathfrak{m}_{\text{unif}})$.
We define the \emph{similarity maps} $F_e: X \to X$ for all $e \in E_{\#}$ as the unique mappings satisfying
\begin{equation}\label{def: similarity maps}
    F_e \circ \chi = \chi \circ \sigma_e. 
\end{equation}
The image of $F_e$ is denoted $X_e := F_e(X)$. For $v \in V_{\#}$, we denote the associated \emph{closed star} by
\[
    X_v := \bigcup_{\substack{e \in E_{\#}\\ v \in e}} X_e.
\]
The metric $d$ on $X$ is chosen to be so that it satisfies the \emph{visual metric property}. This means that there exists a constant $A \geq 1$ and $L_* > 1$ so that for all distinct sequences $(e_i)_{i = 1}^\infty, (f_i)_{i = 1}^\infty \in \Sigma$, if $n \in \N \cup \{0\}$ is the largest non-negative integer so that $e_n \cap f_n \neq \emptyset$, it holds that
\begin{equation}\label{eq: Visual metric}
    A^{-1} L_*^{-n}
        \leq d(\chi((e_i)_{i = 0}^\infty),\chi((f_i)_{i = 0}^\infty)) \leq A L_*^{-n}.
\end{equation}
The metric is explicitly given as follows.
Fix $L_* := \dist(I_+,I_-,d_{G_1})$, which is at least $2$ since the IGS is assumed to be non-degenerate. For $e \in E_\#$ denote $\abs{e} := n$ to be the integer so that $e \in E_n$.
We define a metric $d_n$ on $V_n$ by
\[
    d_{n}(v,w):=\inf_{\Omega}\sum_{e\in\Omega} L_*^{-\abs{e}} 
\]
    where the infimum is taken over all sets $\Omega \subseteq \bigcup_{j = 0}^n E_j$, so that
    \begin{equation}\label{def: Omega}
        \bigcup_{e \in \Omega} e \cdot G_{n - \abs{e}} \subseteq V_n \text{ is connected in $G_n$ and } v,w \in \bigcup_{e \in \Omega}e \cdot G_{n - \abs{e}}.
    \end{equation}
    The set $e \cdot G_m$ is as in Proposition \ref{prop: SM}.
    We now define
    \begin{equation}\label{def: metric}
        d(\chi\left((e_i)_{i = 0}^\infty\right),\chi\left((f_i)_{i = 0}^\infty)\right) := \lim_{n \to \infty} d_n(e_n,f_n).
    \end{equation}
Here, $d_n(e_n,f_n)$ is understood as
\[
    d_n(e_n,f_n) := \min_{\substack{v \in e_n \\ w \in f_n}} d_n(v,w). 
\]
In Proposition \ref{prop: LS well-defined} we show that the limit in \eqref{def: metric} always exists.
\end{definition}

\begin{remark}
    The definition of the metric is slightly modified from the one in \cite{anttila2024constructions}. This is because, in the current work we do not assume, the \emph{$L_*$-uniform scaling property} \cite[Definition 3.17]{anttila2024constructions}, and therefore the definition had to be adjusted.
\end{remark}

\begin{proposition}\label{prop: LS well-defined}
    Let $(X,d,\mu)$ be the limit space of a simple and non-degenerate IGS satisfying the doubling property, and denote its generator $G_1 := (V_1,E_1)$ and $L_* := \dist(I_+,I_-,d_{G_1})$.
    \begin{enumerate}
        \vskip.3cm
        \item \label{item: Visual metric}
        The metric $d$ on $X$ given in \eqref{def: metric} is well-defined and satisfies the visual metric property \eqref{eq: Visual metric}. Moreover, $\chi : \Sigma \to X$ is a continuous surjection.
        \vskip.3cm
        \item \label{item: Similarity maps}
        For all $n \in \N \cup \{0\}$ and $e \in E_n$ the condition \eqref{def: similarity maps} uniquely determines a well-defined mapping $F_e \colon X \to X$. Furthermore, these mappings are injective, $L_*^{-n}$-Lipschitz and their images $F_e(X) = X_e$ satisfy
        \begin{equation}\label{eq: SS of limit space}
            X = \bigcup_{e \in E_n} X_e.
        \end{equation}
    \end{enumerate}
\end{proposition}

\begin{proof}
    \eqref{item: Visual metric}: The arguments in this proof are modified versions of the ones in \cite[Subsection 3.6]{anttila2024constructions}.
    We first show that \eqref{def: metric} produces a well-defined metric. First notice that $d_n \leq 1$ for all $n$ by choosing $\Omega := \{e_0\}$. 
    Now, our first objective is obtain the inequalities
    \begin{equation}\label{eq: Visual metric estimate}
        d_{n}(e_{n},f_{n}) \leq d_{n+m}(e_{n+m},f_{n+m})  \leq 2L_*^{-n} + d_{n}(e_{n},f_{n})
    \end{equation}
    for all pairs of projective sequences $(e_n)_{n = 0}^\infty,(f_n)_{n = 0}^\infty \in \Sigma$.
    These are analogical estimates to \cite[(\textcolor{blue}{DL4}) and (\textcolor{blue}{DL5}) in Corollary 3.19]{anttila2024constructions}.

    We begin with the first one.
    If $e_n \cap f_n \neq \emptyset$ we clearly have $d_n(e_n,f_n) = 0 \leq d_{n+1}(e_{n+1},f_{n+1})$. Assume $e_n \cap f_n = \emptyset$ and let $\Omega$ satisfy \eqref{def: Omega} so that it minimizes the distance $d_{n+1}(e_{n+1},f_{n+1})$. By \ref{SM1}-\ref{SM2}, $\Omega$ also connects a pair of ancestors $v_{n+1},w_{n+1} \in V_{n+1}$ of $v_n \in e_n, w_n\in f_n$ respectively.
    This shows $d_{n+1}(v_{n+1},w_{n+1}) \leq d_{n+1}(e_{n+1},f_{n+1})$.
    By using Lemma \ref{lemma: Distance between ancestors} we see that $d_{n}(v_{n},w_{n}) = d_{n+1}(v_{n+1},w_{n+1})$. We thus have
    \[
    d_{n}(e_{n},f_{n}) \leq d_{n}(v_{n},w_{n}) = d_{n+1}(v_{n+1},w_{n+1}) \leq d_{n+1}(e_n,f_n).
    \]
    
    The latter inequality in \eqref{eq: Visual metric estimate} is obtained by setting $\widehat{\Omega} := \Omega \cup \{ e_n,f_n \}$ where $\Omega$ is a set satisfying \eqref{def: Omega} that minimizes the distance of $d_{n}(e_n,f_n)$. To see that the triangle inequality holds, suppose $e_n,f_n,g_n \in E_n$ and let $\Omega_1,\Omega_2$ satisfy \eqref{def: Omega} and minimize the distances $d_n(e_n,f_n),d_n(f_n,g_n)$ respectively. By setting $\Omega := \Omega_1 \cup \Omega_2 \cup \{f_n\}$, we get
    \[
        d_n(e_n,g_n) \leq d_n(e_n,f_n) + d_n(f_n,g_n) + L_*^{-n}.
    \]
    The triangle inequality now follows by letting $n \to \infty$. We now conclude that $d$ is a well-defined metric on $X$.

    We move on to verifying the visual metric property \eqref{eq: Visual metric}. Fix an arbitrary pair $(e_n)_{n = 0}^\infty,(f_n)_{n = 0}^\infty$ of projective sequences and let $n \in \N$ be the largest $n \in \N \cup \{0\}$ so that $e_n \cap f_n \neq \emptyset$.
    Then $d_{n+m}(e_{n+m},f_{n+m}) \leq 2L_*^{-n}$ (choose $\Omega := \{e_n,f_n\}$). On the other hand, since $e_{n+1} \cap f_{n+1} = \emptyset$, it follows from Lemma \ref{lemma: Distance between ancestors} that $d_{n+1}(e_{n+1},f_{n+1}) \geq L_*^{-(n+1)}$. The visual metric property now follows from \eqref{eq: Visual metric estimate}.
    Consequently, the continuity of $\chi : \Sigma \to X$, where $\Sigma$ is given the word metric, is clear from the visual metric property and $L_* > 1$.

    \eqref{item: Similarity maps}: Fix $e \in E_n$ for $n \in \N \cup \{0\}$.
    The uniqueness of the mapping $F_e$ satisfying \eqref{def: similarity maps} follows from the fact that $\chi$ is a surjection. To see that $F_e$ is well-defined, this follows from the definition of the identifications in \eqref{eq: Identifications in limit space} and \ref{SM1}.
    The condition \eqref{eq: SS of limit space} is then clear from $X_e = \chi(\Sigma_e)$.
    Lastly, to conclude that $F_e$ is $L_*^{-n}$-Lipschitz, this follows from the definition of the metric and \ref{SM1}.
\end{proof}

\begin{proposition}\label{prop: Geom of LS} 
    Let $(X,d,\mu)$ be the limit space of a simple and non-degenerate IGS satisfying the doubling property, and denote its generator $G_1 := (V_1,E_1)$ and $L_* := \dist(I_+,I_-,d_{G_1})$. Then the following hold.
    \begin{enumerate}
    \vskip.3cm
    \item \label{item: diams/measures of cells}
    For every $n \in \N \cup \{0\}$, $e \in E_n$ and $v \in V_n$ we have
        \[
            \diam(X_e) = L_*^{-n} \text{ and } \mu(X_e) = \abs{E_1}^{-n},
        \]
        and
        \[
            \diam(X_v) \leq 2L_*^{-n} \text{ and } \mu(X_v) =  \deg(v)\abs{E_1}^{-n}.
        \]
        Moreover, it holds that
        \begin{equation}\label{eq: Fibers = mu-null}
            \mu(X_e \cap X_f) = 0 \text{ for all distinct edges } e,f \in E_n \text{ and } n \in \N.
        \end{equation}
        \vskip.3cm
        \item \label{item: Stars = balls}
        Let $x \in X$ and $r \in (0,\infty) $ so that $\frac{1}{4}L_*^{-n-1} \leq r < \frac{1}{4}L_*^{-n}$ holds for some  $n \in \N$. Then there is a vertex $v \in V_n$ so that the closed star $X_v$ satisfies
        \[
            B(x,r) \subseteq X_v \subseteq B(x,8L_*r).
        \]
        \vskip.3cm
        \item \label{item: Q-AR}
        $(X,d,\mu)$ is compact, path connected and $Q$-Ahlfors regular for
        \[
            Q := \frac{\log(\abs{E_1})}{\log(L_*)}.
        \]
    \end{enumerate}
\end{proposition}

The proof of Proposition \ref{prop: Geom of LS} is postponed to the end of the next subsection.

\subsection{Fiber sets}\label{subsec: Fiber sets}
According to \ref{SM2}, the self-similar pieces of the replacement graphs intersect exactly at the set of ancestors. The analogous concept for the limit space are the fiber sets.
The precise definition of the fiber set is, again, somewhat technical but the main purpose of this subsection is to simplify the language involving it.

Recall from Lemma \ref{lemma: Phi are bijective} that the mapping $\Phi_{v,m}(a) = v \cdot a$ determines a bijection from $I^{m}$ onto $\pi_{n+m,n}^{-1}(v)$. We shall denote the space of sequences
\[
    I^\N := \{ (a_i)_{i = 1}^\infty : a_i \in I \text{ for all } i \in \N \},
\]
and equip it with the natural word metric. 

For $n \in \N \cup \{0\}$ and any $v \in e \in E_n$, we define the subset $\Sigma_{v,e} \subseteq \Sigma_e$ to be the set obtained as the image of the mapping $T_{v,e} : I^\N \to \Sigma, (a_i)_{i = 1}^\infty \mapsto (e_i)_{i = 0}^{\infty}$, where
    \begin{equation}\label{eq: Prefiber}
        e_i =
        \begin{cases}
            \fe(v \cdot (a_1 \dots a_{i-n}) ,e) & \text{ for all } i > n\\
            \pi_{n,i}(e) & \text{ if } i \leq n.
        \end{cases}
    \end{equation}
Note that by the uniqueness part of Corollary \ref{cor: Neigbours of ancestors}, we have
\[
  \fe(v \cdot a_1\dots a_m a_{m+1},e) = \fe(v \cdot a_1\dots a_m a_{m+1},\fe(v \cdot a_1\dots a_m,e)). 
\]
Combining this with Lemma \ref{lemma: Symbolic associativity}-\eqref{item: Symbolic associativity (4)} and the bijectiveness of the mappings $\Phi_{v,\bullet}$, it follows that $T_{v,e} : I^\N \to \Sigma$ is a homeomorphism onto $\Sigma_{v,e}$.
Then we define the \emph{fiber set} of $v \in V_n$ as the subset $\Fib(v) := \chi(\Sigma_{v,e}) \subseteq X$. We also denote
\[
    \Fib(X) := \bigcup_{v \in V_\#} \Fib(v).
\]

\begin{lemma}\label{lemma: Coding map well-defined}
    For all $v \in V_\#$ the fiber set $\Fib(v)$ does not depend on the choice of the edge $e \in E_\#$ containing $v$. Furthermore, $\chi|_{\Sigma_{v,e}}$ is a homeomorphism onto $\Fib(v)$.
\end{lemma}

\begin{proof}
    Let $e,f \in E_n$ be two edges containing $v$. By the definition of the projections, we have $\pi_{n,i}(e) \cap \pi_{n,i}(f) \neq \emptyset$ for all $0 \leq i \leq n$. 
    For $(a_i)_{i =  1}^\infty \in I^{\mathbb{N}}$, let $v_{n+i}=v\cdot (a_1\dots a_i)$ for $i>0$. Then 
    \[
        v_{n+i} \in \fe(v \cdot a_1\dots a_i,e) \cap \fe(v \cdot a_1\dots a_i,f).
    \]
    Thus, by the identifications in \eqref{eq: Identifications in limit space}, we have $\chi(\Sigma_{v,e}) = \chi(\Sigma_{v,f})$, and we conclude that 
    $\Fib(v)$ does not depend on the choice of the edge containing $v$.

    We next move on to proving that $\chi|_{\Sigma_{v,e}}$ is a homeomorphism onto $\Fib(v)$. Since $\Sigma_{v,e}$ can be regarded as a continuous image of the compact metric space $I^{\N}$ with the word metric, it is clear that $\Sigma_{v,e} \subseteq \Sigma$ is a compact subset. Thus, it is sufficient to verify that $\chi|_{\Sigma_{v,e}}$ is injective.
    Suppose $(a_i)_{i = 1}^{\infty}, (b_i)_{i = 1}^{\infty} \in I^{\N}$ are distinct, and let $j \geq 1$ be the smallest integer for which $a_j \neq b_j$. Since the mapping $I^n \ni a \mapsto v \cdot a = \Phi_{v,m}(a) \in \pi_{n+m,n}^{-1}(v)$ is bijective, we have that $v \cdot (a_1\dots a_j) \neq v \cdot (b_1 \dots b_j)$.
    It now follows from Lemma \ref{lemma: Distance between ancestors} that
    \[
        d_{G_{n + j+2}}(v \cdot (a_1 \dots a_{j+2}),v \cdot (b_1\dots b_{j+2})) \geq L_*^2 \geq 4.
    \]
    In particular, the edges
    $\fe(v \cdot (a_1 \dots a_{j+2}),e)$ and $\fe(v \cdot (b_1 \dots b_{j+2}),e)$ do not share a vertex.
\end{proof}

\begin{proposition}\label{prop: X_e cap X_f}
    Let $n \in \N$ and $e,f \in E_n$ be distinct edges. Then the sets $X_e$ and $X_f$ intersect if and only if $e$ and $f$ have a common vertex. Moreover, their intersection is
    \begin{equation}\label{eq: X_e cap X_f}
        X_e \cap X_f = \bigcup_{v \in e \cap f} \Fib(v).
    \end{equation}
\end{proposition}

\begin{proof}
    The first part of the claim, regarding when the sets $X_e$ and $X_f$ intersect, follows directly from the definition of the identification \eqref{eq: Identifications in limit space}.

    We will now prove \eqref{eq: X_e cap X_f}.
    Since $\Sigma_{v,e} \subseteq \Sigma_e$, we obviously have $\Fib(v) \subseteq X_e$ for all $v \in e \in E_{\#}$. Then assume that $x \in X_e \cap X_f$ for distinct edges $e,f \in E_n$, and write $x = \chi((e_i)_{i = 1}^{\infty}) = \chi((f_i)_{i = 1}^{\infty})$ for $e_n = e$ and $f_n = f$.
    Note that $e_i \cap f_i$ contains exactly one vertex for all $i > n$. To see this, suppose $v_i \in e_i \cap f_i$. By \ref{SM2}, $v_i$ is an ancestor of $v \in e \cap f$. Then we necessarily have
    \[
        e_i = \fe(v_i,e) = \{ v_i,\fn(v_i,e) \} \text{ and } f_i = \fe(v_i,f) = \{ v_i,\fn(v_i,f) \}.
    \]
    By the last part of Corollary \ref{cor: Neigbours of ancestors} the vertices $\fn(v_i,e)$ and $\fn(v_i,f)$ are not ancestors of any vertices in $V_n$. Since $e$ and $f$ are distinct vertices, it now follows from \ref{SM2} that $\fn(v_i,e) \neq \fn(v_i,f)$.

    Thus, we know that there are vertices $v_i \in V_i$ for all $i \geq n$ so that each $v_i$ is an ancestor of a vertex $v \in e \cap f$. By using the uniqueness part in Corollary \ref{cor: Neigbours of ancestors}, we have that each $v_j$ is an ancestor of $v_i$ for all $j \geq i$, and thus, we can construct a sequence $(a_i)_{i = 1}^{\infty}$ inductively so that each $e_i$ can be expressed as in \eqref{eq: Prefiber}. 
\end{proof}

An immediate corollary of the previous lemma is that $\Fib(X) \subseteq X$ is dense. 

\begin{corollary}\label{cor: Fib(X) is dense}
    $\Fib(X)$ is a dense subset of $X$.
\end{corollary}

\begin{proof}
    Let $\varepsilon > 0$ and $x = \chi((e_i)_{i = 1}^{\infty})$. Choose $n \in \N$ large enough so that $L_*^{-n} < \varepsilon$, and an edge $f_n \in E_n$ so that $v_n \in e_n \cap f_n$.
    Then $X_{e_n} \cap X_{f_n} \subseteq B(x,\varepsilon)$ by Proposition \ref{prop: Geom of LS}-\eqref{item: diams/measures of cells}, so it follows from Proposition \ref{prop: X_e cap X_f} that $\Fib(v_n) \cap B(x,\varepsilon) \neq \emptyset$.
\end{proof}

The following lemma will be useful for introducing measures on the fiber sets.

\begin{lemma}\label{lemma: How to measure on fibers}
    The fiber sets satisfy the following properties.
    \begin{enumerate}
        \vskip.3cm
        \item For all $v \in V_\#$ and $n \in \N$ we have
        \begin{equation*}
            \Fib(v) = \bigsqcup_{a \in I^n} \Fib(v \cdot a).
        \end{equation*}
        \vskip.3cm
        \label{item: How to measure on fibers (1)}
        \item For all $e \in E_{\#}$ we have $F_e(\Fib(v)) = \Fib(e \cdot v)$.
        \label{item: How to measure on fibers (2)}
        \vskip.3cm
        \item If $w \in V_{\#}$ is another vertex so that neither $v$ nor $w$ is an ancestor of the other, then
        \[
            \Fib(v) \cap \Fib(w) = \emptyset.
        \]
        \label{item: How to measure on fibers (3)}
    \end{enumerate}
\end{lemma}

\begin{proof}
    \eqref{item: How to measure on fibers (1)}: Let $f \in E_\#$ be an edge containing $v$.
    The coding map $\chi$ restricted to $\Sigma_{v,f}$ is injective by Lemma \ref{lemma: Coding map well-defined}. 
    Thus, the claim follows by noting that
    \[
        \Sigma_{v,f} = \bigsqcup_{a \in I^n} \Sigma_{v \cdot a,\fe(v \cdot a,f)}.
    \]

    \eqref{item: How to measure on fibers (2)}: If $f \in E_{\#}$ is any edge containing $v$, we compute
    \[
        \Fib(e \cdot v) = \chi(\Sigma_{e \cdot v,e \cdot f}) = \chi(\sigma_e(\Sigma_{v,f})) = F_e(\chi(\Sigma_{v,f})) = F_e(\Fib(v)).
    \]

    \eqref{item: How to measure on fibers (3)}:
    According to \eqref{item: How to measure on fibers (1)} of the current lemma, it is sufficient to prove the case where $v,w \in V_{n+2}$ and $\pi_{n+2,n}(v),\pi_{n+2,n}(w) \in V_n$ are distinct vertices for some $n \in \N \cup \{0\}$. If $e,f \in E_n$ are edges containing $\pi_{n+2,n}(v)$ and $\pi_{n+2,n}(w)$ respectively, then it follows from the non-degeneracy of the IGS and Lemma \ref{lemma: Distance between ancestors} that $d_{G_{n+2}}(v,w) \geq 4$. In particular, the edges $\fe(v,e)$ and $\fe(w,f)$ do not share a vertex, and it now follows from Proposition \ref{prop: X_e cap X_f} that
    \[
        \Fib(v) \cap \Fib(w) \subseteq X_{\fe(v,e)} \cap X_{\fe(w,f)} = \emptyset.
    \]
\end{proof}

\begin{definition}\label{def: Fiber measure}
    Let $\nu$ be a probability density on the gluing set $I$.
    For any vertex $v \in V_{\#}$ the \emph{fiber measure} $\nu_v$ on $\Fib(v)$ is defined as the Radon probability measure satisfying
    \begin{equation}\label{eq: Fiber measure}
    \nu_v(\Fib(v \cdot a)) = \prod_{i = 1}^n \nu(a_i) \text{ for all } a = a_1\dots a_n \in I^n \text{ and } n \in \N.
\end{equation}
To see that such fiber measure exist for any given $\nu$, fix $e \in E_{\#}$ with $v \in e$ and recall that $\Fib(v)$ is obtained as a homeomorphic image of $\chi \circ T_{v,e} : I^\N \to \Fib(v)$.
Now, if $\nu_{\infty}$ is the Bernoulli measure on $I^{\mathbb{N}}$ given by,
\[
    \nu_\infty\left( \{ (b_i)_{i = 1}^\infty : a_i = b_i \text{ for all } i = 1,\dots,n \} \right) := \prod_{i = 1}^n \nu(a_i) \text{ for any $n \in \mathbb{N}$, $(a_{i})_{i = 1}^{n} \in I^n$ },
\]
then the push-forward measure $\nu_v := (\chi \circ T_{v,e})_*(\nu_\infty)$ is a Radon probability measure on $\Fib(v)$ satisfying \eqref{eq: Fiber measure}. It is furthermore clear that for all $n \in \N$ and $a \in I^n$ it holds that
    \begin{equation}\label{eq: Fiber measures (Homogeneity)}
    \nu_{v}\restr_{\Fib(v \cdot a)} = \left(\prod_{i = 1}^n \nu(a_i)\right) \nu_{v \cdot a}.
    \end{equation}
\end{definition}

Let us conclude this subsection by proving Proposition \ref{prop: Geom of LS}.

\begin{proof}[Proof of Proposition \ref{prop: Geom of LS}]
    \eqref{item: diams/measures of cells}: Let $e \in E_n$.
    It is clear from the definition of the metric that by choosing $\Omega := \{e\}$, we have $\diam(X_e) \leq L_*^{-n}$. The reverse inequality follows by noting that $\dist(\Fib(e^+),\Fib(e^-),d) = L_*^{-n}$, which follows from Lemma \ref{lemma: Distance between ancestors}. 
    Now $\diam(X_v) \leq 2L_*^{-n}$ follows immediately.
    
    Next, we compute the values of the measure $\mu$.
    Notice that we only need to verify \eqref{eq: Fibers = mu-null}.
    Indeed, we would then have $\mu(X_e) = \mathfrak{m}_{\text{unif}}(\Sigma_e) = \abs{E_1}^{-n}$. According to Proposition \ref{prop: X_e cap X_f} we only need to show that $\mu(\Fib(v)) = 0$ for all $v \in V_n$. The same lemma also implies that
    \[
        \chi^{-1}(\Fib(v)) \subseteq \bigcup_{\substack{e \in E_n \\ v \in e}} \Sigma_{e},
    \]
    so we have $\mu(\Fib(v)) \leq \deg(v) \abs{E_1}^{-n}$. Using Lemma \ref{lemma: How to measure on fibers}-\eqref{item: How to measure on fibers (1)} and the fact that the degrees are uniformly bounded (Lemma \ref{lemma: Degree of ancestors}) by some number $C_{\deg}$, we get
    \[
        \mu(\Fib(v)) = \sum_{a \in I^k} \mu(\Fib(v \cdot a)) \leq C_{\deg}\abs{I}^k \abs{E_1}^{-k} \abs{E_1}^{-n}.
    \]
    Since the IGS is doubling and non-degenerate it holds that $\abs{I} < \abs{E_1}$. Therefore $\mu(\Fib(v)) = 0$ follows by letting $k \to \infty$.

    \eqref{item: Stars = balls}:
    Fix any $e \in E_n$ so that $x \in X_e$.
    In the proof of \eqref{item: diams/measures of cells} of the current proposition, we noted that $\dist(\Fib(e^-),\Fib(e^+)) = L_*^{-n}$.
    Thus, we must have $\dist(x,\Fib(w)) \geq L_*^{-n}/2$ for some $w \in e$. If $v \in e$ is the other vertex, then it follows from the definition of the metric $d$ that
    \[
      \dist(x,X \setminus X_{v}) \geq \dist(x,\Fib(w)) \geq  L_*^{-n}/2.
    \]
    Therefore $B(x,r) \subseteq B(x,3/4\cdot L_*^{-n}) \subseteq X_{v}$. The latter inclusion $X_v \subseteq B(x,8L_*r)$ is a direct corollary of $\diam(X_v) \leq 2L_*^{-n}$.

    \eqref{item: Q-AR}: $X$ is a continuous image of $\chi$ and $\Sigma$ is compact. Thus, $X$ is compact as well. Path connectedness follows from the Arzela-Ascoli theorem using the visual metric property and fact that the graphs $G_n$ are connected.
    The $Q$-Ahlfors regularity of $(X,d,\mu)$ follows from an identical argument as in \cite[Lemma 3.31]{anttila2024constructions} using the previous results of the current proposition.
    The basic idea is to approximate a ball $B(x,r)$ with appropriately chosen set $X_v$ using \eqref{item: Stars = balls} and estimate the measure by using
    \eqref{item: diams/measures of cells}.
\end{proof}

\subsection{Self-similar structure}
For future reference, we have included a discussion on the self-similarity of our construction. Specifically, we clarify that the Laakso-type fractals admit a natural self-similar structure in the sense of \cite{AnalOnFractals}. Another motivation for this discussion is that many works with similar goals use this notion, making it potentially helpful for some readers.

Let $K$ be any compact metrizable topological space, $S$ be any non-empty finite set and $\{ F_i : K \to K \}_{i \in S}$ be a collection of injective continuous maps. Denote
\[
    S^{\N} := \{ \omega_1 \omega_2 \omega_3 \ldots : \omega_i \in S \text{ for all } i \in \N \}
\]
and equip with the natural word metric.
The triplet $(K,S,\{ F_i\}_{i \in S})$ is a \emph{self-similar structure} if there is a surjective continuous map $\xi : S^{\N} \to K$ so that 
\[
    \xi(\omega_1 \omega_2 \omega_3\ldots) = F_{\omega_1}(\xi( \omega_2 \omega_3 \ldots )) \text{ for all } \omega_1 \omega_2 \omega_3\ldots \in S^{\N}.
\]

\begin{proposition}\label{prop: SSS}
    Let $(X,d)$ be a limit space of a simple and non-degenerate IGS satisfying the doubling property, $G_1 := (V_1,E_1)$ be the associated generator and $\{F_e\}_{e \in E_1}$ be the similarity maps. Then $(X,E_1,\{ F_e \}_{e \in E_1})$ is a self-similar structure.
\end{proposition}
\begin{proof}
Recall from the discussion in Definition \ref{def: Measures on sigma} that we have a natural identification $T : E_1^{\N} \to \Sigma$. By setting $\xi := \chi \circ T : E_1^{\N} \to X$ where $\chi : \Sigma \to X$ is the coding map in Definition \ref{def: Limit space},
it follows from Proposition \ref{prop: LS well-defined}-\eqref{item: Similarity maps} and \eqref{def: similarity maps} that $(X,E_1,\{F_e\}_{e \in E_1})$ is a self-similar structure.
\end{proof}

We also clarify the two properties of these self-similar structures.
These easily follow from the definitions so we omit the details. The undefined terms can be found in \cite{AnalOnFractals}.

\begin{corollary}
    Let $(X,d), G_1, \{F_e\}_{e \in E_1}$ be as in Proposition \ref{prop: SSS}. 
    The boundary $V_{0}(\mathcal{L})$ of the self-similar structure $\mathcal{L} \coloneqq (X,E_1,\{F_e\}_{e \in E_1})$ is $\Fib(v_+) \sqcup \Fib(v_-)$.
\end{corollary}

\begin{corollary}
    Let $(X,d), G_1, \{F_e\}_{e \in E_1}$ be as in Proposition \ref{prop: SSS}. 
    The self-similar structure $(X,E_1,\{F_e\}_{e \in E_1})$ is post-critically finite if and only if $\abs{I} = 1$.
\end{corollary}

\section{Discrete Potential theory}\label{sec: DPT}
In this section we study discrete potential theory of IGSs. This is an important preliminary step for constructing the $p$-energy forms to the limit spaces in the later sections because our strategy is to construct them as limit of discrete ones.
The reader may recall the terminology and basic results of discrete potential theory from Subsection \ref{subsec: (preli) DPT}.

Let us begin by recording the required assumptions of our method, and some relevant constants. The first three conditions were introduced in Section \ref{sec:ssLaakso} and the last one, the conductively uniform property, is discussed later in the current section.

\begin{assumption}\label{Assumptions: To make work}
    The IGS consisting of the data $V_1,E_1,I,\{ \phi_{v,e} \}_{v \in e \in E_1}$ and the associated generator $G_1$ satisfy the following conditions.
    \begin{enumerate}
    \vskip.3cm
        \item \textup{\textbf{(Simplicity)}} There are functions $\phi_{+},\phi_- : I \to V_1$ so that for all edges $e = \{v,w\}  \in E_1$ it holds that
        \[
          \{ \phi_{v,e}, \phi_{w,e} \} = \{ \phi_+,\phi_- \}.  
        \]
        We also denote $I_+ := \phi_+(I)$ and $I_- := \phi_-(I)$.
        \vskip.3cm
        \item \textup{\textbf{(Non-degeneracy)}} There is no edge between the sets $I_+$ and $I_-$.
        \vskip.3cm
        \item \textup{\textbf{(Doubling)}} For all $v \in I_+ \cup I_-$ it holds that $\deg(v) = 1$. The unique neighbor of such $v$ is denoted $\fn(v)$.
        \vskip.3cm
        \item \textup{\textbf{(Conductively uniform)}} For all $p \in (1,\infty)$ and $a \in I$ it holds that
        \[
            \mathcal{J}_{p,+}(\phi_+(a),\fn(\phi_+(a))) = -\mathcal{J}_{p,+}(\phi_-(a),\fn(\phi_-(a))).
        \]
    Here, $\mathcal{J}_{p,+}$ is the $p$-energy minimizing unit flow from $I_+$ to $I_-$ in $G_1$.
    \end{enumerate}
    We also fix the following constants that depend on at most the initial data associated to the IGS and the exponent $p \in (1,\infty)$.
    \begin{enumerate}[(i)]
    \vskip.3cm
        \item \textup{\textbf{(Geometric constants)}}
        \[
        C_{\diam} := \diam(V_1,d_{G_1}), \quad C_{\deg} := \max_{v \in V_1} \deg(v), \quad  L_* := \dist(I_+,I_-,d_{G_1}).
        \]
        \item \textup{\textbf{(Hausdorff dimension)}}
        \[
            Q := \frac{\log(\abs{E_1})}{\log(L_*)}.
        \]
        \item \textup{\textbf{($p$-capacity constant)}}
        \begin{equation*}\label{eq: Capacity constant}
            \cM_p := \cCap_p(I_+,I_-,G_1).
        \end{equation*}
        \item
        \textup{\textbf{($p$-walk dimension)}}
        \begin{equation*}
            \pwalk := \frac{\log(\abs{E_1}\cM_p^{-1})}{\log(L_*)}.
        \end{equation*}
    \end{enumerate}
    Lastly, during the computations, we will frequently use the notation $A \lesssim B$ to indicate the existence of a constant $C \geq 1$ so that $A \leq C \cdot B$, where $C$ depends on some inessential parameters.
    For the most part, these parameters are the constants mentioned above and $p$.
\end{assumption}

\subsection{Conductive uniform property}
We introduce the final general assumption of our framework, the \emph{conductively uniform property}. This was a key ingredient for establishing the combinatorial Loewner property in \cite{anttila2024constructions}.

We need some potential theoretic notation first. It may be helpful for the reader to revisit Figures \ref{intro: potential/flow} and \ref{intro: flow}.
Suppose we are given a simple IGS: Let $n \in \N \cup\{0\}$ and $p \in (1,\infty)$. We denote $U_{p,+,n}$ the solution to the $p$-capacity problem
\begin{equation}\label{eq: capacity problem}
    \cCap_p\left(I_{+}^{(n)},I_{-}^{(n)},G_n \right).
\end{equation}
The \emph{$p$-capacity constant} 
$\cM_p$ is the value of the $p$-capacity problem \eqref{eq: capacity problem} for $n = 1$, i.e., $\cM_p := \mathcal{E}_p(U_{p,+,1})$.
Recall from Lemma \ref{lemma: Potential function exists} that $U_{p,+,n}$ exists and is unique.
The solution to the $p$-capacity problem \eqref{eq: capacity problem} where roles of the signs $+$ and $-$ are interchanged is denoted $U_{p,-,n}$.
These functions are collectively referred as the \emph{optimal potential functions}.

We denote $\mathcal{J}_{p,+,n}$ the solution to the $p$-resistance problem
\begin{equation}\label{eq: resistance problem}
    \cRes_p\left(I_{+}^{(n)},I_{-}^{(n)},G_n \right).
\end{equation}
Similarly, we denote $\mathcal{J}_{p,-,n}$
the solution to \eqref{eq: resistance problem} where the signs are interchanged. According to Lemma \ref{lemma:duality}, the flow $\mathcal{J}_{p,\pm,n}$ exists and is unique. These are called \emph{optimal flows}.

By the uniqueness of the optimal potential function, we have $U_{p,\pm,n} = 1 - U_{p,\mp,n}$.
In particular, the gradient $\abs{\nabla U_{p,\pm,n}}$ does not depend on the choice of the sign. We thus denote $\abs{\nabla U_{p,n}} := \abs{\nabla U_{p,\pm,n}}$ and $\abs{\nabla U_{p}} := \abs{\nabla U_{p,1}}$. Similarly, the energies of the optimal potentials are denoted as $\mathcal{E}_p(U_{p,n})$ and $\mathcal{E}_p(U_{p})$.
In the case of the flows, since $\mathcal{J}_{p,\pm,n} = -\mathcal{J}_{p,\mp,n}$, we slightly abuse the notation by writing $\mathcal{J}_{p,n} :=  \mathcal{J}_{p,\pm,n}$ and $\mathcal{J}_{p} := \mathcal{J}_{p,\pm} := \mathcal{J}_{p,\pm,1}$ whenever the sign has no role during the computations or is understood from the context.

\begin{example}
    Consider the IGSs in Examples \ref{example: Laakso space} and \ref{example: Counterexample}, and let $U : V_1 \to \R$ be the function in Figure \ref{fig: Optimal potential/flows Laakso/Kleiner} (note that the function is the same in both graphs). To see that $U$ solves the $p$-capacity problems $\cCap_p(I_+,I_-,G_1)$ for both IGSs, it is direct to verify that $U$ is $p$-harmonic in $V_{1} \setminus (I_+ \cup I_-)$ for all $p \in (1,\infty)$. 
    \begin{figure}[!ht]
    \includegraphics[width=315pt]{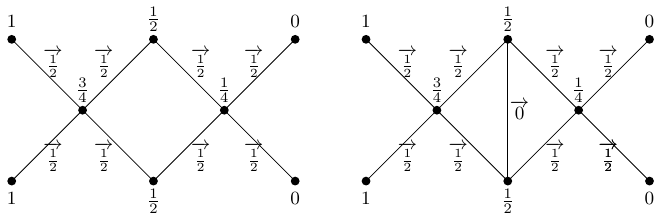}
    \caption{Figure of the optimal potentials/flows on $G_1$ of the IGSs in  Examples \ref{example: Laakso space} and \ref{example: Counterexample}. The values next to vertices indicate the value of the optimal potential $U_{p,+}$, and the values next to the edges indicate the values of the optimal flow $\mathcal{J}_{p,+}$.}
    \label{fig: Optimal potential/flows Laakso/Kleiner}
    \end{figure}
    According to Lemma \ref{lemma: Potential function exists}, $U$ is the $p$-energy minimizer. The optimal flows are clearly the ones in Figure \ref{fig: Optimal potential/flows Laakso/Kleiner} by the symmetry of these graphs. Alternatively, the flow can be computed by duality \eqref{eq:duality (potentials and flows)}.
    See Figures \ref{intro: potential/flow} and \ref{fig: singularities not equiv} for other examples.
\end{example}

\begin{definition}\label{def: Cond.unif}
    We say that a simple and non-degenerate IGS satisfying the doubling property is \emph{conductively uniform}
    if for all $a \in I$ and $p \in (1,\infty)$
    \begin{equation}\label{eq: Flows agree on boundary}
        \mathcal{J}_{p,+}(\phi_{+}(a),\fn(\phi_+(a))) = -\mathcal{J}_{p,+}(\phi_{-}(a),\fn(\phi_-(a))).
    \end{equation}
    Here $\fn(z)$ denotes the unique neighbor of $z \in I_+ \cup I_-$ as defined in Definition \ref{def: doubling}. 
    Whenever this condition is satisfied, the value in \eqref{eq: Flows agree on boundary} is denoted as $\mathcal{J}_p(a)$.
\end{definition}

\begin{remark}\label{Remark: Measure on gluing set}
    $\mathcal{J}_p(a)$ is a probability density on $a \in I$. Indeed, by the doubling property, every vertex in $I_+$ has degree 1, so we have
    \[
        1 = \sum_{a \in I} \divr(\mathcal{J}_{p,+})(\phi_+(a)) = \sum_{a \in I} \mathcal{J}_{p,+}(\phi_+(a),\fn(\phi_+(a))) =
        \sum_{a \in I} \mathcal{J}_p(a).
    \]
The strong maximum principle of $p$-harmonic functions (Lemma \ref{lemma: Strong max principle}) and the duality \eqref{eq:duality (potentials and flows)} (note that here we use the version that does note involve absolute values) implies that  $\mathcal{J}_p(a) > 0$.
Consequently, $\mathcal{J}_{p}(\cdot)$ produces a probability density on $I^n$, denoted as $\mathcal{J}_{p,n}(\cdot)$, which is given by
\begin{equation}\label{eq: defn-measureIn}
    \mathcal{J}_{p,n}(a_1\dots a_n) := \prod_{i = 1}^n \mathcal{J}_p(a_i).
\end{equation}
For $n = 0$, we intepret $I^0$ to consist of a single (empty) element, and simply set $\mathcal{J}_{p,0}(\cdot) := 1$.
\end{remark}

If we view the $p$-capacity problem \eqref{eq: capacity problem} as a discrete optimization problem, the number of variables grow exponentially with respect to $n$.
The main motivation of conductively uniform property is to reduce this problem to the case $n = 1$.
This property is very helpful when we take a limit of discrete energies.
Nevertheless, the conductively uniform property is by far the most restrictive condition in our framework. In the following theorem, we review a few sufficient conditions for it.

\begin{theorem}\label{thm: Symmetric IGS}
    Suppose that the IGS is simple, non-degenerate and satisfies the doubling property.
    Then any one of the following conditions are sufficient for the IGS to satisfy the conductive uniform property.
    \begin{enumerate}[label={\textup{(\textcolor{blue}{CUP-\arabic*})}}, widest=a, leftmargin=*]
    \vskip.3cm
    \item \label{item: Sufficient conditions (PCF)}
    The gluing set $I$ only contains one element.
    \vskip.3cm
    \item \label{item: Sufficient conditions (Symmetry)} 
    There is a graph isomorphism $\eta : V_1 \to V_1$ of $G_1$ satisfying
    \begin{equation*}
        \eta(\phi_{\pm}(a)) = \phi_{\mp}(a) \text{ for all } a \in I.
    \end{equation*}
    \vskip.3cm
    \item \label{item: Sufficient conditions (dwp = p)}
    There exists a family of paths $\Theta$ in $G_1$ so that each $\theta \in \Theta$ connects $I_+$ to $I_-$ and is of length $L_* = \dist(I_+,I_-,d_{G_1})$, and each edge in $e \in E_1$ is contained in exactly one of the paths $\theta \in \Theta$.
\end{enumerate}
\end{theorem}

\begin{proof}
    Assume first that \ref{item: Sufficient conditions (PCF)} holds and $I = \{ a \}$. Since $\mathcal{J}_{p,+}$ is a unit flow from $I_+ = \{\phi_+(a)\}$ to $I_-=\{\phi_-(a)\}$, we have
    \[
    \mathcal{J}_{p,+}(\phi_{+}(a),\fn(\phi_{+}(a))) = 1 = \mathcal{J}_{p,+}(\phi_{-}(a),\fn(\phi_{-}(a))).
    \]

    Next assume \ref{item: Sufficient conditions (Symmetry)}, and let $\eta$ be a graph isomorphism as in the claim.
    Since $\eta$ satisfies $\eta(I_\pm) = I_\mp$, it holds that
    \[
        \eta(\mathcal{J}_{p,+})(v,w) := \mathcal{J}_{p,+}(\eta(v),\eta(w))
    \]
    is a unit flow from $I_{-}$ to $I_+$ satisfying $\mathcal{E}_q(\eta(\mathcal{J}_{p,\pm})) = \mathcal{E}_q(\mathcal{J}_{p,\mp})$. It follows from the uniqueness of the energy minimizing flow that
    $\eta(\mathcal{J}_{p,+}) = \mathcal{J}_{p,-}$. 
    By the doubling property, we necessarily have $\eta(\fn(\phi_-(a))) = \fn(\phi_+(a))$.
    Thus, we compute
    \begin{align*}
        \mathcal{J}_{p,+}(\phi_-(a),\fn(\phi_-(a))) & = \eta(\mathcal{J}_{p,+})(\phi_+(a),\fn(\phi_+(a))\\
        & = \mathcal{J}_{p,-}(\phi_+(a),\fn(\phi_+(a))\\
        & = -\mathcal{J}_{p,+}(\phi_+(a),\fn(\phi_+(a)).
    \end{align*}

    Lastly, we assume \ref{item: Sufficient conditions (dwp = p)}, and let $\Theta$ be a family of paths satisfying the therein conditions.
    For $a \in I$ let $\theta_a \in \Theta$ be the path containing the edge $\fe(\phi_+(a))$. Note that $a \mapsto \theta_a, I \to \Theta$ is a bijection by the properties of $\Theta$. Thus, it holds that
    \begin{equation}\label{eq: E_1 = glued intervals}
        \abs{E_1} = \abs{I}\cdot L_*.
    \end{equation}

    Now, for every $a \in I$, consider the unit flow $\mathcal{J}_a$ from $I_+$ to $I_-$ along the path $\theta_a$. Then we define the unit flow from $I_+$ to $I_-$, 
    \[
        \mathcal{J} := \abs{I}^{-1} \sum_{a \in I} \mathcal{J}_a.
    \]
    Using the fact that the paths in $\Theta$ are edge-wise disjoint, we get
    \[
    \mathcal{J}(\phi_+(a),\fn(\phi_+(a))) = \abs{I}^{-1} = \mathcal{J}(\fn(\phi_-(a)),\phi_-(a)).
    \]
    Therefore, the conductively uniform property follows after we show that $\mathcal{J}$ is the $p$-energy minimizing unit flow. To see this, it follows from the properties of $\Theta$ that $\abs{\mathcal{J}} \equiv \abs{I}^{-1}$.
    By the duality \eqref{eq:duality (potentials and flows)} and \eqref{eq: E_1 = glued intervals},
    \[
        \cM_p \geq \mathcal{E}_q(\mathcal{J})^{\frac{-p}{q}} = \left(\abs{E_1} \abs{I}^{-q} \right)^{\frac{-p}{q}} = \frac{\abs{E_1}}{L_*^p}.
    \]
    In Lemma \ref{lemma: Properties of cond.unif.}-\eqref{item: Properties of cond.unif. (3)} below we show that the reverse inequality $\cM_p \leq \abs{E_1}/L_*^p$ always holds. By using the duality one last time, we see that $\mathcal{J}$ is the $p$-energy minimizing unit flow.
\end{proof}

Let us assume for the rest of the section that the conditions in Assumption \ref{Assumptions: To make work} are satisfied by the IGSs.
We shall now review a few important consequences of the conductively uniform property.

\begin{lemma}\label{lemma: Properties of cond.unif.}
    Let $p \in (1,\infty)$. Then the following hold.
    \begin{enumerate}
    \vskip.3cm
    \item \label{item: Properties of cond.unif. (1)}
    For all $n \in \N \cup \{0\}$ we have
    \[
        \cCap_p\left(I_{+}^{(n)},I_{-}^{(n)},G_n\right) = \cM_p^n.
    \]
    \vskip.3cm
        \item \label{item: Properties of cond.unif. (2)}
        For all $n \in \N$ and $a \in I^n$ we have
        \begin{equation*}
        \divr\left(\mathcal{J}_{p,\pm,n}\right)(\phi_{\pm,n}(a)) = \mathcal{J}_{p,n}(a) = -\divr\left(\mathcal{J}_{p,\pm,n}\right)(\phi_{\mp,n}(a)).
    \end{equation*}
    Here $\phi_{\pm,n}$ are the higher order gluing maps in Definition \ref{def: n-Gluing maps} and $\mathcal{J}_{p,n}(\cdot)$ is the probability density given in Remark \ref{Remark: Measure on gluing set}.
    \vskip.3cm
    \item \label{item: Properties of cond.unif. (3)} The $p$-capacity constant has the upper bound $\cM_p \leq \abs{E_1}/L_*^p$, or equivalently, $d_{w,p} \geq p$. These inequalities are equalities if and only if $\abs{\nabla U_p} \equiv L_*^{-1}$. Here, $d_{w,p}$ is the $p$-walk dimension in Assumption \ref{Assumptions: To make work}.
    \vskip.3cm
    \item \label{item: Properties of cond.unif. (4)}
    The optimal potential functions $U_{p,\pm,n}$ satisfy the following strong maximum principle: If $z \in V_{n} \setminus \left(I_{+}^{(n)} \cup I_-^{(n)}\right)$ then $0 < U_{p,\pm,n}(z) < 1$.
    \end{enumerate}
\end{lemma}

\begin{proof}
    \eqref{item: Properties of cond.unif. (1)}: This was proven in \cite[Corollary 4.32]{anttila2024constructions} for discrete (edge) $p$-modulus, which is always equal to discrete $p$-capacity (see e.g. \cite[Lemma 4.13]{anttila2024constructions}).
    %Let us sketch the argument. The inequality ''$\leq$'' essentially follows from the proof of Theorem \ref{thm: Expression of optimal potential}. To prove the other inequality ''$\geq$'', we construct a flow suitable flow and apply the duality Lemma \ref{lemma:duality}. This construction heavily depends on the conductively uniform property.
    
    \eqref{item: Properties of cond.unif. (2)}: This was verified in the proof of \cite[Proposition 4.30]{anttila2024constructions}.

    \eqref{item: Properties of cond.unif. (3)}:
    A variant of this claim was proven in \cite[Proposition 6.11]{anttila2024constructions}.
    Let us present the details.
    Consider the potential function $U : V_1 \to \R$ given by
    \[
        U(v) := \frac{\dist(v,I_-,d_{G_1})}{L_*} \land 1.
    \]
    Since $L_* = \dist(I_+,I_-,d_{G_1})$, we 
    clearly have $U|_{I_+} \equiv 1$ and $U|_{I_-} \equiv 0$. Moreover, for any edge $\{ v,w \} \in E_1$ we have
    \[
        \abs{U(v)-U(w)} \leq \frac{1}{L_*}\abs{\dist(v,I_-) - \dist(w,I_-)} \leq \frac{1}{L_*}.
    \]
    Thus, we obtain the inequality
    \[
        \cM_p \leq \mathcal{E}_p(U) \leq \frac{\abs{E_1}}{L_*^p}.
    \]
    The latter part of the claim now follows from the uniqueness of the optimal potential function.

    \eqref{item: Properties of cond.unif. (4)}: Since the graphs $G_n$ are connected and the IGS is doubling, it follows from Lemma \ref{lemma: Degree of ancestors} that $V_n \setminus (I_{+}^{(n)} \cup I_-^{(n)})$ is a connected subset. Thus, the claim follows from the strong maximum principle of $p$-harmonic functions (Lemma \ref{lemma: Strong max principle}).
\end{proof}

\begin{theorem}\label{thm: Expression of optimal potential}
    For all $n,m \in \N \cup \{ 0 \}$ the optimal potential function $U_{p,\pm,n+m}$ can be expressed as
    \begin{equation}\label{eq: Expression of optimal potential}
        U_{p,\pm,n+m}(e \cdot v) = U_{p,\pm,n}(e^-) + (U_{p,\pm,n}(e^+) - U_{p,\pm,n}(e^-)) \cdot U_{p,+,m}(v)
    \end{equation}
    for all $e \in E_n$ and $v \in V_m$.
    In particular, for all $v \in V_n$ it holds that
    \begin{equation}\label{eq: Optimal potential on Fibers}
        U_{p,\pm,n+m}(w) = U_{p,\pm,n}(v) \text{ for all } w \in \pi_{n+m,n}^{-1}(v).
    \end{equation}
\end{theorem}

\begin{proof}
Let $\widehat{U}_{p,\pm,n+m}$ be the function given by the right-hand side of \eqref{eq: Expression of optimal potential}. First, we verify that it is well-defined.
Assume that $e \cdot v = e' \cdot v'$ for $e,e' \in E_n$ and $v,v' \in V_m$. By \ref{SM2} the edges $e,e'$ have a common vertex $w \in e \cap e'$ so that $v \in I_{w,e}^{(n)}$ and $v' \in I_{w,e'}^{(n)}$. It then follows by the definition of $\widehat{U}_{p,\pm,n+m}$
\begin{equation}\label{eq: to prove expression of optimal potential}
    \widehat{U}_{p,\pm,n+m}(e \cdot v) = U_{p,\pm,n}(w) = \widehat{U}_{p,\pm,n+m}(e' \cdot v'),
\end{equation}
and $\widehat{U}_{p,\pm,n+m}$ is thus well-defined. It is also clear from \eqref{eq: to prove expression of optimal potential} that $\widehat{U}_{p,\pm,n+m}$ is potential function to the $p$-capacity problem.
To conclude that $\widehat{U}_{p,\pm,n+m} = U_{p,\pm,n+m}$, we thus only need to show that $\mathcal{E}_{p}(\widehat{U}_{p,n+m}) = \mathcal{E}_{p}(U_{p,n+m})$. The desired equality then follows from the uniqueness of the optimal potential function (Lemma \ref{lemma: Potential function exists}).
Consequently, \eqref{eq: Optimal potential on Fibers} would follow from \eqref{eq: to prove expression of optimal potential}.

Using \ref{SM3} (or alternatively Lemma \ref{lemma: Discrete energy (Self-similarity)}) and Lemma \ref{lemma: Properties of cond.unif.}-\eqref{item: Properties of cond.unif. (1)}, we compute
    \begin{align*}
        \cE_p(\widehat{U}_{p,\pm,n+m}) & = \sum_{ e_n \in E_n } \sum_{ e_m \in E_m } \abs{\nabla \widehat{U}_{p,\pm,n+m}(e_n \cdot e_m) }^p\\
        & = \sum_{ e_n \in E_n } \abs{\nabla U_{p,n}(e_n)}^p \sum_{ e_m \in E_m } \abs{\nabla U_{p,m}(e_m)}^p\\
        & = \cM_p^{n+m} =  \cE_p(U_{p,n+m}).
    \end{align*}
\end{proof}

\begin{corollary}\label{cor: Gradient product}
    Let $e_1,e_2,\dots,e_n \in E_1$ and $e := e_1e_2\ldots e_n \in E_n$. Then
    \begin{equation}\label{eq: Gradient product}
        \abs{\nabla U_{p,n}(e)} = \prod_{i = 1}^n \abs{\nabla U_{p}(e_i)} \, \text{ and } \, \abs{\mathcal{J}_{p,n}(e)} = \prod_{i = 1}^n \abs{ \mathcal{J}_{p}(e_i)}.
    \end{equation}
    Here the sequence $e_1e_2\ldots e_n$ is understood as in Remark \ref{remark: Identify as words}.
\end{corollary}

\begin{proof}
    The first equality in \eqref{eq: Gradient product} follows from a direct iteration of \eqref{eq: Expression of optimal potential}, and the second would then follow from the duality \eqref{eq:duality (potentials and flows)} and Lemma \ref{lemma: Properties of cond.unif.}-\eqref{item: Properties of cond.unif. (1)}.
\end{proof}
    
\begin{corollary}\label{Cor: Discrete hölder continuity}
    Let $p \in (1,\infty)$ and set
    \[
        c_p := \max_{ e \in E_1 } \, \abs{\nabla U_p(e)}.
    \]
    For all $n,m \in \N \cup \{0\}$ we have
    \begin{equation}
        \abs{U_{p,\pm,n+m} (e \cdot v_1) - U_{p,\pm,n+m}(e \cdot v_2)} \leq c_p^n
    \end{equation}
    for all $e \in E_n$ and $v_1,v_2 \in V_m$.
\end{corollary}

\begin{proof}
    It follows from Corollary \ref{cor: Gradient product} that
    \[
        \abs{\nabla U_{p,n}(e)} \leq c_p^n \text{ for all } e \in E_n.
    \]
    By using the strong maximum principle (Lemma \ref{lemma: Properties of cond.unif.}-\eqref{item: Properties of cond.unif. (4)}), we see that
    \[
        |U_{p,\pm,m}(v_1) - U_{p,\pm,m}(v_2) | \leq 1 \text{ for all } v_1,v_2 \in V_m.
    \]
    By combining the previous estimates and using \eqref{eq: Expression of optimal potential}, we obtain
    \begin{align*}
        &  \, \left|U_{p,\pm,n+m}(e \cdot v_1)-U_{p,\pm,n+m}(e \cdot v_2)\right|\\
        = & \, \, \abs{\nabla U_{p,n}(e)} \cdot | U_{p,\pm,m}(v_1) - U_{p,\pm,m}(v_2) |\\
        \leq & \, \, \abs{\nabla U_{p,n}(e)}
        \leq  c_p^n.
    \end{align*}
\end{proof}

We finish the subsection by discussing the equality $d_{w,p} = p$ from Lemma \ref{lemma: Properties of cond.unif.}-\eqref{item: Properties of cond.unif. (3)}. We show that it has a geometric characterization that is \emph{independent} of $p$.
\begin{proposition}\label{prop: dwp = p}
    The equality $d_{w,p} = p$ is equivalent to the existence of a family of paths $\Theta$ as in \ref{item: Sufficient conditions (dwp = p)}. In particular, whether $d_{w,p} = p$ holds is independent of the exponent $p$.
\end{proposition}

\begin{proof}
    In the proof of Theorem \ref{thm: Symmetric IGS} we showed that \ref{item: Sufficient conditions (dwp = p)} implies $\cM_p = \abs{E_1}/L_*^p$, which by definition means $d_{w,p} = p$.
    Let us then assume the equality $d_{w,p} = p$.
    For any vertex $v \in V_1$ we divide the edges containing $v$ into
    \[
    E_+(v) := \{ \{v,w\} \in E_1 : U_{p,+}(v) > U_{p,+}(w) \},
    \]
    and the rest are denoted $E_-(v)$.
    Since $\abs{\nabla U_p} \equiv L_*^{-1}$ holds by Lemma \ref{lemma: Properties of cond.unif.}-\eqref{item: Properties of cond.unif. (3)}, it follows from the $p$-harmonicity of $U_{p,+}$,
    \begin{equation}\label{eq: dwp=p same amount}
        \abs{E_+(v)} = \abs{E_-(v)} \text{ for any } v \in V_1 \setminus (I_+ \cup I_-)
    \end{equation}
    
    We now have enough ingredient to describe an algorithm that constructs the family of paths $\Theta := \{ \theta_1,\dots,\theta_{\abs{I}} \}$.
    First, fix $v_0 \in I_+$ and choose any path $\theta_1 = [v_0,v_1\dots,v_k]$ so that $U_{p,+}(v_i)$ is strictly decreasing in $i$ and $v_k \in I_-$. Such choice is possible due to the strong maximum principle (Lemma \ref{lemma: Properties of cond.unif.}-\eqref{item: Properties of cond.unif. (4)}) and $p$-harmonicity of $U_{p,+}$.
    Further, since $\abs{\nabla U_p} \equiv L_*^{-1}$, we have $k = L_*$.
    
    Next, to construct $\theta_2$, we choose any vertex $v_0' \in I_+ \setminus \{v_0\}$ and a path $\theta_2 = [v_0',v_1'\dots,v_k']$ satisfying the same conditions as $\theta_1$.
    Additionally, we require that $\theta_1$ and $\theta_2$ are edge-wise disjoint. Such choice is possible thanks to \eqref{eq: dwp=p same amount}. By repeating, we obtain the paths $\Theta := \{ \theta_1,\dots, \theta_{\abs{I}} \}$. It follows from \eqref{eq: dwp=p same amount} and the doubling property of the IGS that the paths in $\Theta$ contain all edges.
\end{proof}

\begin{remark}\label{rem:CK}
    This remark continues the discussion in Subsection \ref{subsec:CK}
    In what follows from Proposition \ref{prop: dwp = p}, the current IGS framework and Cheeger--Kleiner \cite{CK_PI} intersect precisely at the IGSs admitting a family of paths $\Theta$ where \ref{item: Sufficient conditions (dwp = p)} holds. This is because the Hausdorff measure would satisfy the axioms of \cite{CK_PI}.
    \begin{figure}[!ht]
    \centering\includegraphics[width=235pt]{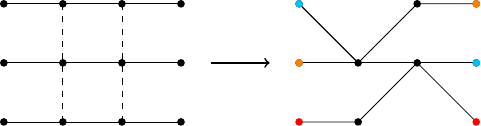}
    \caption{Construction that produces examples satisfying \ref{item: Sufficient conditions (dwp = p)}.
    }
    \label{fig: genlaakso}
    \end{figure}
    The frameworks agree also in Sobolev spaces; see Theorem \ref{Thm: NS-equivalence}.
    Moreover, all Examples satisfying \ref{item: Sufficient conditions (dwp = p)} can be constructed by taking interval graphs of same length and applying ``vertical identifications''; see Figure \ref{fig: genlaakso}. This procedure closely resembles Laakso's construction in \cite{La00} in the sense that the vertical identifications correspond to the ``wormholes'' in Laakso's work.
\end{remark}

\subsection{Strong monotonicity}\label{Subsec: Strong monotonicity}
For our goals, it is crucial to understand the interrelation between the discrete energies at different levels.
In other words, we need a convenient way to translate functions between graphs of different levels, and effectively compare their energies. To this end, we introduce the \emph{averaging operators} $V_{p,n,m}[\cdot] : \R^{V_n} \to \R^{V_m}$ for $n \geq m$ defined as
\begin{equation}\label{eq: Discrete average}
    V_{p,n,m}[f](v) := \sum_{a \in I^{n - m}} \mathcal{J}_{p,n-m}(a) \cdot f(v \cdot a),
\end{equation}
where $\mathcal{J}_{p,k}(\cdot)$ is the probability density on $I^k$ defined in Remark \ref{Remark: Measure on gluing set}. Note that these operators in general depend on $p$. However, we omit the exponent from the notation for the most part, and instead write $V_{n,m}[\cdot]$.

While the expression \eqref{eq: Discrete average} is clearly valid for any probability density on the gluing set $I$, our specific choice is one of the most fundamental pieces of our framework.
We will prove in Theorem \ref{thm: Discrete energy monotonicity} that these operators satisfy a highly useful analytic condition we call \emph{strong monotonicity}. First, we review a few lemmas.  

\begin{lemma}\label{lemma: Average compatible-discrete} 
    Let $n,m,l \in \mathbb{N} \cup \{ 0 \}$ and $m \le l \le n$. Then $V_{n,m}[\cdot] = V_{l,m}[\cdot] \circ V_{n,l}[\cdot]$. 
\end{lemma}

\begin{proof}
    This is immediate from \eqref{eq: defn-measureIn}. Indeed, for any $f \in \mathbb{R}^{V_{n}}$ and $v \in V_{n}$, by using the associativity property, Lemma \ref{lemma: Symbolic associativity}-\eqref{item: Symbolic associativity (4)}, 
    \begin{align*}
        V_{n,m}[f](v) 
        &= \sum_{a \in I^{n-l}}\sum_{b \in I^{l - m}}\mathcal{J}_{p,n-l}(a)\mathcal{J}_{p,l-m}(b)f(v \cdot (a \cdot b)) \\
        &= \sum_{a \in I^{n-l}}\sum_{b \in I^{l - m}}\mathcal{J}_{p,n-l}(a)\mathcal{J}_{p,l-m}(b)f((v \cdot a) \cdot b) \\
        &= \sum_{a \in I^{n-l}}\mathcal{J}_{p,n-l}(a)V_{n,l}[f](v \cdot a) 
        = V_{l,m}[V_{n,l}[f]](v), 
    \end{align*}
    proving $V_{n,m} = V_{l,m} \circ V_{n,l}$. 
\end{proof}

\begin{lemma}\label{lemma: Average commutes with similarity}
    Let $n,m,k \in \N \cup \{0\}$ so that $n \geq m$. Then for all $f : V_{n+k} \to \R$ and $e \in E_k$ we have
    \begin{equation}\label{eq: Average commutes with similarity}
        V_{n,m}[f \circ \sigma_{e,n}] = V_{n+k,m+k}[f] \circ \sigma_{e,m}
    \end{equation}
\end{lemma}

\begin{proof}
    By Lemma \ref{lemma: Symbolic associativity}-\eqref{item: Symbolic associativity (3)} we have $\sigma_{e,n}(v \cdot a) = \sigma_{e,m}(v)\cdot a$ for all $v \in V_m$ and $a \in I^{n-m}$. Thus
    \begin{align*}
         V_{n,m}[f \circ \sigma_{e,n}](v) & = \sum_{a \in I^{n-m}}\mathcal{J}_{p,n-m}(a) f(\sigma_{e,n}(v \cdot a))\\
         & = \sum_{a \in I^{n-m}} \mathcal{J}_{p,n-m}(a)  f(\sigma_{e,m}(v)\cdot a)\\
         & = (V_{n+k,m+k}[f] \circ \sigma_{e,m})(v).
    \end{align*}
\end{proof}

\begin{lemma}\label{lemma: Discrete energy (Self-similarity)}
    Let $n,m \in \N \cup \{0\}$ so that $n \geq m$. Then for all $f : V_n \to \R$
    \begin{equation}\label{eq: Discrete energy (Self-similarity)}
        \cE_p(f) = \sum_{e \in E_m} \cE_p(f \circ \sigma_{e,n-m}). 
    \end{equation}
\end{lemma}

\begin{proof}
    Directly follows from \ref{SM3}.
\end{proof}

\begin{theorem}[Strong monotonicity]\label{thm: Discrete energy monotonicity}
    Let $n,m \in \N \cup \{ 0 \}$ so that $n \geq m$ and $f : V_n \to \R$. Then
    \begin{equation}\label{eq: Discrete energy monotonicity}
        \cE_p(V_{n,m}[f]) \leq \cM_p^{-(n-m)}\cE_p(f),
    \end{equation}
    where $\cM_p$ is the $p$-capacity constant in Assumption \ref{Assumptions: To make work}.
\end{theorem}

\begin{proof}
    There is nothing to prove if $n = m$, so we assume $n > m$. First, choose any edge $e \in E_{m}$. Using Lemma \ref{lemma: Properties of cond.unif.}-\eqref{item: Properties of cond.unif. (2)} and the fact that $\mathcal{J}_{p,n-m} := \mathcal{J}_{p,+,n-m}$ is a unit flow from $I_+^{(n-m)}$ to $I_-^{(n-m)}$, we have for all $v \in V_{n-m}$ that
    \begin{equation}\label{eq: Divergences}
        \divr(\mathcal{J}_{p,n-m})(v) =
        \begin{cases}
            \mathcal{J}_{p,n-m}(a) & \text{ if } v = \phi_{+,n-m}(a) \text{ for } a \in I^{n-m}\\
            -\mathcal{J}_{p,n-m}(a) & \text{ if } v = \phi_{-,n-m}(a) \text{ for } a \in I^{n-m}\\
            0 & \text{ otherwise}.
        \end{cases}
    \end{equation}
    By applying this, we have
    \begin{align*}
        & \, \quad \abs{V_{n,m}[f](e^+) - V_{n,m}[f](e^-)}\\
        = & \quad \left| \sum_{a \in I^{n-m}} \mathcal{J}_{p,n-m}(a) f(e^+ \cdot a) - \sum_{a \in I^{n-m}} \mathcal{J}_{p,n-m}(a) f(e^- \cdot a) \right|\\
        = & \quad \left|\sum_{v \in V_{n-m}} \divr(\mathcal{J}_{p,n-m})(v) (f \circ \sigma_{e,n-m})(v)\right| && \text{(Equation \eqref{eq: Divergences})} \\
        = & \quad \left| \sum_{ \{v,w\} \in E_{n-m} } \mathcal{J}_{p,n-m}(v,w) \nabla (f \circ \sigma_{e,n-m})(v,w) \right| && \text{(Divergence theorem \eqref{eq:divthm})} \\
        \leq & \quad \cE_q(\mathcal{J}_{p,n-m})^{\frac{1}{q}} \cE_p(f \circ \sigma_{e,n-m})^{\frac{1}{p}} && \text{(Hölder's IE)} \\
        = & \quad (\cM_p)^{\frac{-(n-m)}{p}} \cE_p(f \circ \sigma_{e,n-m})^{\frac{1}{p}}. && \text{(Lemma \ref{lemma: Properties of cond.unif.}-\eqref{item: Properties of cond.unif. (1)} and \eqref{eq:duality})}
    \end{align*}
    By summing over all edges in $E_n$ and using Lemma \ref{lemma: Discrete energy (Self-similarity)}, we obtain the desired estimate
    \begin{align*}
        \cE_p(V_{n,m}[f]) & =
        \sum_{ e \in E_m  } \abs{V_{n,m}[f](e^-) - V_{n,m}[f](e^+)}^p \\
        & \leq \cM_p^{-(n-m)} \sum_{e \in E_m} \cE_p(f \circ \sigma_{e,n-m})\\
        & = \cM_p^{-(n-m)} \cE_p(f).
    \end{align*}
\end{proof}

\subsection{Poincar\'e inequality}\label{Subsec: Poincare}
The final addition to our toolbox is the following $(p,p)$-\emph{Poincar\'e-type inequality}. The proof is strongly inspired by \cite[Section 5]{CK_PI}.

\begin{theorem}[Poincar\'e inequality]\label{thm: Discrete PI}
    There is a constant $C > 0$ so that for all $n \in \N \cup \{0\}$ and $f : V_n \to \R$
    \begin{equation}\label{eq: Discrete PI}
        \frac{1}{\abs{E_n}}\sum_{v \in V_{n}} \abs{f(v) - V_{n,0}[f](v_{\pm})}^p \leq C\cM_p^{-n}\cE_p(f).
    \end{equation}
    Here, $v_{\pm}$ indicates that the statement is true for both $v_+$ and $v_-$.
\end{theorem}

 We first prove the following lemma.

\begin{lemma}\label{lemma: Differences of averages}
    There is a constant $C > 0$ so that the following holds.
    For every $n,i \in \N \cup \{0\}$ with $n > i$, and any $v_i \in e_i \in E_{i}$ and $v_{i+1} \in e_i \cdot G_1$, we have
    \begin{equation}
        \abs{V_{n,i+1}[f](v_{i+1}) - V_{n,i}[f](v_i)}^p \leq C \cM_p^{-(n-i)} \cE_p(f \circ \sigma_{e_i,n-i}).
    \end{equation}
\end{lemma}

\begin{proof}
    Note that for any $a \in I$ we have $v_i \cdot a \in e_i \cdot G_1$. By the self-similarity \ref{SM1} there is a path $[w_1,\dots w_k]$ from $v_{i+1}$ to $v_i \cdot a$ contained in $e_i \cdot G_1$ and of length at most $\diam(V_1,d_{G_1})$. Then we estimate using Hölder's inequality
    \begin{align*}
        \abs{V_{n,i+1}[f](v_{i+1}) - V_{n,i+1}[f](v_{i} \cdot a)}^p
        & \lesssim \sum_{l = 0}^{k-1}  \abs{V_{n,i+1}[f](w_{l+1})- V_{n,i+1}[f](w_{l})}^p\\
        & \leq  \mathcal{E}_p(V_{n,i+1}[f] \circ \sigma_{e_i,1}).
    \end{align*}
    We finish the proof by estimating
    \begin{align*}
        & \, \quad \abs{V_{n,i+1}[f](v_{i+1}) -V_{n,i}[f](v_i)}^p\\
        = & \quad \left| \sum_{a \in I} \mathcal{J}_p(a) (V_{n,i+1}[f](v_{i+1}) - V_{n,i+1}[f](v_{i} \cdot a)) \right|^p && \text{(Lemma \ref{lemma: Average compatible-discrete})} \\
        \leq & \quad \sum_{a \in I} \mathcal{J}_p(a) \abs{V_{n,i+1}[f](v_{i+1}) - V_{n,i+1}[f](v_{i} \cdot a)}^p && \text{(Jensen's IE)} \\
        \lesssim & \quad \sum_{a \in I} \mathcal{J}_p(a)\cE_p(V_{n,i+1}[f] \circ \sigma_{e_i,1} )\\
        = & \quad \cE_p(V_{n,i+1}[f] \circ \sigma_{e_i,1})\\
        = & \quad \cE_p(V_{n-i,1}[f \circ \sigma_{e_i,n-i}]) && \text{(Lemma \ref{lemma: Average commutes with similarity})} \\
        \lesssim & \quad   \cM_p^{-(n-i)}\cE_p(f \circ \sigma_{e_i,n-i}) && \text{(Theorem \ref{thm: Discrete energy monotonicity}).}
    \end{align*}
\end{proof}

\begin{proof}[Proof of Theorem \ref{thm: Discrete PI}]
    Since the IGS is non-degenerate, it follows from Lemma \ref{lemma: Properties of cond.unif.}-\eqref{item: Properties of cond.unif. (3)} that there is $\alpha > 1$ satisfying $\alpha \cdot\cM_p < \abs{E_1}$.

    First fix some $v := v_n \in V_n$. For each $i = 0,\dots,n-1$ we choose the vertices $v_0,\dots,v_{n-1}$ and edges $e_0,\dots,e_{n-1}$ so that $v_i \in e_{i-1} \cdot G_1$ and $v_i \in e_i \in E_i$. Then we estimate
    \begin{align*}
        \abs{f(v) - V_{n,0}[f](v_{\pm})}^p & \leq \left(\sum_{i = 0}^{n-1}  \abs{V_{n,i+1}[f](v_{i+1}) - V_{n,i}[f](v_i)}\alpha^{\frac{i}{p}} \alpha^{\frac{-i}{p}}\right)^p\\
        & \leq \left(\sum_{i = 0}^{n-1} \left(\alpha^{\frac{q}{p}}\right)^{-i} \right)^{\frac{p}{q}}\left(\sum_{i = 0}^{n-1} \abs{V_{n,i+1}[f](v_{i+1}) - V_{n,i}[f](v_i)}^p \alpha^{i} \right)\\
        & \leq \left(\sum_{i = 0}^{\infty} \left(\alpha^{\frac{q}{p}}\right)^{-i} \right)^{\frac{p}{q}} \sum_{i = 0}^{n-1} \abs{V_{n,i+1}[f](v_{i+1}) - V_{n,i}[f](v_i)}^p \alpha^{i}\\
        & \lesssim \sum_{i = 0}^{n-1} \cM_p^{-(n-i)}\cE_p(f \circ \sigma_{e_i,n-i}) \alpha^{i}.
    \end{align*} 
    We used triangle inequality in the first row, Hölder's inequality in the second and the last row follows from Lemma \ref{lemma: Differences of averages} and the fact that $\alpha > 1$.
    
    Next we do the previous choices of vertices and edges for all $v \in V_n$, which we denote $v_i(v) \in e_i(v) \in E_i$. Notice that $v \in e_i(v) \cdot G_{n-i}$. Since $\abs{V_{n-i}} \leq 2\abs{E_1}^{n-i}$, we conclude the following:
    \begin{equation}\label{Theorem: Poincare number of choices}
        \text{Each edge in $E_i$ is chosen at most } \abs{E_1}^{n-i} \text{ many times.} \tag{$\star$}
    \end{equation}
    We finish the proof by estimating
    \begin{align*}
        \sum_{v \in V_n} \abs{f(v) - V_{n,0}[f](v_-)}^p & \lesssim \sum_{v \in V_n} \sum_{i = 0}^{n-1} \cM_p^{-(n-i)}\cE_p(f \circ \sigma_{e_i(v),n-i}) \alpha^{i}\\
        & \stackrel{\eqref{Theorem: Poincare number of choices}}{\leq} \sum_{i = 0}^{n-1} \left( \cM_p^{-(n-i)}\alpha^{i} \cdot \abs{E_1}^{n-i}\sum_{e \in E_{i}}\cE_p(f \circ \sigma_{e,n-i})\right)\\
        & = \cE_p(f) \cdot \sum_{i = 0}^{n-1} \cM_p^{-(n-i)}\alpha^{i} \cdot \abs{E_1}^{n-i}\\
        & = \abs{E_1}^n\cM_p^{-n}\cE_p(f)\sum_{i = 1}^{n-1} \left(\alpha \frac{\cM_p}{\abs{E_1}}\right)^i \\
        & \leq \abs{E_1}^n\cM_p^{-n}\cE_p(f)\sum_{i = 1}^{\infty} \left(\alpha \frac{\cM_p}{\abs{E_1}}\right)^i \\
        & \lesssim \abs{E_1}^n\cM_p^{-n}\cE_p(f). 
    \end{align*}
    The equality in the third line follows from Lemma \ref{lemma: Discrete energy (Self-similarity)}, and the geometric series in the second last line converges due to the choice of the constant $\alpha$.
\end{proof}

\section{Discretizations and Mollifiers}\label{Sec: Discretizations and mollifiers}
We are now ready to take the first step towards analysis. Using the tools of discrete potential theory developed in Section \ref{sec: DPT}, we introduce analytic objects that will support us in later sections.
Unless indicated otherwise, we shall throughout the paper consider a fixed exponent $p \in (1,\infty)$ and an IGS satisfying Assumption \ref{Assumptions: To make work}. 
For its limit space $(X,d,\mu)$, we use the same notation as in Definition \ref{def: Limit space}.

\subsection{Discretization operators}
The preliminary step in our framework is discretization. We introduce linear operators $V_n[\cdot] : C(X) \to \R^{V_n}$ that discretizes continuous functions.
To properly define them, we need to recall the definition of the fiber measures in Definition \ref{def: Fiber measure}, and the probability density $\mathcal{J}_p(\cdot)$ on $I$ from Remark \ref{Remark: Measure on gluing set}.

\begin{definition}\label{def: Discretization operator}
    Let $n \in \N \cup\{0\}$ and $p \in (1,\infty)$. We define the \emph{discretization operator} $V_{p,n}[\cdot] : C(X) \to \R^{V_n}$ as the linear operator which maps any continuous function $f \in C(X)$ to the function $V_{p,n}[f] : V_n \to \R$ given by
    \begin{equation*}
        V_{p,n}[f](v) := \int_{\Fib(v)} f \, d\mathcal{J}_{p,v}.
    \end{equation*}
    We then define the \emph{discrete $p$-energy} of $f \in C(X)$ as the normalized discrete $p$-energy
    \begin{equation*}
        \cE_p^{(n)}(f) := \cM_p^{-n} \cE_p( V_{p,n}[f]).
    \end{equation*}
\end{definition}

Since the underlying fiber measure depends on $p$, so do the above discretization operators. For simplicity we omit the exponent from the notation and write $V_n[\cdot]$ instead of $V_{p,n}[\cdot]$.

The operators $V_n[\cdot]$ satisfy the following \emph{tower rule}, which is essentially thanks to the self-similarity of the fiber measures.

\begin{lemma}[Tower rule]\label{lemma: Tower rule}
    Let $n,m\in \N \cup \{0\}$ so that $n \geq m$.
    For all $f \in C(X)$,
    \begin{equation}\label{eq: Tower rule}
        V_{n,m}[[V_{n}[f]] = V_m[f].
    \end{equation}
    Here $V_{n,m}[\cdot]$ is the averaging operator in \eqref{eq: Discrete average}.
\end{lemma}

\begin{proof}
    Let $v \in V_m$. Using \eqref{eq: Fiber measures (Homogeneity)} and Lemma \ref{lemma: How to measure on fibers}-\eqref{item: How to measure on fibers (1)}, we compute
    \begin{align*}
        V_{n,m}[V_n[f]](v) & = \sum_{a \in I^{n-m}} \mathcal{J}_{p,n-m}(a) V_n[f](v \cdot a)\\
        & = \sum_{a \in I^{n-m}} \mathcal{J}_{p,n-m}(a)\int_{\Fib(v\cdot a)}f \, d\mathcal{J}_{p,v \cdot a}\\
        & = \sum_{a \in I^{n-m}} \int_{\Fib(v\cdot a)} f \, d \mathcal{J}_{p,v}\\
        & = \int_{\Fib(v)}f \, d\mathcal{J}_{p,v} \\
        & = V_m[f](v).
    \end{align*}
\end{proof}

\begin{lemma}\label{lemma: Similarity maps and average}
    Let $f \in C(X)$ and $n,k \in \N \cup \{0\}.$ Then for all edges $e \in E_k$ it holds that
    \begin{equation}\label{eq: Similarity maps and average}
        V_{n}[f \circ F_e] = V_{n+k}[f] \circ \sigma_{e,n}.
    \end{equation}
\end{lemma}

\begin{proof}
    By the definition of fiber measures \eqref{eq: Fiber measure} and Lemmas \ref{lemma: Symbolic associativity}-\eqref{item: Symbolic associativity (3)} and \ref{lemma: How to measure on fibers}-\eqref{item: How to measure on fibers (2)}, the fiber measure $\mathcal{J}_{p,e \cdot v}$ is equal to the push-forward measure $\mathcal{J}_{p,e \cdot v} = (F_e)_*(\mathcal{J}_{p,v})$. Given this, we compute
    \begin{align*}
        V_{n}[f \circ F_e](v) & = \int_{\Fib(v)} f \circ F_e \, d\mathcal{J}_{p,v} = \int_{\Fib(e \cdot v)} f \, d(F_e)_*(\mathcal{J}_{p,v})\\
        & = \int_{\Fib(e \cdot v)} f \, d\mathcal{J}_{p,e \cdot v} = V_{n+k}[f](\sigma_{e,n}(v)).
    \end{align*}
\end{proof}

\subsection{Optimal potential functions}\label{subsec: OPF}
Next, we discuss the most important family of functions in our analysis. For $p \in (1,\infty)$ we denote
$\mathscr{U}_{p,+}: X \to \R$ the unique continuous function satisfying
\begin{equation}\label{eq: Smooth function on fibers}
    \mathscr{U}_{p,+}|_{\Fib(v)} \equiv U_{p,+,n}(v) \text{ for all } n \in \N \cup \{0\} \text{ and } v \in V_n.
\end{equation}
Similarly, we define $\mathscr{U}_{p,-}$ by replacing the sign $+$ by $-$ in \eqref{eq: Smooth function on fibers}. Frequently, when the computation at the moment works for both signs, we write $\mathscr{U}_{p,\pm}$.
These functions are collectively referred as the \emph{(continuous) optimal potential functions}.

\begin{proposition}\label{prop: Smooth function on fibers}
The condition \eqref{eq: Smooth function on fibers} uniquely determines a continuous function $\mathscr{U}_{p,\pm}: X \to \R$, and it satisfies the following properties.
\begin{enumerate}
    \vskip.3cm
    \item \label{item: Smooth function on fibers (1)} For all $v \in e \in E_\#$ the function
    $\mathscr{U}_{p,\pm}$ is $\delta_p$-Hölder continuous for
    \[
        \delta_p := -\frac{\log(c_p)}{\log(L_*)},
    \]
    where $c_p$ is as in Corollary \ref{Cor: Discrete hölder continuity}.
    \vskip.3cm
    \item \label{item: Smooth function on fibers (2)}
    It holds that
    \[
        \mathscr{U}_{p,\pm}|_{\Fib(v_{\pm})} = 1 \text{ and } \mathscr{U}_{p,\pm}|_{\Fib(v_{\mp})} = 0.
    \]
    \vskip.3cm
    \item \label{item: Smooth function on fibers (3)}
    It holds that
    \[
        0 \leq \mathscr{U}_{p,\pm} \leq 1 \text{ and } \mathscr{U}_{p,+} = 1 - \mathscr{U}_{p,-}  .
    \]
    \vskip.3cm
    \item \label{item: Smooth function on fibers (4)}
    For all $n \in \N \cup \{0\}$ it holds that $\cE_p^{(n)}(\mathscr{U}_{p,\pm}) = 1$.
\end{enumerate}
\end{proposition}

\begin{proof}
    It follows from Lemma \ref{lemma: How to measure on fibers}-\eqref{item: How to measure on fibers (3)} and \eqref{eq: Optimal potential on Fibers} that the condition \eqref{eq: Smooth function on fibers} determines a well-defined function on the fiber set $\mathscr{U}_{p,\pm} : \Fib(X) \to \R$.
    Since $\Fib(X) \subseteq X$ is a dense subset according to Corollary \ref{cor: Fib(X) is dense}, the existence and uniqueness of the continuous extension of $\mathscr{U}_{p,\pm}$ to $X$, and consequently \eqref{item: Smooth function on fibers (1)} of the current lemma, follow once we prove that $\mathscr{U}_{p,\pm}$ is $\delta_p$-Hölder continuous in $\Fib(X)$.

    \eqref{item: Smooth function on fibers (1)}:
    Let $x,y \in \Fib(X)$ be arbitrary distinct points. By Lemma \ref{lemma: How to measure on fibers}-\eqref{item: How to measure on fibers (1)} there is $m \in \N \cup \{0\}$ and vertices $v,w \in V_m$ so that $x \in \Fib(v)$ and $y \in \Fib(w)$. Then we choose any edges $v \in e_v$ and $w \in e_w$. By the definition of the fiber sets, we can express $x = \chi((e_i)_{i = 1}^{\infty})$ and $y = \chi((e_i')_{i = 1}^{\infty})$ where $e_m = e_v$ and $e_m' = e_w$. Lastly, we fix $n \in \N \cup \{0\}$ to be the largest non-negative integer for which $e_n \cap e_n' \neq \emptyset$. By the visual metric property (Proposition \ref{prop: LS well-defined}-\eqref{item: Visual metric}) we have
    \begin{equation}\label{eq: For Hölder}
        d(x,y) \gtrsim L_*^{-n}.
    \end{equation}
    Also, thanks to Lemma \ref{lemma: How to measure on fibers}-\eqref{item: How to measure on fibers (1)}, we may assume that $m \geq n$.

    Since $v \in e_v = e_m$ and $\pi_{m,n}(e_v) = e_n$, it follows from \ref{SM3} that $v \in e_n \cdot G_{m-n}$. Similarly, $w \in e_n' \cdot G_{m-n}$. Using \ref{SM2} we can choose a vertex $v_m \in e_n \cdot G_{m-n} \cap e_n' \cdot G_{m-n}$. We now use Corollary \ref{Cor: Discrete hölder continuity} and \eqref{eq: For Hölder} to estimate
    \begin{align*}
        \abs{\mathscr{U}_{p,\pm}(x) - \mathscr{U}_{p,\pm}(y)} & = \abs{U_{p,\pm,m}(v) - U_{p,\pm,m}(w)}\\
        & \leq \abs{U_{p,\pm,m}(v) - U_{p,\pm,m}(v_m)} + \abs{U_{p,\pm,m}(v_m) - U_{p,\pm,m}(w)}\\
        & \leq 2c_p^n = 2(L_*^{-n})^{\delta_p} \lesssim d(x,y)^{\delta_p}.
    \end{align*}

    \eqref{item: Smooth function on fibers (2)}:
    By using the definition \eqref{eq: Smooth function on fibers}, we have
    \[
        \mathscr{U}_{p,\pm}|_{\Fib(v_{\pm})} = U_{p,0}(v_\pm) = 1 \text{ and } \mathscr{U}_{p,\pm}|_{\Fib(v_{\mp})} = U_{p,0}(v_\mp) = 0.
    \]
    
    \eqref{item: Smooth function on fibers (3)}:
    Note that, by Lemma \ref{lemma: Properties of cond.unif.}-\eqref{item: Properties of cond.unif. (4)}, we have $0 \leq U_{p,\pm,m} \leq 1$ for all $m \in \N \cup \{0\}$. Also $U_{p,+,m} = 1 - U_{p,-,m}$ follows from the uniqueness of the optimal potential.
    By the density $\Fib(X) \subseteq X$, these properties directly translate to the analogous properties of $\mathscr{U}_{p,\pm}$.

    \eqref{item: Smooth function on fibers (4)}: Since, by the definition \eqref{eq: Smooth function on fibers}, $\mathscr{U}_{p,\pm}$ is constant on all fibers.
    Thus, we have $V_n[\mathscr{U}_{p,\pm}] = U_{p,\pm,n}$ for all $n \in \N \cup \{0\}$. By using Lemma \ref{lemma: Properties of cond.unif.}-\eqref{item: Properties of cond.unif. (1)}, we have
    \[
        \mathcal{E}_p^{(n)}(\mathscr{U}_{p,\pm}) = \cM_p^{-n}\mathcal{E}_p(U_{p,\pm,n}) = 1.
    \]
\end{proof}

\subsection{Mollifiers}
Next, we introduce the mollifier operators. As the name suggests, their main purpose is to produce ``smoother'' approximations. In this subsection, the mollifiers are defined only on continuous functions, but later in Subsection \ref{subsec: Extending operators} we extend them for general Sobolev functions. We prove in Theorem \ref{thm: Mollifiers converge} that the mollifiers enjoy convergence properties in the Sobolev norm.
    
\begin{definition}\label{Def: Definition of U_p}
    For $n \in \N \cup \{0\}$ we define the \emph{interpolation operator} $\mathscr{U}_{p,n}[\cdot] : \R^{V_n} \to C(X)$ to be the linear operator so that for any given function $g : V_n \to \R$, $\mathscr{U}_{p,n}[g] : X \to \R $ is the continuous function satisfying 
    \begin{equation}\label{eq: Definition of U_p}
    (\mathscr{U}_{p,n}[g] \circ F_e)(x) := g(e^-) + (g(e^+) - g(e^-))\mathscr{U}_{p,+}(x)
\end{equation}
    for all $e \in E_n$ and $x \in X$. The \emph{mollifier} is then defined as the composition $\Psi_{p,n}[\cdot] := \mathscr{U}_{p,n}[\cdot] \circ V_n[\cdot] : C(X) \to C(X)$.
\end{definition}

\begin{lemma}\label{lemma: Properties of U_p}
    For all $n \in \N \cup \{0\}$ and $g : V_n \to \R$ the condition \eqref{eq: Definition of U_p} determines a well-defined $\delta_p$-Hölder continuous function $\mathscr{U}_{p,n}[g] \in C(X)$, where $\delta_p$ is as in Proposition \ref{prop: Smooth function on fibers}, satisfying the following conditions.
    \vskip.3cm
    \begin{enumerate}
        \item \label{item: Properties of U_p (1)}
    $\cE_p^{(m)}(\mathscr{U}_{p,n}[g]) = \cM_p^{-n}\cE_p(g)$ for all $m \geq n$.
        \vskip.3cm
        \item \label{item: Properties of U_p (2)} For all $n \in \N \cup \{0\}$ and $e \in E_n$ it holds that
        \begin{align*}
            \sup_{x,y \in X_e} & \abs{\mathscr{U}_{p,n}[g](x) - \mathscr{U}_{p,n}[g](y)} = \abs{\nabla g(e)},\\
        \sup_{x \in X_e} & \abs{\mathscr{U}_{p,n}[g](x)} = \max\{ \abs{g(e^-)},\abs{g(e^+)} \}. 
        \end{align*}
    \end{enumerate}
\end{lemma}

\begin{proof}
    First, we verify that $\mathscr{U}_{p,n}[g]$ is well-defined. Since $e \cdot v_\pm = e^{\pm}$ it follows from Lemma \ref{lemma: How to measure on fibers}-\eqref{item: How to measure on fibers (2)}, Proposition \eqref{prop: Smooth function on fibers}-\eqref{item: Smooth function on fibers (2)} and the definition of the interpolation \eqref{eq: Definition of U_p} that
    \[
        \mathscr{U}_{p,n}[g]|_{\Fib(e^\pm)} = g(e^\pm) \text{ for all } e \in E_n.
    \]
    Hence, $\mathscr{U}_{p,n}[g]$ is well-defined by Proposition \ref{prop: X_e cap X_f}. Furthermore, it is $\delta_p$-Hölder continuous by Proposition \ref{prop: Smooth function on fibers}-\eqref{item: Smooth function on fibers (1)}.
    
    \eqref{item: Properties of U_p (1)}: Let $n,m \in N \cup \{0\}$ so that $m \geq n$. We compute
    \begin{align*}
        & \, \quad \mathcal{E}_p^{(m)}(\mathscr{U}_{p,n}[g])\\
        = & \quad \cM_p^{-m} \mathcal{E}_p(V_m[\mathscr{U}_{p,n}[g]])\\
        = & \quad \cM_p^{-m}\sum_{e \in E_n} \mathcal{E}_p(V_m[\mathscr{U}_{p,n}[g]] \circ \sigma_{e,m-n}) && \text{(Lemma \ref{lemma: Discrete energy (Self-similarity)})} \\
        = & \quad \cM_p^{-m}\sum_{e \in E_n} \mathcal{E}_p(V_{m-n}[\mathscr{U}_{p,n}[g] \circ F_e]) && \text{(Lemma \ref{lemma: Similarity maps and average})}\\
        = & \quad  \cM_p^{-m}\sum_{e \in E_n} \abs{\nabla g(e)}^p \mathcal{E}_p(V_{m-n}[\mathscr{U}_{p,+}]) && \text{(Definition \ref{Def: Definition of U_p})} \\
        = & \quad \cM_p^{-n} \mathcal{E}_p(g) && \text{(Proposition \ref{prop: Smooth function on fibers})-\eqref{item: Smooth function on fibers (4)}).}
    \end{align*}

    \eqref{item: Properties of U_p (2)}:
    Let $e \in E_n$ and $y \in X_e$. By \eqref{item: Smooth function on fibers (2)} and \eqref{item: Smooth function on fibers (3)} of Lemma \ref{prop: Smooth function on fibers}, $\mathscr{U}_{p,n}[g](y)$ is a convex combination of $g(e^+)$ and $g(e^-)$, and $\mathscr{U}_{p,n}[g](y) \in \{ g(e^-), g(e^+)\}$ if $y \in \Fib(e^{\pm})$. Thus we have $\sup_{x \in X_{e}}\abs{\mathscr{U}_{p,n}[g](x)} = \max\{ \abs{g(e^-)}, \abs{g(e^+)} \}$,
    \[
        \abs{\mathscr{U}_{p,n}[g](x) - \mathscr{U}_{p,n}[g](y)} \leq \abs{\nabla g(e)} \text{ for all } x,y \in X_e,
    \]
    and this inequality is equality for $x \in \Fib(e^-)$ and $y \in \Fib(e^+)$.
\end{proof}

\begin{proposition}\label{prop: Properties of mollifier}
    Let $f \in C(X)$ and $n \in \N \cup \{0\}$. Then the following holds.
    \begin{enumerate}
        \item \label{item: Properties of mollifier (1)}
        $\cE_p^{(m)}(\Psi_{p,n}[f]) = \cE_p^{(n)}(f)$ for all $m \geq n$. 
        \vskip.3cm
        \item \label{item: Properties of mollifier (2)} For all $m \geq n$ and $v \in V_m$ the function $\Psi_{p,n}[f]$ is constant on $\Fib(v)$.
        \vskip.3cm
        \item \label{item: Properties of mollifier (3)}
        For all $e \in E_n$ we have
        \begin{equation*}
            \sup_{x,y \in X_e} \abs{\Psi_{p,n}[f](x)-\Psi_{p,n}[f](y)} = \abs{\nabla V_n[f](e)}
        \end{equation*}
        and
        \begin{equation*}
        \sup_{x \in X_e}
        \abs{\Psi_{p,n}[f](x)} = \max\{ \abs{V_n[f](e^-)}, \abs{V_n[f](e^+)} \}.
        \end{equation*}
        \vskip.3cm
        \item \label{item: Properties of mollifier (4)}
        $\Psi_{p,n}[f] \to f$ in $L^{\infty}(X)$ as $n \to \infty$.
    \end{enumerate}
\end{proposition}

\begin{proof}
    \eqref{item: Properties of mollifier (1)}: Directly follows from Lemma \ref{lemma: Properties of U_p}-\eqref{item: Properties of U_p (1)}.

    \eqref{item: Properties of mollifier (2)}:
    Let $m \geq n$ and $v \in V_m$. Choose an edge $e \in E_n$ and $v' \in V_{m-n}$ so that $v = e \cdot v'$. Then we have
    \[
        \Psi_{p,n}[f]|_{\Fib(v)} = V_{n}[f](e^-) + (V_{n}[f](e^+) - V_{n}[f](e^-))\mathscr{U}_{p,+}|_{\Fib(v')}.
    \]
    It is clear from \eqref{eq: Smooth function on fibers} that $\mathscr{U}_{p,+}$ is constant on $\Fib(v')$, so we are done.

    \eqref{item: Properties of mollifier (3)}: Directly follows from Lemma \ref{lemma: Properties of U_p}-\eqref{item: Properties of U_p (2)}.

    \eqref{item: Properties of mollifier (4)}: If $v \in e \in E_n$, it follows from Proposition \ref{prop: Geom of LS}-\eqref{item: diams/measures of cells} that
    \[
        \diam(\Fib(v)) \leq \diam(X_e) = L_*^{-n}.
    \]
    It is now routine to verify that $\Psi_{p,n}[f] \to f$ in $L^{\infty}(X)$ as $n \to \infty$ using the uniform continuity of $f$ and \eqref{item: Properties of mollifier (3)} of the current proposition.
\end{proof}

\section{Construction of the energy forms}\label{Sec: Sobolev spaces}
We now have enough ingredients to define and study the self-similar $p$-energy forms $\mathscr{E}_p : L^p(X,\mu) \to [0,\infty]$.
We state the general results in Subsection \ref{subsec: Main results}, and prove them in Subsections \ref{subsec: Closability} and \ref{subsec: contraction}. In Subsection \ref{subsec: Extending operators}, we extend the discretizations and mollifiers introduced in Section \ref{Sec: Discretizations and mollifiers} to general Sobolev functions. Subsection \ref{subsec:EM} studies the natural counterparts of the measures $A \mapsto \int_A \abs{\nabla f}^pdx$.
Lastly, let us recall that Assumption \ref{Assumptions: To make work} holds.

\subsection{Main theorems}\label{subsec: Main results}
We first consider the \emph{pre-energy form} $\mathsf{E}_p : C(X) \to [0,\infty]$,
\begin{equation}
    \mathsf{E}_p(f) := \lim_{n \to \infty} \cM_p^{-n} \mathcal{E}_p(V_n[f]) = \lim_{n \to \infty} \mathcal{E}_p^{(n)}(f),
\end{equation}
where $\cM_p$ is the $p$-capacity constant in Assumption \ref{Assumptions: To make work} and $V_n[\cdot]$ are the discretization operators in Definition \ref{def: Discretization operator}.
We show in Theorem \ref{thm: Strong monotonicity} that $\mathsf{E}_p(f)$ always exists in $[0,\infty]$.
We also define the \emph{core} as the collection
\begin{equation}\label{eq:CORE}
    \mathscr{C}_p := \left\{ f \in C(X) : \mathsf{E}_p(f) < \infty \right\}.
\end{equation}
Lastly, we say that a sequence $\{ f_i \}_{i \in \N} \subseteq \mathscr{C}_p$ is \emph{$\mathsf{E}_p$-Cauchy} if 
\begin{equation}\label{eq: Cauchy seq}
    \lim_{i,j \to \infty} \mathsf{E}_p(f_i - f_j) = 0.
\end{equation}

Now, our initial target is to establish the following \emph{closability} of $\mathsf{E}_p$.

\begin{theorem}\label{thm: Closability}
    The pre-energy form $\mathsf{E}_p : \mathscr{C}_p \to [0,\infty)$ is closable, i.e., for every $\mathsf{E}_p$-Cauchy sequence $\{ f_i \}_{i = 1}^{\infty} \subseteq \mathscr{C}_p$ so that $f_i \to 0$ in $L^p(X,\mu)$ we have $\mathsf{E}_p(f_i) \to 0$.
\end{theorem}

Having the closability, $\mathsf{E}_p$ has the following canonical extension to $L^p(X,\mu)$.

\begin{definition}\label{def: Sobolev}
We define the \emph{self-similar $p$-energy form} $\mathscr{E}_p : L^p(X,\mu) \to [0,\infty]$ as follows.
For an arbitrary $f \in L^p(X,\mu)$ we set
\begin{equation*}\label{def: E_p for Lp}
    \mathscr{E}_p(f) := \lim_{i \to \infty} \mathsf{E}_p(f_i)
\end{equation*}
where $\{f_i\}_{i = 1}^{\infty} \subseteq \mathscr{C}_p$ is any $\mathsf{E}_p$-Cauchy sequence so that $f_i \to f$ in $L^p(X,\mu)$. 
Note that it follows from Minkowski's inequality that $\mathscr{E}_p(f)$ does not depend on the choice of the $\mathsf{E}_p$-Cauchy sequence as soon as we have verified closability of $\mathsf{E}_p$.
If such sequence does not exist, we define $\mathscr{E}_p(f):= \infty$. Finally, we define the associated \emph{Sobolev space},
\[
    \mathscr{F}_p := \{ f \in L^p(X,\mu) : \mathscr{E}_p(f) < \infty \}
\]
and equip it with the \emph{Sobolev norm} $\norm{\cdot}_{\mathscr{F}_p} := \norm{\cdot}_{L^p} + \mathscr{E}_p(\cdot)^{\frac{1}{p}}$.
\end{definition}

We deal with the closability in Subsection \ref{subsec: Closability}. For now, assume that we already have it, so that we can state the important analytic properties of $\mathscr{E}_p$.

\begin{theorem}\label{Thm: Existence of p-Energy}
    The $p$-energy form $\mathscr{E}_p$ and the associated Sobolev space $(\mathscr{F}_p,\norm{\cdot}_{\mathscr{F}_p})$ satisfy the following properties.
    \begin{enumerate}
        \vskip.3cm
        \item \label{item: Energy is explicit}
        $\mathscr{E}_p(f) = \mathsf{E}_p(f)$ for all $f \in \mathscr{C}_p$.
        \vskip.3cm
        \item \label{item: Refl-Sep}
        The Sobolev space $(\mathscr{F}_p,\norm{\cdot}_{\mathscr{F}_p})$ is a reflexive and separable Banach space.
        \vskip.3cm
        \item \label{item: Regularity}
        \textup{\textbf{(Regularity)}} 
        $\mathscr{C}_p$ is dense in both $(C(X), \norm{\cdot}_{L^{\infty}})$ and $(\mathscr{F}_p,\norm{\cdot}_{\mathscr{F}_p})$.
        \vskip.3cm
        \item \label{item: p-energySS}
        \textup{\textbf{(Self-similarity)}} 
        For all $f \in \mathscr{F}_p$,
        \begin{equation*}
            \mathscr{E}_p(f) = \cM_p^{-1} \sum_{e \in E_1}\mathscr{E}_p(f \circ F_e).
        \end{equation*}
        \vskip.3cm
        \item \label{item: Contraction}
        \textup{\textbf{(Lipschitz contractivity)}}  If $\varphi : \R \to \R$ is 1-Lipschitz and $f \in \mathscr{F}_p$ then $\mathscr{E}_p(\varphi \circ f) \leq \mathscr{E}_p(f)$.
    \end{enumerate}
\end{theorem}

\begin{remark}
    The equality in Theorem \ref{Thm: Existence of p-Energy}-\eqref{item: Energy is explicit} is extended from $\mathscr{C}_p$ to $C(X)$ in Corollary \ref{cor: Self-similarity of domain}. This turns out to be a surprisingly delicate detail.
\end{remark}

\begin{theorem}\label{thm: Energy analytic properties}
    We also have the following properties.
    \begin{enumerate}
        \item \label{item: Poincare}
        \textup{\textbf{(Poincar\'e inequality)}}  There is a constant $C > 0$ so that for all $n \in \N \cup \{0\}$, $e \in E_n$ and  $f \in \mathscr{F}_p$ we have the $(p,p)$-Poincar\'e-type inequality
        \begin{align}
            \label{eq: PrePI (1)} \kint_{X_e} \abs{f - f_{X_e}}^p \, d\mu & \leq C \mathscr{E}_p(f \circ F_e).
        \end{align}
        In particular, $\mathscr{E}_p(f) = 0$ if and only if $f$ is constant $\mu$-almost everywhere.
        \vskip.3cm
        \item \label{item: p-Clarkson}
        \textup{\textbf{($p$-Clarkson's inequality)}}  For all $f,g \in \mathscr{F}_p$ it holds that\begin{equation}\label{e:Cp} \begin{cases}\mathscr{E}_{p}(f+g) + \mathscr{E}_{p}(f-g) \ge 2\bigl(\mathscr{E}_{p}(f)^{\frac{1}{p-1}} + \mathscr{E}_{p}(g)^{\frac{1}{p-1}}\bigr)^{p - 1} \quad &\text{if $p \in (1,2]$,} \\ \mathscr{E}_{p}(f+g) + \mathscr{E}_{p}(f-g) \le 2\bigl(\mathscr{E}_{p}(f)^{\frac{1}{p-1}} + \mathscr{E}_{p}(g)^{\frac{1}{p-1}}\bigr)^{p - 1} \quad &\text{if $p \in [2,\infty)$.}\end{cases}\end{equation}
        \vskip.3cm
        \item \label{item: strong-locality}
        \textup{\textbf{(Strong locality)}}  If $f,g \in \mathscr{F}_p$ so that $\supp_{\mu}(f - a) \cap \supp_{\mu}(g - b) = \emptyset$ for some $a,b \in \R$, then for all $h \in \mathscr{F}_p$ we have
        \begin{equation*}
        \mathscr{E}_p(f + g + h) + \mathscr{E}_p(h) = \mathscr{E}_p(f + h) + \mathscr{E}_p(g + h).\end{equation*}
    \end{enumerate}
\end{theorem}

\subsection{Closability
}\label{subsec: Closability}
In order to verify the closability of $\mathsf{E}_p$ we first show that $\mathsf{E}_p$ is well-defined in the first place.

\begin{theorem}\label{thm: Strong monotonicity}
    For all $f \in C(X)$ the sequence of $p$-energies $\cE_p^{(n)}(f)$ as defined in Definition \ref{def: Discretization operator} is non-decreasing with respect to $n \in \N \cup \{0\}$. In particular, $\mathsf{E}_p(f) \in [0,\infty]$ is well-defined.
\end{theorem}

\begin{proof}
    Let $n,m \in \N \cup \{0\}$ so that $m \geq n$. The monotonicity now follow as
    \begin{align*}
        \cE_p^{(n)}(f) & = \cM_p^{-n} \cdot \cE_p(V_n[f])\\
        & = \cM_p^{-n} \cdot \cE_p(V_{m,n}[V_m[f]]) && \text{(Tower rule)} \\
        & \leq \cM_p^{-n} \cdot \cM_p^{-(m-n)} \cE_p(V_m[f]) && \text{(Strong monotonicity)} \\
        & = \cM_p^{-m} \cdot \cE_p(V_m[f])\\
        & = \cE_p^{(m)}(f).
    \end{align*}
\end{proof}

The main ingredients to verify the closability are the following two lemmas. The first one regards the self-similarity of $\mathsf{E}_p$.

\begin{lemma}\label{lemma: Energy self-similarity}
   For all $f \in C(X)$ and $m \in \N$
   \begin{equation}\label{eq: Energy self-similarity}
       \mathsf{E}_p(f) = \cM_p^{-m}  \sum_{e \in E_m} \mathsf{E}_p(f \circ F_e).
   \end{equation}
\end{lemma}

\begin{proof}
    By applying Lemmas \ref{lemma: Discrete energy (Self-similarity)} and \ref{lemma: Similarity maps and average} we simply compute
    \begin{align*}
        \mathsf{E}_p(f) 
        & = \lim_{n \to \infty}  \cM_p^{-(n+m)} \cE_p(V_{n+m}[f])\\
        & = \cM_p^{-m} \sum_{e \in E_m} \lim_{n \to \infty} \cM_p^{-n} \cE_p(V_{n+m}[f] \circ \sigma_{e,n})\\
        & = \cM_p^{-m} \ \sum_{e \in E_m} \lim_{n \to \infty} \cM_p^{-n} \cE_p(V_{n}[f \circ F_e])\\
        & = \cM_p^{-m}  \sum_{e \in E_m} \mathsf{E}_p(f \circ F_e).
    \end{align*}
\end{proof}

The second ingredient is the $(p,p)$-Poincar\'e-type inequality.
\begin{lemma}\label{lemma: Poincare for C(X)}
    There is a constant $C > 0$ so that for all continuous functions $f \in C(X), m \in \N \cup \{0\}$ and $e \in E_{m}$, it holds that
    \begin{equation}\label{eq: Poincare for C(X)}
        \kint_{X_e} \abs{f - V_m[f](e^{\pm})}^p \, d\mu \leq C\mathsf{E}_p(f \circ F_e).
    \end{equation}
\end{lemma}

\begin{proof}
    We will only prove the case $m = 0$ and $v = v_{\pm}$. To justify this, it follows from the self-similarity properties (Proposition \ref{prop: Geom of LS}-\eqref{item: diams/measures of cells} and Lemma \ref{lemma: Similarity maps and average}) that
    \[
        \kint_{X_e} \abs{f - V_m[f](e^{\pm})}^p \,d\mu = \int_{X} \abs{f \circ F_e - V_0[f \circ F_e](v_{\pm})}^p \, d\mu.
    \]
    
    Fix $\varepsilon > 0$.
    By Proposition \ref{prop: Properties of mollifier}-\eqref{item: Properties of mollifier (4)} we can choose a large enough $n \in \N$ so that
    \[
        \int_X \abs{f}^p d\mu - \varepsilon \leq \int_X \abs{\Psi_{p,n}[f]}^p \, d\mu.
    \]
    Using this and Proposition \ref{prop: Properties of mollifier}-\eqref{item: Properties of mollifier (3)}, we estimate
    \begin{align*}
        \int_{X} \abs{f}^p\, d\mu - \varepsilon & \leq \int_{X} \abs{\Psi_{p,n}[f]}^p \, d\mu \\
        & = \sum_{e \in E_n} \int_{X_{e}} \abs{\Psi_{p,n}[f]}^{p} \, d\mu\\
        & \leq \sum_{\substack{e \in E_n}} \int_{X_e} \abs{V_n[f](e^-)}^p + \abs{V_n[f](e^+)}^p \, d\mu \\
        & = \sum_{v \in V_n} \sum_{v \in e \in E_n} \mu(X_e)\abs{V_n[f](v)}^p\\
        & \leq \frac{C_{\deg}}{\abs{E_1}^n} \sum_{v \in V_n}\abs{V_n[f](v)}^p.
    \end{align*}
    In the last line we used $\mu(X_e)=\abs{E_1}^{-n}$, which is proven in Proposition \ref{prop: Geom of LS}-\eqref{item: diams/measures of cells}.
    Next, we replace $f$ with $f - V_0[f](v_\pm)$ in the previous computation. Since 
    \[
        V_n[f - V_0[f](v_\pm)] = V_n[f] - V_0[f](v_\pm)    
    \]
    we obtain the estimate
    \begin{align*}
        & \, \quad\int_{X} \abs{f - V_0[f](v_\pm)}^p\, d\mu - \varepsilon\\
        \lesssim & \quad \frac{1}{\abs{E_1}^n} \sum_{v \in V_n}\abs{V_n[f](v) - V_0[f](v_\pm)}^p\\
        = & \quad \frac{1}{\abs{E_1}^n}\sum_{v \in V_n}\abs{V_n[f](v) - V_{n,0}[V_n[f]](v_\pm)}^p && \text{(Tower rule)} \\
        \lesssim & \quad \cM_p^{-n}\mathcal{E}_p(V_n[f]) && \text{(Poincar\'e IE \eqref{eq: Discrete PI})} \\
        \leq & \quad \mathsf{E}_p(f). && \text{(Theorem \ref{thm: Discrete energy monotonicity})}
    \end{align*}
    The result now follows by letting $\varepsilon \to 0$.
\end{proof}

Lastly we derive the following careful estimate.

\begin{proposition}\label{prop: Equivalent change}
    There is a constant $C_0 > 0$ so that for all $f \in C(X)$, $n,m \in \N$ and $v \in V_n$ we have
    \begin{equation}\label{eq: Equivalent change}
        \abs{V_n[f](v)}^p \leq C_0 \left(\frac{\cM_p^n}{2^{m(p-1)}} \mathsf{E}_p(f) + \frac{\abs{I}^m\abs{E_1}^{n+m}}{\cM_p^{m}2^{m(p-1)}}\norm{f}_{L^p}^p\right).
    \end{equation}
\end{proposition}

\begin{proof}
    Fix $k \in \N \cup \{0\}$ and $w \in e \in E_k$. We first derive a preliminary inequality,
    \begin{align}
        \label{eq: Equivalent change (1)}
        \abs{V_k[f](w)}^p & \lesssim (\abs{V_k[f](w) - f_{X_e}}^p + \abs{f_{X_e}}^p)\\
        & \leq \kint_{X_e} \abs{f - V_k[f](w)}^p \,d\mu + \kint_{X_e} \abs{f}^p \, d\mu && \text{(Jensen's IE)}\nonumber\\
        & \lesssim \mathsf{E}_p(f \circ F_e) + \frac{1}{\mu(X_e)} \int_{X_e} \abs{f}^p \, d\mu && \text{(Poincar\'e IE \eqref{eq: Poincare for C(X)})} \nonumber\\
        & \leq \mathsf{E}_p(f \circ F_e) + \abs{E_1}^k\norm{f}_{L^p}^p. && \text{(Proposition \eqref{prop: Geom of LS}-\eqref{item: diams/measures of cells})} \nonumber
    \end{align}
    Notice that \eqref{eq: Equivalent change (1)} is almost the objective \eqref{eq: Equivalent change}. The difference is that the coefficients do not match. To fix this, we show that the inequality \eqref{eq: Equivalent change (1)} has a \emph{self-improving property}, which means that we can alter the coefficients of the two terms to more suitable ones.
    
    Since the gluing sets are independent, the edges $\fe(v_\pm \cdot a), \fe(v_\pm \cdot b) \in E_1$ are distinct whenever $a,b \in I$ are distinct. Also the edges $\fe(v_+ \cdot a), \fe(v_- \cdot b) \in E_1$ are distinct for all $a,b \in I$ due to the non-degenerateness. Thus, using the duality relation \eqref{eq:duality} and the conductive uniform property, we estimate
    \begin{align*}
        \left(\cM_p^{-1}\right)^{\frac{q}{p}} & = \sum_{e \in E_1} \abs{\mathcal{J}_p(e)}^q \geq \sum_{a \in I} \left(\abs{\mathcal{J}_p(\fe(v_- \cdot a))}^q + \abs{\mathcal{J}_p(\fe(v_+ \cdot a))}^q\right) = 2 \cdot \sum_{a \in I} \mathcal{J}_p(a)^q.
    \end{align*}
    Using Lemma \ref{lemma: Properties of cond.unif.}-\eqref{item: Properties of cond.unif. (2)} we obtain
    \begin{align}
    \label{eq: Equivalent change (2)}
        \sum_{a \in I^m} \mathcal{J}_{p,m}(a)^q = \left(\sum_{a \in I} \mathcal{J}_p(a)^q\right)^m \leq \frac{\left(\cM_p^{-m}\right)^{\frac{q}{p}}}{2^m}.
    \end{align}
    Now fix any edge $e \in E_n$ containing $v$.
    By combining the inequalities \eqref{eq: Equivalent change (1)} and \eqref{eq: Equivalent change (2)} we obtain the estimate
    \begin{align*}
        \abs{V_n[f](v)}^p & = \abs{V_{n+m,n}[V_{n+m}[f]](v)}^p && \text{(Tower rule)}\\
       & = \left|\sum_{a \in I^m}\mathcal{J}_{p,m}(a) V_{n+m}[f](v \cdot a) \right|^p\\
        & \leq \left(\sum_{a \in I^m} \mathcal{J}_{p,m}(a)^q \right)^{\frac{p}{q}} \left( \sum_{a \in I^{m}}\abs{V_{n+m}[f](v \cdot a)}^p \right) && \text{(Hölder's IE)}\\
        & \lesssim \frac{\cM_p^{-m}}{2^{m(p-1)}} \sum_{a \in I^m} \left(\mathsf{E}_p(f \circ F_{\fe(v \cdot a,e)}) + \abs{E_1}^{n+m}\norm{f}_{L^p}^p  \right)&& \text{\eqref{eq: Equivalent change (1)}-\eqref{eq: Equivalent change (2)}} \\
        & \leq \frac{\cM_p^{n}}{2^{m(p-1)}} \mathsf{E}_p(f) + \frac{\abs{I}^m\abs{E_1}^{n+m}}{\cM_p^{m}2^{m(p-1)}}\norm{f}_{L^p}^p. && \text{(Lemma \ref{lemma: Energy self-similarity})}
    \end{align*}
    In the last line, we used that the mapping $a \mapsto \fe(v \cdot a,e)$, where $\fe(v \cdot a,e)$ is the unique edge in $e \cdot G_m$ containing $v \cdot a$, is injective. This holds by Corollary \ref{cor: Neigbours of ancestors}.
\end{proof}

We are now ready to prove the closability.

\begin{proof}[Proof of Theorem \ref{thm: Closability}]
    Let $\{ f_i \}_{i = 1}^\infty \subseteq \mathscr{C}_p$ be $\mathsf{E}_p$-Cauchy satisfying $f_i \to 0$ in $L^p(X,\mu)$.
    First, we show for any fixed $n \in \N \cup \{0\}$ and $v \in V_n$ that we have
    \begin{equation}\label{eq: V_n[f] to 0}
        V_n[f_i](v) \xrightarrow[]{i \to \infty} 0.
    \end{equation}
    Fix $\varepsilon > 0$.
    Since $\{ f_i \}_{i = 1}^{\infty} \subseteq \mathscr{C}_p$ is $\mathsf{E}_p$-Cauchy, there is a uniform upper bound for the energies
    \[
        \mathsf{E}_p(f_i) \leq M < \infty.
    \]
    Next, let $C_0 > 0$ be as in Proposition \ref{prop: Equivalent change}, and choose $m \in \N$ so large that
    \[
        \frac{M  \cM_p^n }{2^{m(p-1)}} < \frac{\varepsilon}{2C_0}.
    \]
    Since $f_i \xrightarrow[]{} 0$ in $L^p(X\mu)$, there is $k \in \N$ so that for any $i \geq k$ we have
    \[
        \frac{\abs{I}^m\abs{E_1}^{n+m}}{\cM_p^{m}2^{m(p-1)}}\norm{f_i}^p_{L^p} < \frac{\varepsilon}{2C_0}.
    \]
    For any such $i \in \N$, by Proposition \ref{prop: Equivalent change}, we have
    \[
        \abs{V_n[f_i](v)}^p < \varepsilon.
    \]
    This proves \eqref{eq: V_n[f] to 0}, and consequently, we have
    \[
        \mathcal{E}_p^{(n)}(f_i) = \cM_p^{-n}\mathcal{E}_p(V_n[f_i]) \xrightarrow[]{i \to \infty} 0.
    \]

    Next, we show that $\mathsf{E}_p(f_i) \to 0$. Fix $N(\varepsilon) \in \N$ so that
    \[
        \sup_{i,j \geq N(\varepsilon)} \mathsf{E}_p(f_i - f_j)^{\frac{1}{p}} < \frac{\varepsilon}{2}.
    \]
    By Theorem \ref{thm: Strong monotonicity}, it is sufficient to prove that $\mathcal{E}_p^{(n)}(f_i) < \varepsilon$ for all $n \in \N$ and $i \geq N(\varepsilon)$.
    We achieve this by choosing $j \geq N(\varepsilon)$ so that $\cE_p^{(n)}(f_j)^{\frac{1}{p}} < \frac{\varepsilon}{2}$, and using Minkowski's inequality and Theorem \ref{thm: Strong monotonicity} to estimate
    \begin{align*}
        \cE_p^{(n)}(f_i)^{\frac{1}{p}}  \leq \cE_p^{(n)}(f_i - f_j)^{\frac{1}{p}} + \cE_p^{(n)}(f_j)^{\frac{1}{p}}
         \leq \mathsf{E}_p(f_i - f_j)^{\frac{1}{p}} + \cE_p^{(n)}(f_j)^{\frac{1}{p}}
        < \varepsilon.
    \end{align*}
\end{proof}

\begin{remark}
    Hereafter, up to the end of the paper, $\mathscr{E}_p : L^p(X,\mu) \to [0,\infty]$ denotes the $p$-energy form and $\mathscr{F}_p$ the associated Sobolev space, as given in Definition \ref{def: Sobolev}.
\end{remark}

We now move on to establishing the main theorems of the section, Theorems \ref{Thm: Existence of p-Energy} and \ref{thm: Energy analytic properties}.
Proving Theorem \ref{Thm: Existence of p-Energy}-\eqref{item: Contraction} regarding the Lipschitz contractivity requires more careful arguments so we postpone it to the next subsection.

\begin{proof}[Proof of Theorem \ref{Thm: Existence of p-Energy}; \eqref{item: Energy is explicit}-\eqref{item: p-energySS}]
    \eqref{item: Energy is explicit}: If $f \in \mathscr{C}_p$, then $\mathscr{E}_p(f) \leq \mathsf{E}_p(f)$ follows by taking the constant sequence $\{f_i\}_{i = 1}^\infty$, $f_i = f$ for all $i \in \N$. The opposite inequality follows from Theorem \ref{thm: Closability} and by sending $n\to\infty$ in the estimate $\cE_p^{(n)}(f)^{\frac{1}{p}} \leq \cE_p^{(n)}(f-f_i)^{\frac{1}{p}}+ \cE_p^{(n)}(f_i)^{\frac{1}{p}}$, which follows from Minkowski's inequality.
    
    \eqref{item: Refl-Sep}: The Sobolev space $\mathscr{F}_p$ can be regarded as the completion of $\mathscr{C}_p$ with respect to the Sobolev norm. Also note that $\mathscr{E}_p(\cdot)^{\frac{1}{p}}$ satisfies the triangle inequality because the discrete ones $\mathcal{E}_p^{(n)}(\cdot)^{\frac{1}{p}}$ do. Thus, we have a Banach space.
    Reflexivity and separability of $\mathscr{F}_p$ follow once we have established $p$-Clarkson's inequality \eqref{e:Cp} (see \cite[Proposition 3.13 and Corollary 3.16]{KajinoContraction} for details), which is covered in Theorem \ref{thm: Energy analytic properties}. Separability can also be directly seen from Theorem \ref{thm: Mollifiers converge} below.

    \eqref{item: Regularity}: The density of $\mathscr{C}_p \subseteq \mathscr{F}_p$ is clear from the construction. The density $\mathscr{C}_p \subseteq C(X)$ in the uniform norm $\norm{\cdot}_{L^\infty}$follows from \eqref{item: Properties of mollifier (1)} and \eqref{item: Properties of mollifier (4)} of Proposition \ref{prop: Properties of mollifier}.

    \eqref{item: p-energySS}:
    It is sufficient to show that whenever $\{f_i\}_{i = 1}^{\infty} \subseteq \mathscr{C}_p$ is an $\mathsf{E}_p$-Cauchy sequence then so is $\{f_i \circ F_e \}_{i = 1}^{\infty}$ for all $e \in E_{n}$. The self-similarity of $\mathscr{E}_p$ is then be inherited from the self-similarity of the pre-energy form $\mathsf{E}_p$ established in Lemma \ref{lemma: Energy self-similarity}.
    
    Let $\{f_i\}_{i = 1}^{\infty} \subseteq \mathscr{C}_p$ be $\mathsf{E}_p$-Cauchy and $e \in E_n$.
    By Theorem \eqref{item: Energy is explicit} of the current theorem and Lemma \ref{lemma: Energy self-similarity} we get
    \[
    \mathscr{E}_p((f_i - f_j) \circ F_e) = \mathsf{E}_p((f_i - f_j) \circ F_e) \leq \cM_p^n \mathsf{E}_p(f_i - f_j) \xrightarrow[]{i,j\to \infty} 0.
    \]
\end{proof}

\begin{proof}[Proof of Theorem \ref{thm: Energy analytic properties}]
    \eqref{item: Poincare}:
    Since $\mathscr{C}_p \subseteq \mathscr{F}_p$ is dense, it is sufficient to prove \eqref{eq: PrePI (1)} only for the functions $f \in \mathscr{C}_p$. Now, by Minkowski's inequality,
    \begin{equation}\label{eq: const in PI}
        \left(\kint_A \abs{f - f_A}^p \, d\mu \right)^{\frac{1}{p}} \leq 2 \left(\inf_{c \in \R} \kint_A \abs{f - c}^p \, d\mu \right)^{\frac{1}{p}},
    \end{equation}
    where $A \subseteq X$ is any Borel set of positive and finite measure.
    Thus \eqref{eq: PrePI (1)} follows from Lemma \ref{lemma: Poincare for C(X)} and \eqref{eq: const in PI}.

    \eqref{item: p-Clarkson}: Since the operators $V_n[\cdot]$ are linear, for any $n \in \mathbb{N} \cup \{0\}$ and $f,g \in \mathscr{C}_{p}$ we have \eqref{e:Cp} with $\mathcal{E}_{p}^{(n)}$ in place of $\mathscr{E}_{p}$ by $p$-Clarkson's inequality on $\ell^p(V_n)$. We obtain \eqref{e:Cp} in the case $f,g \in \mathscr{C}_{p}$ by letting $n \to \infty$.
    For general $f \in \mathscr{F}_p$, the claim follows by taking approximating sequences of continuous functions.

    \eqref{item: strong-locality}:
    Since the closed sets $\supp_{\mu}(f - a)$ and $\supp_{\mu}(g - b)$ are disjoint, it follows from Proposition \ref{prop: Geom of LS}-\eqref{item: diams/measures of cells} that we can choose a large $n \in \N$ and a subset of edges $H_n \subseteq E_n$ so that
    \[
        \supp_{\mu}(f - a) \subseteq \bigcup_{e \in H_n} X_e \text{ and } \supp_{\mu}(f - b) \subseteq \bigcup_{e \in E_n \setminus H_n} X_e.
    \]
    It is clear from the construction of $\mathscr{E}_p$ that $\mathscr{E}_p(f + c1_X) = \mathscr{E}_p(f)$ for all $c \in \R$. This fact combined with Theorem \ref{Thm: Existence of p-Energy}-\eqref{item: p-energySS} yield
    \begin{align*}
        & \, \quad \mathscr{E}_p(f + g + h) \\
        = & \quad \cM_p^{-n}\sum_{e \in H_n} \mathscr{E}_p((f + g + h) \circ F_e)
         + \cM_p^{-n}\sum_{e \in E_n \setminus H_n} \mathscr{E}_p((f + g + h) \circ F_e) \\
        = & \quad \cM_p^{-n}\sum_{e \in H_n} \mathscr{E}_p((f + h) \circ F_e)
         + \cM_p^{-n}\sum_{e \in E_n \setminus H_n} \mathscr{E}_p((g + h) \circ F_e).
    \end{align*}
    Similarly we compute
    \begin{align*}
        \mathscr{E}_p(h) & = \cM_p^{-n}\sum_{e \in H_n} \mathscr{E}_p(h \circ F_e) + \cM_p^{-n} \sum_{e \in E_n \setminus H_n} \mathscr{E}_p(h \circ F_e)\\
        & = \cM_p^{-n}\sum_{e \in H_n} \mathscr{E}_p((g + h) \circ F_e) + \cM_p^{-n} \sum_{e \in E_n \setminus H_n} \mathscr{E}_p((f + h) \circ F_e).
    \end{align*}
    The desired claim now follows by combining the previous two equalities using Theorem \ref{Thm: Existence of p-Energy}-\eqref{item: p-energySS} one last time.
\end{proof}

\subsection{Contraction properties}\label{subsec: contraction}
Next we finish the proof of Theorem \ref{Thm: Existence of p-Energy} by establishing the Lipschitz contraction property.
% We also record the proof of the \emph{generalized $p$-contraction property} introduced in \cite{KajinoContraction} for future references.
The proof turned out to be quite delicate unlike, for instance, the proof of $p$-Clarkson's inequality, which was easy to prove for the functions in $\mathscr{C}_p$ and the general case immediately followed from the density $\mathscr{C}_p \subseteq \mathscr{F}_p$.
The main difficulty comes with the fact that the averaging operators $V_n[\cdot]$ and the composition $f \mapsto \varphi \circ f$ do not commute in general. More precisely, the equality of functions
\begin{equation}\label{eq: Average and contraction commute}
    \varphi \circ V_n[f] = V_n[\varphi \circ f]
\end{equation}
does not need to hold for an arbitrary $f \in \mathscr{C}_p$. The key observation for our argument is that \eqref{eq: Average and contraction commute} still holds for a large family of $f \in \mathscr{C}_p$, which we will prove in the following lemma.

\begin{lemma}\label{lemma: Average and contraction commute}
    Let $m \in \N \cup \{0\}$, and $f := \mathscr{U}_{p,m}[g]$ for any $g : V_m \to \R$. Then \eqref{eq: Average and contraction commute} holds for all $n \geq m$ and 1-Lipschitz functions $\varphi : \R \to \R$. Moreover, it holds that $\mathscr{E}_p(\varphi \circ f) \leq \mathscr{E}_p(f)$.
\end{lemma}

\begin{proof}
    It follows from Proposition \ref{prop: Properties of mollifier}-\eqref{item: Properties of mollifier (2)} that $f$ is constant on each fiber $\Fib(v)$ for $v \in V_n$ and $n \geq m$.
    Thus, if $x \in \Fib(v)$ is any point, we have
    \[
        V_n[\varphi \circ f](v) = \varphi(f(x)) = \varphi(V_n[f](v)).
    \]
    Note that the discrete graph energies satisfy the Lipschitz contraction property, i.e. for all $n \in \N \cup \{0\}$ and $h : V_n \to \R$ we have $\mathcal{E}_p(\varphi \circ h) \leq \mathcal{E}_p(h)$. Thus we finish the proof by computing
    \begin{align*}
        \mathscr{E}_p(\varphi \circ f) & = \lim_{n \to \infty} \cM_p^{-n}\mathcal{E}_p (V_n[\varphi \circ f])\\
        & = \lim_{n \to \infty} \cM_p^{-n}\mathcal{E}_p(\varphi \circ V_n[f])\\
        & \leq \lim_{n \to \infty} \cM_p^{-n}\mathcal{E}_p(V_n[f]) = \mathscr{E}_p(f).
    \end{align*}
\end{proof}

The next step is the lower-semicontinuity of $\mathscr{E}_p$, which is a direct corollary of $p$-Clarkson's inequality \eqref{e:Cp} and \cite[Proposition 3.18-(a)]{KajinoContraction}. We omit the details.

\begin{lemma}\label{lemma: Lower-semicontinuity of energy}
    Let $f \in L^p(X,\mu)$ and $\{f_i\}_{i = 1}^{\infty}$ be a sequence of $L^p(X,\mu)$-functions so that $f_i \to f$ in $L^p(X,\mu)$. Then
    \[
        \mathscr{E}_p(f) \leq \liminf_{i \to \infty} \mathscr{E}_p(f_i).
    \]
    If additionally $\{ f_i \}_{i = 1}^\infty \cup \{ f \} \subseteq \mathscr{F}_p$ and
    \[
        \lim_{i \to \infty} \mathscr{E}_p(f_i) = \mathscr{E}_p(f),
    \]
    then $f_i \to f$ in $\mathscr{F}_p$.
\end{lemma}

We are now ready to finish the proof of Theorem \ref{Thm: Existence of p-Energy}.

\begin{proof}[Proof of Theorem \ref{Thm: Existence of p-Energy}-\eqref{item: Contraction}]
    We first assume that $f \in \mathscr{C}_p$, and consider the approximating functions $\varphi \circ \Psi_{p,n}[f]$ for $n \in \N \cup \{0\}$. It follows from Proposition \ref{prop: Properties of mollifier}-\eqref{item: Properties of mollifier (4)} that $\varphi \circ \Psi_{p,n}[f] \to \varphi \circ f$ in $L^p(X,\mu)$.
    Hence, we obtain the desired contraction property by computing
    \begin{align}
    \label{comp:Contraction}
        \mathscr{E}_p(\varphi \circ f) & \leq \liminf_{n \to \infty} \mathscr{E}_p(\varphi \circ \Psi_{p,n}[f]) && \text{(Lemma \ref{lemma: Lower-semicontinuity of energy})} \\
        & \leq \liminf_{n \to \infty} \mathscr{E}_p(\Psi_{p,n}[f]) && \text{(Lemma \ref{lemma: Average and contraction commute})} \nonumber \\
        & = \mathscr{E}_p(f). && \text{(Proposition \ref{prop: Properties of mollifier}-\eqref{item: Properties of mollifier (1)})} \nonumber
    \end{align}
    For a general $f \in \mathscr{F}_p$, we argue by taking an approximating sequence of continuous functions $\{ f_n \}_{n = 1}^{\infty}$ in $\mathscr{F}_p$. Since $\varphi \circ f_n \to \varphi \circ f$ in $L^p(X,\mu)$ and $\varphi \circ f_n \in \mathscr{F}_p$ are continuous, it follows from Lemma \ref{lemma: Lower-semicontinuity of energy} and the fact that we have verified the contraction property for continuous Sobolev functions that
    \begin{equation}\label{eq: Contraction to general sobolev}
        \mathscr{E}_p(\varphi \circ f) \leq\liminf_{n \to \infty} \mathscr{E}_p(\varphi \circ f_n) \leq \liminf_{n \to \infty} \mathscr{E}_p(f_n) = \mathscr{E}_p(f).
    \end{equation}
\end{proof}

\begin{remark}
    The previous proof can be modified to establish a stronger notion of contractivity, the \emph{generalized $p$-contraction property}, introduced by Kajino and the third author \cite{KajinoContraction}.
\end{remark}

We shall finish the discussion on the contraction property by recording that, for $p = 2$, our construction produces a Dirichlet form; we refer to \cite{fukushima2011dirichlet} for background.

\begin{corollary}
    Consider the two variable functional $\widetilde{\mathscr{E}}_2 : \mathscr{F}_2 \times \mathscr{F}_2 \to \R$ given by
    \[
        \widetilde{\mathscr{E}}_2(f,g) := \frac{1}{4} (\mathscr{E}_2(f + g) - \mathscr{E}_2(f - g)).
    \]
    Then $\widetilde{\mathscr{E}}_2$ is a strongly local, regular Dirichlet form on $L^2(X,\mu)$.
\end{corollary}

\begin{proof}
    By Theorems \ref{Thm: Existence of p-Energy} and \ref{thm: Energy analytic properties}, the two variable functional $\widetilde{\mathscr{E}}_2$ is strongly loca, regular, Markovian and closed. Thus, we only need to show that it is bilinear.
    To see this, observe that the bilinearity clearly holds if we replace $\mathscr{E}_2$ with $\mathcal{E}_2^{(n)}$.
    Thus, the bilinearity for $f,g \in \mathscr{C}_p$ follows from Theorem \ref{Thm: Existence of p-Energy}-\eqref{item: Energy is explicit}. The general case then follows by taking approximating sequences of continuous functions. 
\end{proof}

\begin{remark}\label{rem:PHI}
    It follows from \cite{ResistanceConjecture} and  that the strongly local, regular Dirichlet form $(\widetilde{\mathcal{E}}_2,\mathscr{F}_2)$ on $L^2(X,\mu)$ satisfies sub-Gaussian heat kernel estimates. The required hypotheses to apply \cite{ResistanceConjecture} are provided by Propositions \ref{prop: Geom of LS}-\eqref{item: Q-AR}, \ref{prop: EM-Poincare} and \ref{prop: EM-Ucap}.
\end{remark}

\subsection{Extensions of operators}\label{subsec: Extending operators}
Next we extend the operators of Section \ref{Sec: Discretizations and mollifiers} to general Sobolev functions. This is important because many arguments in later sections heavily rely on approximation arguments.

First, we introduce for all $n \in \N \cup \{0\}$ the convenient norm on $\R^{V_n}$,
\[
  \norm{\cdot}_{\mathcal{E}_p} := \abs{E_1}^{-\frac{n}{p}}\norm{\cdot}_{\ell^p} + (\cM_p)^{-\frac{n}{p}}\mathcal{E}_p(\cdot)^{\frac{1}{p}}.  
\]
\begin{lemma}\label{lemma: How to extend}
    There is a constant $C \geq 1$ so that for all $n \in \N \cup \{0\}, f \in \mathscr{C}_p$ and $g \in \R^{V_n}$ we have
    \begin{equation}
        \norm{V_n[f]}_{\mathcal{E}_p} \leq C\norm{f}_{\mathscr{F}_p} \text{ and } \norm{\mathscr{U}_{p,n}[g]}_{\mathscr{F}_p} \leq C \norm{g}_{\mathcal{E}_p}.
    \end{equation}
\end{lemma}

\begin{proof}
    We first derive the estimate for $\norm{V_n[f]}_{\ell^p}$.
    It follows from Theorem \ref{thm: Strong monotonicity} that $\cM_p^{-n}\mathcal{E}_p(V_n[f]) \leq \mathscr{E}_p(f)$.
    The $\ell^p$-norm of $V_n[f]$ can be estimated  by using the Poincar\'e inequality \eqref{eq: Poincare for C(X)} and Proposition \ref{prop: Geom of LS}-\eqref{item: diams/measures of cells},
    \begin{align*}
        \norm{V_n[f]}_{\ell^p} = \left(\sum_{v \in V_n}\abs{V_n[f](v)}^p\right)^{\frac{1}{p}} & \leq \left(\sum_{v \in V_n}\abs{V_n[f](v) - f_{X_v}}^p\right)^{\frac{1}{p}} + \left(\sum_{v \in V_n} \abs{f_{X_v}}^p\right)^{\frac{1}{p}} \\
        & \lesssim (\cM_p)^{\frac{n}{p}} \mathscr{E}_p(f)^{\frac{1}{p}} + \abs{E_1}^{\frac{n}{p}} \norm{f}_{L^p}.
    \end{align*}
    Since $\cM_p < \abs{E_1}$ by Lemma \ref{lemma: Properties of cond.unif.}-\eqref{item: Properties of cond.unif. (3)}, we have
    \[
        \norm{V_n[f]}_{\mathcal{E}_p} \lesssim \left(\frac{\cM_p^n}{\abs{E_1}^{n}}\mathscr{E}_p(f)\right)^{\frac{1}{p}} + \mathscr{E}_p(f)^{\frac{1}{p}} + \norm{f}_{L^p} \lesssim \norm{f}_{\mathscr{F}_p}.
    \]

    We move on to estimating $\norm{\mathscr{U}_{p,n}[g]}_{\mathscr{F}_p}$.
    It follows from Lemma \ref{lemma: Properties of U_p}-\eqref{item: Properties of U_p (1)} that $\mathscr{E}_p(\mathscr{U}_{p,n}[g]) = \cM_p^{-n}\mathcal{E}_p(g)$. Using Lemma \ref{lemma: Properties of U_p}-\eqref{item: Properties of U_p (2)} we estimate
    \begin{align*}
        \left( \int_{X} \abs{\mathscr{U}_{p,n}[g]}^p \, d\mu \right)^{\frac{1}{p}} & \leq \sum_{e \in E_n} \left( \int_{X_e} \abs{g(e^+)}^p + \abs{g(e^-)}^p \, d\mu \right)^{\frac{1}{p}} \lesssim \abs{E_1}^{-\frac{n}{p}} \norm{g}_{\ell^p}.
    \end{align*}
\end{proof}

Thus, the linear operators $V_n[\cdot]$ restricted to the core $\mathscr{C}_p$ are all bounded. Because $\mathscr{F}_p$ is a Banach space and $\mathscr{C}_p \subseteq \mathscr{F}_p$ is dense, the domains of $V_n[\cdot]$ can be uniquely extended to the whole Sobolev space. However, there is a certain ambiguity involving these extensions.
Indeed the extended operator, let us denote it $\widetilde{V}_n[\cdot] : \mathscr{F}_p \to \R^{V_n}$, on principle could differ from $V_n[\cdot]$, the operator defined in Section \ref{Sec: Discretizations and mollifiers}, at the functions in $f \in (\mathscr{F}_p \cap C(X))\setminus \mathscr{C}_p$. 
Fortunately, no such $f$ exists, but the proof is quite delicate because we use reflexivity of $\mathscr{F}_p$; see Corollary \ref{cor: Self-similarity of domain}. But first we need additional tools.

\begin{definition}\label{def: G_n operators}
Let $n \in \N \cup \{0\}$. We define $M_n[\cdot] : L^p(X,\mu) \to \R^{V_n}$ to be the linear operator
\[
    M_n[f](v) := \kint_{X_v} f\, d\mu \text{ for all } v \in V_n.
\]
We then define the \emph{rough mollifier} $\Xi_{p,n}[\cdot]: L^p(X,\mu) \to \mathscr{F}_p$ as the composition $\mathscr{U}_{p,n}[\cdot] \circ M_n[\cdot]$. 
\end{definition}

The operators $M_n[\cdot]$ and $\Xi_{p,n}[\cdot]$ have one major advantage to the operators in Section \ref{Sec: Discretizations and mollifiers}. Namely these can be defined for a general $L^p$ functions and they do not require the additional extension step to $\mathscr{F}_p$.

\begin{lemma}\label{lemma: G_n mollifier}
    The following conditions hold for the operators $M_n[\cdot]$ and $\Xi_{p,n}[\cdot]$.
    \begin{enumerate}
        \vskip.3cm
        \item \label{item: G_n mollifier bounded}
        The restrictions of $M_n[\cdot]$ and $\Xi_{p,n}[\cdot]$ for $n \in \N \cup \{0\}$ to the Sobolev space $\mathscr{F}_p$ are uniformly bounded with respect to the norms $\norm{\cdot}_{\mathcal{E}_p}$ and $\norm{\cdot}_{\mathscr{F}_p}$.
        \vskip.3cm
        \item \label{item: G_n mollfier converge}
        For all $f \in L^p(X,\mu)$ it holds that $\Xi_{p,n}[f] \to f$ as $n \to \infty$ in $L^p(X,\mu)$. If additionally $f \in C(X)$, then $\Xi_{p,n}[f] \to f$ uniformly.
    \end{enumerate}
\end{lemma}

\begin{proof}
    \eqref{item: G_n mollifier bounded}: This follows from an identical argument as in Lemma \ref{lemma: How to extend}.
    
    \eqref{item: G_n mollfier converge}: 
    Recall from Proposition \ref{prop: Geom of LS}-\eqref{item: diams/measures of cells} that $\diam(X_v) \leq 2L_*^{-n}$.
    It is now routine to check that for all continuous functions $f \in C(X)$ we have $\Xi_{p,n}[f] \to f$ uniformly using Lemma \ref{lemma: Properties of U_p}-\eqref{item: Properties of U_p (2)} and the uniform continuity of $f$.
    Then assume $f \in L^p(X,\mu)$. By noting that the linear operators $\Xi_{p,n}[\cdot] : L^p(X,\mu) \to L^p(X,\mu)$ are uniformly bounded, the desired convergence $\Xi_{p,n}[f] \to f$ follows by taking approximating sequences of continuous functions.
\end{proof}

\begin{corollary}\label{cor: Self-similarity of domain}
    It holds that $\mathscr{C}_p = C(X) \cap \mathscr{F}_p$. In particular, the equality $\mathscr{E}_p(f) = \mathsf{E}_p(f)$ holds for all $f\in C(X)$.
    Furthermore, we also have
    \begin{equation}\label{eq: Self-similarity of domain}
        C(X) \cap \mathscr{F}_p = \{ f \in C(X) : f \circ F_e \in \mathscr{F}_p \text{ for all } e \in E_1 \}.
    \end{equation}
\end{corollary}

\begin{proof}
    It is trivial that $\mathscr{C}_p \subseteq C(X) \cap \mathscr{F}_p$. To prove the other inclusion, fix an arbitrary continuous Sobolev function $f \in C(X) \cap \mathscr{F}_p$.

    Let us denote $\widetilde{V}_n[\cdot] : \mathscr{F}_p \to \R^{V_n}$ the extension of $V_n[\cdot]|_{\mathscr{C}_p}$, which exists due to Lemma \ref{lemma: How to extend},
    and also recall the averaging operators $V_n[\cdot] : C(X) \to \R^{V_n}$ given in Definition \ref{def: Discretization operator}.
    We will show that $V_n[f] = \widetilde{V}_n[f]$ for all $n \in \N \cup \{0\}$. It would then follow from Lemmas \ref{lemma: Properties of U_p}-\eqref{item: Properties of U_p (1)} and \ref{lemma: How to extend},
    \begin{align*}
        \mathsf{E}_p(f) & = 
        \lim_{n \to \infty} \cM_p^{-n}\mathcal{E}_p(V_n[f]) = \lim_{n \to \infty} \cM_p^{-n} \mathcal{E}_p(\widetilde{V}_n[f])\\
        & \lesssim \limsup_{n \to \infty} \norm{\widetilde{V}_n[f]}_{\mathcal{E}_p}^p \lesssim \norm{f}_{\mathscr{F}_p}^p < \infty.
    \end{align*}
    By the definition of $\mathscr{C}_p$ given in \eqref{eq:CORE}, we would have $f \in \mathscr{C}_p$.
    
    In order to verify the desired equality $V_n[f]=\widetilde{V}_n[f]$, we consider the rough mollifiers $\Xi_{p,n}[\cdot]$.
    Since $\Xi_{p,n}[f] \to f$ uniformly by Lemma \ref{lemma: G_n mollifier}-\eqref{item: G_n mollfier converge}, it follows from the definition of $V_{n}[\cdot]$ that $V_n[\Xi_{p,m}[f]] \to V_n[f]$ as $m \to \infty$.
    Using the fact that $\Xi_{p,m}[f] \in \mathscr{C}_p$,
    \begin{equation}\label{eq: Discs agree}
    \widetilde{V}_n[\Xi_{p,m}[f]] = V_n[\Xi_{p,m}[f]].
    \end{equation}
    On the other hand, by Lemma \ref{lemma: G_n mollifier}-\eqref{item: G_n mollifier bounded}, the sequence $\{\Xi_{p,m}[f]  \}_{m = 0}^\infty$ is bounded in $\mathscr{F}_p$.
    Since $\mathscr{F}_p$ is reflexive by Theorem \ref{Thm: Existence of p-Energy}-\eqref{item: Refl-Sep}, it follows from Mazur's lemma (see \cite[Page 19 and Theorem 2.41]{HKST}) that there is a convergent sequence of convex combinations,
\[
    \sum_{j = m}^{N_m} \lambda_{j,m}\Xi_{p,m}[f] \xrightarrow[]{m \to \infty} f \text{ in } \mathscr{F}_p.
\]
By combining \eqref{eq: Discs agree}, the linearity and continuity of $\widetilde{V}_{n}[\cdot]$ along with the fact that $V_n[\Xi_{p,m}[f]] \to V_n[f]$,
    \[
    \widetilde{V}_{n}[f] = \lim_{m \to \infty} \sum_{j = m}^{N_m} \lambda_{j,m}\widetilde{V}_{n}[\Xi_{p,m}[f]] = \lim_{m \to \infty} \sum_{j = m}^{N_m} \lambda_{j,m}V_{n}[\Xi_{p,m}[f]] = V_{n}[f].
    \]
Thus, we conclude that $f \in \mathscr{C}_p$.

It is now clear from Theorem \ref{Thm: Existence of p-Energy}-\eqref{item: Energy is explicit} that the equality $\mathscr{E}_p(f) = \mathsf{E}_p(f)$ holds for all $f \in C(X)$. Indeed, for $f \in C(X) \setminus \mathscr{C}_p$, both values read $\infty$.
Thus, \eqref{eq: Self-similarity of domain} now follows from Lemma \ref{lemma: Energy self-similarity}.
\end{proof}

Since we have now cleared the ambiguity involving the discretization operator, we now extend the operators from Section \ref{Sec: Discretizations and mollifiers} to the whole Sobolev space using Lemma \ref{lemma: How to extend}. Let us denote them accordingly $V_n[\cdot] : \mathscr{F}_p \to \R^{V_n}$ and $\Psi_{p,n} = \mathscr{U}_{p,n}[\cdot] \circ V_n[\cdot] : \mathscr{F}_p \to \mathscr{C}_p$.
We now finish the subsection to establishing a useful convergence result.

\begin{theorem}\label{thm: Mollifiers converge}
    Let $V_n[\cdot]$ and $\Psi_{p,n}[\cdot]$ be the extended linear operators as discussed above.
    Then for all $f \in \mathscr{F}_p$, 
    \begin{equation}\label{eq: Mollifications converge}
        \Psi_{p,n}[f] \xrightarrow[]{n \to \infty} f \text{ in } \mathscr{F}_p.
    \end{equation}
\end{theorem}

\begin{proof}
    We first derive \eqref{eq: Mollifications converge} for the functions in the core $f \in \mathscr{C}_p$.
    By \eqref{item: Properties of mollifier (1)} and \eqref{item: Properties of mollifier (4)} of Proposition \ref{prop: Properties of mollifier} we have
    \[
        \mathscr{E}_p(\Psi_{p,n}[f]) = \cE_p^{(n)}(f) \xrightarrow[]{n \to \infty} \mathscr{E}_p(f) \quad \text{and} \quad \Psi_{p,n}[f] \xrightarrow[]{n \to \infty} f \text{ in } L^p(X,\mu).
    \]
    It now follows from Lemma \ref{lemma: Lower-semicontinuity of energy} that \eqref{eq: Mollifications converge} holds for all functions in the core $f \in \mathscr{C}_p$. The general case follows from the density $\mathscr{C}_p \subseteq \mathscr{F}_p$, and that the operators $\Psi_{p,n}[\cdot]$ for $n \in \N \cup \{0\}$ are bounded by a constant independent of $n$.
    The latter follows from Lemma \ref{lemma: How to extend}.
\end{proof}

\begin{remark}\label{remark:Extension}
    Hereafter up to the end of the work, $V_n[\cdot]$ and $\Psi_{p,n}[\cdot]$ denote the extended operators. It follows from Lemma \ref{lemma: How to extend} and Theorem \ref{thm: Mollifiers converge} that many analogous properties proved in Section \ref{Sec: Discretizations and mollifiers} are inherited to the extensions.
\end{remark}

\subsection{Energy measures}\label{subsec:EM}
We finish the section in the construction of energy measures which are the natural counterparts of the measures $A \mapsto \int_{A} \abs{\nabla f}^p \, dx$ in the classical setting. Our method in the construction uses the self-similarity of $\mathscr{E}_p$, and similar approach was used in e.g. \cite{hino2005singularity,shimizu,murugan2023first}. See also \cite{Sasaya2025} for a construction that does not rely on the self-similarity.
These two approaches yield the same measures by \cite[Corollary 5.13]{KajinoContraction} and \cite[Corollary 5.6]{Sasaya2025}.

Given a Sobolev function $f \in \mathscr{F}_p$ we define its \emph{pre-energy measure} as the Radon measure $\mathfrak{m}_p\Span{f}$ on the symbol space $\Sigma$ given by the conditions
\[
    \mathfrak{m}_p\Span{f}(\Sigma_e) := \cM_p^{-n} \mathscr{E}_p(f \circ F_e) \, \text{ for all } e \in E_n.
\]
As it was discussed in Definition \ref{def: Measures on sigma}, these conditions uniquely determine a Radon measure on $\Sigma$ by Theorem \ref{Thm: Existence of p-Energy}-\eqref{item: p-energySS} and Kolmogorov's extension theorem.

\begin{definition}
    The \emph{$p$-energy measure} of $f \in \mathscr{F}_p$ is the Radon measure $\Gamma_p\Span{f}$ given by the push-forward $\Gamma_p\Span{f} :=  \chi_*(\mathfrak{m}_p\Span{f})$.
\end{definition}

We collect some useful properties of the $p$-energy measures.

\begin{lemma}\label{lemma: Properties of energy measures}
    The $p$-energy measures satisfy the following properties.
    \begin{enumerate}
        \vskip.3cm
        \item \label{item: Properties of energy measures (1)}
        For every $f \in \mathscr{F}_p$ we have $\Gamma_p\Span{f}(X) = \mathscr{E}_p(f)$.
        \vskip.3cm
        \item \label{item: Properties of energy measures (2)}
        For any Borel set $A \subseteq X$ and $f,g \in \mathscr{F}_p$,
        \[
            \Gamma_p\Span{f + g}(A)^{\frac{1}{p}} \leq \Gamma_p\Span{f}(A)^{\frac{1}{p}} + \Gamma_p\Span{g}(A)^{\frac{1}{p}}.
        \]
        \vskip.3cm
        \item \label{item: Properties of energy measures (3)}
        Let $A \subseteq X$ be a Borel set and $f,g \in \mathscr{F}_{p} \cap C(X)$. If $(f - g)|_{A}$ is a constant function on $A$, then $\Gamma_{p}\langle f \rangle(A) = \Gamma_{p}\langle g \rangle(A)$.
    \end{enumerate}
\end{lemma}

\begin{proof}
    \eqref{item: Properties of energy measures (1)}: This follows from $\Gamma_p\Span{f}(X) = \mathfrak{m}_p\Span{f}(\Sigma) = \mathscr{E}_p(f)$.

    \eqref{item: Properties of energy measures (2)}: 
    Because $\mathscr{E}_p(\cdot)^{\frac{1}{p}}$ is a semi-norm we have
    \begin{equation}\label{eq: TriagEQ Sigma}
        \mathfrak{m}_p\Span{f+g}(\Sigma_e)^{1/p}\leq \mathfrak{m}_p\Span{f}(\Sigma_e)^{1/p} + \mathfrak{m}_p\Span{g}(\Sigma_e)^{1/p} \text{ for all } e \in E_\#.
    \end{equation}
    Since the sets $\Sigma_e$ generate the topology of $\Sigma$, we have \eqref{eq: TriagEQ Sigma} for all Borel sets of $\Sigma$. Applying this to the Borel set $\chi^{-1}(A) \subseteq \Sigma$ yields the claim.

    \eqref{item: Properties of energy measures (3)}: This would follow from \cite[Corollary 5.15-(b)]{KajinoContraction}. For the convenience of the reader, we sketch a short argument.
    By \eqref{item: Properties of energy measures (2)} of the current lemma, it suffices to consider the case $g=0$. If $A\subseteq X$ is open, then for all $e\in E_\#$ satisfying $X_e\subseteq A$ we have $f\circ F_e = 0$. Thus, we have $\mathfrak{m}_p\Span {f}(\Sigma_e)=0$. Therefore, for the open set $\Omega=\chi^{-1}(A)$ we get:
    \[
    \Gamma_p\Span{f}(A)=\mathfrak{m}_p\Span{f}(\Omega) \leq \sum_{e\in E_\#, \Sigma_e\subset \Omega} \mathfrak{m}_p\Span{f}(\Sigma_e)= 0.
    \]
    Now, if $A \subseteq X$ is a general Borel set, then $f=0$ on $\overline{A}$. Consider the approximations
    \[
    f_n:=\max(f,n^{-1})+\min(f,-n^{-1}).
    \]
    Clearly $f_n\to f$ in $L^p(X,\mu)$. Since the map $t \mapsto \max\{ t,c \} + \min \{ t,-c \}$ is a contraction for all $c \in \R$, it follows from the Lipschitz contraction property of $\mathscr{E}_p$ (Theorem \ref{Thm: Existence of p-Energy}-\eqref{item: Contraction}) that $\mathscr{E}_p(f_n)\leq \mathscr{E}_p(f)$ for all $n\in \N$. These together with Lemma \ref{lemma: Lower-semicontinuity of energy} imply that $f_n \to f$ in $\mathscr{F}_p$. Thus, $\lim_{n\to\infty}\Gamma_p\langle f_n \rangle(A)=\Gamma_p\langle f \rangle(A)$ by \eqref{item: Properties of energy measures (2)} of the current lemma (see also Lemma \ref{lemma: Abs.cont. of energy measures} below). But $\Gamma_p\langle f_n \rangle(A)=0$, since $f_n=0$ in a neighborhood of $A$.
\end{proof}

\begin{lemma}\label{lemma: Abs.cont. of energy measures}
    Let $f \in\mathscr{F}_p$ and $\{f_i\}_{i = 1}^{\infty}$ be any sequence of Sobolev functions so that $f_i \to f$ in $\mathscr{F}_p$. Then 
    \[
        \Gamma_p\Span{f_i}(A) \xrightarrow[]{i \to \infty} \Gamma_p\Span{f}(A) \text{ for all Borel sets } A \subseteq X.  
    \]
    In particular, if $\nu$ is a Borel measure of $X$ so that $\Gamma_p\Span{f_i} \ll \nu$ for all $n \in \N$, then $\Gamma_p\Span{f} \ll \nu$.
\end{lemma}

\begin{proof}
    Let $A \subseteq X$ be a Borel set. Then it follows from \eqref{item: Properties of energy measures (1)}-\eqref{item: Properties of energy measures (2)} of Lemma \ref{lemma: Properties of energy measures} that
    \begin{align*}
        \Bigl\lvert\Gamma_p\Span{f}(A)^{\frac{1}{p}} - \Gamma_p\Span{f_{n}}(A)^{\frac{1}{p}}\Bigr\rvert  
        \leq \Gamma_p\Span{f-f_n}(A)^{\frac{1}{p}}  
        \leq \mathscr{E}_p(f - f_n)^{\frac{1}{p}}. 
    \end{align*}
    Since $\mathscr{E}_p(f - f_n) \to 0$ as $n \to \infty$, we are done.
\end{proof}

The final objective of the section is Theorem \ref{thm: Energy measure optimal potential} where we derive nice expression for the energy measures. The optimal potential functions $\mathscr{U}_{p,\pm}$ from Section \ref{Sec: Discretizations and mollifiers} have an important role here. We first need the following lemma.

\begin{lemma}\label{lemma: Optimal potential self-similar}
    For all $n \in \N \cup \{0\}$ and $e \in E_n$ 
    \begin{equation}\label{eq: Optimal potential self-similar}
        \mathscr{U}_{p,\pm} \circ F_e = U_{p,\pm,n}(e^-) + (U_{p,\pm,n}(e^+) - U_{p,\pm,n}(e^-))\mathscr{U}_{p,+}.
    \end{equation}
\end{lemma}

\begin{proof}
    Let $m\in \N$, $v \in V_m$ be any vertex and $x \in \Fib(v)$. It follows from the definition of $\mathscr{U}_{p,\pm}$, \eqref{eq: Smooth function on fibers}, and Theorem \ref{thm: Expression of optimal potential} that
    \begin{align*}
        (\mathscr{U}_{p,\pm} \circ F_e)|_{\Fib(v)} & =  \mathscr{U}_{p,\pm}|_{\Fib(e \cdot v)}
        = U_{p,\pm,n+m}(e \cdot v)\\
        & = U_{p,\pm,n}(e^-) + (U_{p,\pm,n}(e^+) - U_{p,\pm,n}(e^-))U_{p,+,m}(v)\\
        & = U_{p,\pm,n}(e^-) + (U_{p,\pm,n}(e^+) - U_{p,\pm,n}(e^-))\mathscr{U}_{p,+}|_{\Fib(v)}.
    \end{align*}
    So the equality of functions \eqref{eq: Optimal potential self-similar} holds on the dense subset $\Fib(X)$. We are done because the functions here are continuous.
\end{proof}

Observe that $\mathscr{U}_{p,-} = 1 - \mathscr{U}_{p,+}$, so it holds that $\mathfrak{m}_p\Span{\mathscr{U}_{p,+}} = \mathfrak{m}_p\Span{\mathscr{U}_{p,-}}$ and $\Gamma_p\Span{\mathscr{U}_{p,+}} = \Gamma_p\Span{\mathscr{U}_{p,-}}$.
Thus we simplify the notation by writing
$\mathfrak{m}_p\Span{\mathscr{U}_{p}} := \mathfrak{m}_p\Span{\mathscr{U}_{p,\pm}} $ and $ \Gamma_p\Span{\mathscr{U}_{p}} := \Gamma_p\Span{\mathscr{U}_{p,\pm}}$.

We recall the definition of Bernoulli measures from Definition \ref{def: Measures on sigma}.

\begin{proposition}\label{prop: pre-energy measure optimal potential}
$\mathfrak{m}_p\Span{\mathscr{U}_p}$ is a Bernoulli measure given by the weights
    \begin{equation}\label{eq: pre-energy optimal (1)}
        \mathfrak{m}_p\Span{\mathscr{U}_p}(\Sigma_{e}) =  \abs{\nabla U_p(e)} \abs{\mathcal{J}_p(e)} = \cM_p^{-1}\abs{\nabla U_p(e)}^p \text{ for all } e \in E_1.
    \end{equation}
    Furthermore, if $f := \mathscr{U}_{p,n}[g]$ for $g : V_n \to \R$, then the pre-energy measure $\mathfrak{m}_p\Span{f}$ satisfies
    \begin{equation}\label{eq: pre-energy optimal (2)}
        \mathfrak{m}_p\Span{f}\restr_{\Sigma_{e}} = \cM_p^{-n}\abs{\nabla g(e)}^p \cdot (\sigma_e)_*(\mathfrak{m}_p\Span{\mathscr{U}_p}) \text{ for all } e \in E_n.
    \end{equation}
    \end{proposition}

\begin{proof}
    It follows from Proposition \ref{prop: Smooth function on fibers}-\eqref{item: Smooth function on fibers (4)} and Lemma \ref{lemma: Properties of energy measures}-\eqref{item: Properties of energy measures (1)} that $\mathfrak{m}_p\Span{\mathscr{U}_p}$ is a probability measure. For $e \in E_n$ we compute using Lemma \ref{lemma: Optimal potential self-similar},
    \begin{align*}
        \mathfrak{m}_p\Span{\mathscr{U}_p}(\Sigma_e) & = \cM_p^{-n} \mathscr{E}_p(\mathscr{U}_p \circ F_e) = \cM_p^{-n} \abs{\nabla U_{p,n}(e)}^p \mathscr{E}_p(\mathscr{U}_{p})\\
        & = \cM_p^{-n} \abs{\nabla U_{p,n}(e)}^p = \abs{\nabla U_{p,n}(e)}\abs{\mathcal{J}_{p,n}(e)}.
    \end{align*}
    The last equality follows from the duality relation \eqref{eq:duality (potentials and flows)}.
    Thus, it follows from Corollary \ref{cor: Gradient product} that $\mathfrak{m}_p\Span{\mathscr{U}_p}$ is the Bernoulli measure with weights \eqref{eq: pre-energy optimal (1)}.

    We move on to \eqref{eq: pre-energy optimal (2)}.
    Let $e \in E_n$ and $e' \in E_{m}$.
    From \eqref{eq: Definition of U_p} we see that
    \begin{align*}
         \mathfrak{m}_p\Span{f}(\Sigma_{e \cdot e'}) & = \cM_p^{-(n+m)} \mathscr{E}_p(f \circ F_{e \cdot e'} )\\
        & = \cM_p^{-(n+m)} \mathscr{E}_p(f \circ F_{e} \circ F_{e'} )\\
        & = \cM_p^{-n}\abs{\nabla g(e)}^p \cM_p^{-m}\mathscr{E}_p(\mathscr{U}_p \circ F_{e'})\\
        & = \cM_p^{-n}\abs{\nabla g(e)}^p (\sigma_{e})_*(\mathfrak{m}_p \Span{\mathscr{U}_p})(\Sigma_{e \cdot e'}) 
    \end{align*}
    Since the sets $\Sigma_{e \cdot e'}$ generate the topology of $\Sigma_e$, the equality \eqref{eq: pre-energy optimal (2)} follows.
\end{proof}

%Measures on self-similar sets that are obtained as push-forwards of Bernoulli measures are typically called \emph{self-similar measures}. According to Proposition \ref{prop: pre-energy measure optimal potential}, the $p$-energy measure of the optimal potential function $\Gamma_p\Span{\mathscr{U}_p} = \chi_*(\mathfrak{m}_p\Span{\mathscr{U}_p})$ is, by definition, a self-similar measure. Such behavior of energy measures does not hold in broader context of self-similar fractals. See Remark \ref{rem: ssmeasures vs EM} for further discussions.

We aim to replace the pre-energy measures with the energy measures in Proposition \ref{prop: pre-energy measure optimal potential}. To this end, we show that the energy measures do not charge fibers.

\begin{proposition}\label{prop: Energy measure of fibers}
    For all $f \in \mathscr{F}_p$ we have $\Gamma_p \Span{f}(\Fib(X)) = 0$. In particular,
    \[
    \Gamma_p\Span{f}(X_e \cap X_{e'}) = 0 \text{ for all distinct edges } e,e'\in E_n \text{ and } n \in \N.
    \]
\end{proposition}

\begin{proof}
    Because the set $V_\#$ is countable, it is sufficient to verify $\Gamma_p \Span{f}(\Fib(v)) = 0$ for all $v \in V_{\#}$ .
    Moreover, we only need to prove this for a dense subset of Sobolev functions thanks to Lemma \ref{lemma: Abs.cont. of energy measures}.
    Thus, by Theorem \ref{thm: Mollifiers converge}, we may assume that $f = \mathscr{U}_{p,n}[g]$ for any $g : V_n \to \R$ and $n \in \N \cup \{0\}$.
    Lastly, it is sufficient to verify the case $f = \mathscr{U}_p$ and $v = v_{\pm}$ due to Lemma \ref{lemma: Optimal potential self-similar} and \eqref{eq: pre-energy optimal (2)}.

    We take $v = v_+$ for simplicity.
    Following these reductions, note that
    \[
        \chi^{-1}(\Fib(v_+)) \subseteq \bigcup_{a \in I^k} \Sigma_{\fe(v_+ \cdot a)} \text{ for all } k \in \N.
    \]
    For any $a \in I^k$, we have
    \begin{align*}
        \mathfrak{m}_p\Span{\mathscr{U}_p} (\Sigma_{\fe(v_+ \cdot a)}) & = \mathfrak{m}_p\Span{\mathscr{U}_p} (\Sigma_{\fe(v_+ \cdot a)}) \sum_{e \in E_1} \mathfrak{m}_p\Span{\mathscr{U}_p}(\Sigma_{e})\\
        & \geq \mathfrak{m}_p\Span{\mathscr{U}_p} (\Sigma_{\fe(v_+ \cdot a)}) \sum_{b \in I} \left(\mathfrak{m}_p\Span{\mathscr{U}_p}(\Sigma_{\fe(\phi_+(b))}) + \mathfrak{m}_p\Span{\mathscr{U}_p}(\Sigma_{\fe(\phi_-(b))})\right)\\
        & = 2\mathfrak{m}_p\Span{\mathscr{U}_p} (\Sigma_{\fe(v_+ \cdot a)}) \sum_{b \in I} \mathfrak{m}_p\Span{\mathscr{U}_p}(\Sigma_{\fe(\phi_+(b))})\\
        & = 2 \sum_{b \in I} \mathfrak{m}_p\Span{\mathscr{U}_p}(\Sigma_{\fe(v_+ \cdot (a \cdot b))}).
    \end{align*}
    The third row follows from the conductive uniform property and the last row from the fact that $\mathfrak{m}_p\Span{\mathscr{U}_p}$ is a Bernoulli measure.
    By iterating the previous estimate,
    \begin{align*}
        \Gamma_p\Span{\mathscr{U}_p}(\Fib(v_+)) \leq \sum_{a \in I^k} \mathfrak{m}_p\Span{\mathscr{U}_p}(\Sigma_{\fe(v_+ \cdot a)})
        \leq \frac{1}{2} \sum_{a \in I^{k-1}} \mathfrak{m}_p\Span{\mathscr{U}_p}(\Sigma_{\fe(v_- \cdot a)}) \leq
        \dots
        \leq \frac{1}{2^k}.
    \end{align*}
    Thus, $\Gamma_p\Span{\mathscr{U}_p}(\Fib(v_+)) = 0$ follows by letting $k \to \infty$.
\end{proof}

\begin{theorem}\label{thm: Energy measure optimal potential}
    For all $f \in \mathscr{F}_{p}$ and $n \in \N$ we have
    \[
      \Gamma_p \Span{f} (X_e) = \cM_p^{-n}\mathscr{E}_p(f \circ F_e) \text{ for all } e \in E_n.  
    \]
    In particular, the energy measure of the optimal potential function satisfy
    \begin{equation}\label{eq: energy measure optimal (1)}
        \Gamma_p\Span{\mathscr{U}_p}(X_{e}) =  \abs{\nabla U_{p,n}(e)} \abs{\mathcal{J}_{p,n}(e)} = \cM_p^{-n}\abs{\nabla U_{p,n}(e)}^p \text{ for all } e \in E_n.
    \end{equation}
    The energy measure of the Sobolev function $f := \mathscr{U}_{p,n}[g]$ where $g : V_n \to \R$ satisfies
    \begin{equation}\label{eq: energy measure optimal (2)}
        \Gamma_p\Span{f}\restr_{X_{e}} = \cM_p^{-n}\abs{\nabla g(e)}^p \cdot (F_e)_*(\Gamma_p\Span{\mathscr{U}_p}) \text{ for all } e \in E_n.
    \end{equation}
\end{theorem}

\begin{proof}
    This directly follows from Propositions \ref{prop: pre-energy measure optimal potential} and \ref{prop: Energy measure of fibers}.
\end{proof}

\section{Comparisons with other energies}\label{sec:otherSob}
This section compares our self-similar $p$-energy form $\mathscr{E}_p$ that was constructed in Section \ref{Sec: Sobolev spaces} to some other frequently considered $p$-energies in the literature. 

\subsection{Equivalence with Korevaar--Schoen}
We first verify that $\mathscr{E}_p$ is always equivalent to the \emph{Korevaar--Schoen} type $p$-energy \cite{korevaar1993sobolev}, where the $L^p$ Besov critical exponent (See \cite[Definition 4.1]{Baudoin_2024}) is equal to $\alpha_p := d_{w,p}/p$.
Let us state it formally.
\begin{theorem}\label{Thm: KS-identify}
    There exists $C \geq 1$ so that for all $f \in L^{p}(X,\mu)$,
\begin{align}\label{eq:KorevaarSchoen.equivalence}
    C^{-1}\mathscr{E}_{p}(f) 
		& \leq \liminf_{r \to 0^+} \int_{X}\kint_{B(x,r)}\frac{\abs{f(x) - f(y)}^{p}}{r^{\pwalk}}\,\mu(dy)\,\mu(dx)\\
        & \leq 
        \sup_{r > 0}\int_{X}\kint_{B(x,r)}\frac{\abs{f(x) - f(y)}^{p}}{r^{\pwalk}}\,\mu(dy)\,\mu(dx)\leq 
		C\mathscr{E}_{p}(f).\nonumber
\end{align}
\end{theorem}

Let us first derive some useful general properties.

\begin{proposition}\label{prop: EM-Poincare}
    There exist constants $A,C \geq 1$ so that for every Sobolev function $f \in \mathscr{F}_p$ and any ball $B := B(x,r) \subseteq X$ with radius $r \in (0,\infty)$ we have
    \begin{equation}\label{eq: EM-Poincare}
        \int_B \abs{f - f_B}^p \, d\mu \leq C r^{\pwalk} \Gamma_p\Span{f}(B(x,A r)).
    \end{equation}
\end{proposition}

\begin{proof}
    We assume that $r < \frac{1}{4} L_*^{-1}$ because otherwise the inequality follows from Theorem \ref{thm: Energy analytic properties}-\eqref{item: Poincare} by making the constant $A$ so large that $B(x,A r) = X$.

    Let $n \in \N$ be so that $\frac{1}{4}L_*^{-(n+1)} \leq r < \frac{1}{4}L_*^{-n}$. By choosing $A := 8L_*$, by Proposition \ref{prop: Geom of LS}-\eqref{item: Stars = balls}, we can choose a vertex $v \in V_n$ satisfying $B \subseteq X_v \subseteq B(x,Ar)$. We can now estimate
    \begin{align*}
        \kint_{B} \abs{f - f_{B}}^p \, d\mu & \lesssim \kint_{B} \abs{f -f_{X_v}}^p \, d\mu && \text{(Equation \eqref{eq: const in PI})}\\
        & \lesssim \kint_{X_v} \abs{f - f_{X_v}}^p \, d\mu && \text{(Ahlfors regularity of $\mu$)}\\
        & \lesssim \cM_p^n \sum_{\substack{e \in E_m\\v \in e}} \cM_p^{-n} \mathscr{E}_p(f \circ F_e) && \text{(Poincar\'e IE \eqref{eq: PrePI (1)})}\\
        & = (L_*^{-n})^{d_{w,p} - Q} \Gamma_p\Span{f}(X_v) && \text{(Proposition \ref{prop: Energy measure of fibers})} \\
        & \lesssim r^{d_{w,p} - Q} \Gamma_p\Span{f}(B(x,Ar)).
    \end{align*}
    The desired inequality \eqref{eq: EM-Poincare} follows by using the Ahlfors regularity of $\mu$,
    \[
        \int_{B} \abs{f - f_{B}}^p \, d\mu \lesssim r^Q \kint_{B} \abs{f - f_{B}}^p \, d\mu \lesssim r^{d_{w,p}} \Gamma_p\Span{f}(B(x,Ar)).
    \]
\end{proof}

We also verify the natural capacity upper bound.

\begin{proposition}\label{prop: EM-Ucap}
    There exist constants $A,C > 1$ so that for every open ball $B(x,R) \subseteq X$ there is a $\varphi \in \mathscr{C}_p$ satisfying 
    \begin{equation}\label{eq: EM-Ucap}
        \varphi|_{B(x,R)} \equiv 1, \, \supp[\varphi] \subseteq B(x,AR) \text{ and } \mathscr{E}_p(\varphi) \leq C\frac{\mu(B(x,R))}{R^{d_{w,p}}}.
    \end{equation}
\end{proposition}

\begin{proof}
    We may assume that $R < \frac{1}{4}L_*^{-1}$. Because otherwise we could choose the constant $A$ so that $B(x,Ar) = X$ and $\varphi \equiv 1$.

    Thus, let $n \in \N$ be so that
    $\frac{1}{4}L_*^{-(n+1)} \leq R < \frac{1}{4}L_*^{-n}$.
    Now it follows from Proposition \ref{prop: Geom of LS}-\eqref{item: Stars = balls} that there is a vertex $v \in V_n$ so that $B \subseteq X_v$.
    We will prove that the Sobolev function $\varphi := \mathscr{U}_{p,n}[1_N] \in \mathscr{C}_p$, where $ 1_{N} : V_n \to \R$ is the characteristic function of the set 
    \[
        N := \{ w \in V_n : d_{G_n}(v,w) \leq 1\},
    \]
    satisfies \eqref{eq: EM-Ucap}.

    It is clear that $\varphi|_{X_v} \equiv 1$. 
    According to Proposition \ref{prop: Geom of LS}-\eqref{item: diams/measures of cells}, we have
    \[
        \supp[\varphi] \subseteq \bigcup_{w \in N} X_w \subseteq B(x,4L_*^{-n}) \subseteq B(x,16L_*R).
    \]
    Thus, we may choose $A := 16L_*$, and what is left is to verify the energy estimate.
    Consider the cut-set of edges $C_N := \{ \{ w,w' \} : w \in N \text{ and } w' \notin N \}$, which clearly satisfies $\abs{C_N} \leq C_{\deg}^2$.
    It follows from \eqref{eq: Definition of U_p} that $\varphi$ is constant on each set $X_e$ for $e \in E_n \setminus C_N$. Thus, we compute
    \begin{align*}
        \mathscr{E}_p(\varphi) & = \cM_p^{-n}\sum_{e \in E_n} \mathscr{E}_p(\varphi \circ F_e) =  \cM_p^{-n}\sum_{e \in C_N} \mathscr{E}_p(\varphi \circ F_e)\\
        & = \cM_p^{-n}\sum_{e \in C_N} \mathscr{E}_p(\mathscr{U}_{p}) = \abs{C_N}\cM_p^{-n}\\
        & \lesssim R^{Q - d_{w,p}} \lesssim \frac{\mu(B(x,R))}{R^{d_{w,p}}}.
    \end{align*}
    In the last row we used $Q$-Ahlfors regularity of $\mu$ (Proposition \ref{prop: Geom of LS}-\eqref{item: Q-AR}).
\end{proof}

\begin{remark}\label{rem:CS}
    The proof of Proposition \ref{prop: EM-Ucap} can be used to establish the cutoff energy condition introduced in \cite{anttila2025approach}.
    Indeed, let $B(x,R) \subseteq X$ and $\varphi$ be as in the proof.
    If $B(y,r) \subseteq X$ is another open ball with $r \leq R$ then we have
    \[
        \Gamma_p\Span{\varphi}(B(y,r)) \leq C \left(\frac{r}{R}\right)^{p\delta_p}\frac{\mu(B(y,r))}{r^{d_{w,p}}}
    \]
    where $\delta_p$ is as in Proposition \ref{prop: Smooth function on fibers}.
    To see this, let $n \in \N$ so that $\frac{1}{4}L_*^{-(n+1)} \leq R < \frac{1}{4}L_*^{-n}$. Choose any $m \in \N$ and $e = e_n \cdot e_m \in E_{n+m}$. Then we apply the formulas \eqref{eq: energy measure optimal (1)} and \eqref{eq: energy measure optimal (2)} in Theorem \ref{thm: Energy measure optimal potential},
    \begin{align*}
        \Gamma_p\Span{\varphi}(X_{e_n e_m}) & = \cM_p^{-n}\abs{\nabla 1_N(e_n)}\cdot \Gamma_p\Span{\mathscr{U}_p}(X_{e_m})\\
        & \leq \cM_p^{-n}\Gamma_p\Span{\mathscr{U}_p}(X_{e_m}) = \cM_p^{-(n+m)}\abs{\nabla U_{p,m}(e_m)}^p\\
        & \leq \left(L_*^{-m}\right)^{p \delta_p} \left(L_*^{-(n+m)}\right)^{Q - d_{w,p}}
    \end{align*}
    The desired energy estimate now follows by covering the ball $B(y,r)$ with suitable sets $X_{e_n e_m}$ using Proposition \ref{prop: Geom of LS}-\eqref{item: diams/measures of cells}.
\end{remark}

\begin{proof}[Proof of Theorem \ref{Thm: KS-identify}] 
The theorem is implied by \cite[Corollaries 1.14 and 1.15]{shimizu2024characterizationssobolevfunctionsbesovtype} since therein assumptions follow from Propositions \ref{prop: Geom of LS}, \ref{prop: EM-Poincare}, \ref{prop: EM-Ucap}, Theorems \ref{Thm: Existence of p-Energy} and \ref{thm: Energy analytic properties} and Lemma \ref{lemma: Properties of energy measures}. For the convenience of the reader, we also sketch a short argument. For any $x \in X$, $r \in (0,\infty)$ and $f \in L^p(X,\mu)$, 
    \begin{align*}
        & \, \quad\int_{B(x,r)}\kint_{B(y,r)}\frac{\abs{f(y) - f(z)}^{p}}{r^{\pwalk}}\,\mu(dz)\,\mu(dy) \\
         \lesssim & \quad \int_{B(x,2r)}\kint_{B(x,2r)}\frac{\abs{f(y) - f(z)}^{p}}{r^{\pwalk}}\,\mu(dz)\,\mu(dy) && \text{(Ahlfors regularity of $\mu$)}\\
        \lesssim &\quad \int_{B(x,2r)}\frac{\abs{f(y) - f_{B(x,2r)}}^{p}}{r^{\pwalk}}\,\mu(dy) \\
        \lesssim &\quad  \Gamma_{p}\Span{f}(B(x,2Ar)) && \text{(Proposition \ref{prop: EM-Poincare}).}
    \end{align*}
    (Here we put $\Gamma_{p}\Span{f}(B(x,2Ar)) = \infty$ if $f \not\in \mathscr{F}_p$.) Taking the summation over a suitable family of points $x$, we get last estimate in \eqref{eq:KorevaarSchoen.equivalence}.

    The middle estimate is trivial so we are left with the first one.
    Now, for $n \in \mathbb{N}$ and $e \in E_{n}$ it follows from the definition of $M_n[\cdot]$ and Hölder's inequality, 
    \[
        \abs{M_{n}[f](e^+) - M_{n}[f](e^-)}^p 
        \lesssim \abs{E_1}^n\int_{X_{e^+}}\kint_{B(x,4L_{*}^{-n})}\abs{f(x) - f(y)}^{p}\,\mu(dy)\,\mu(dx). 
    \]
    We proceed using \eqref{eq: energy measure optimal (2)},
    \begin{align*}
        \Gamma_p\Span{\Xi_{p,n}[f]}(X_e) 
        &\le \cM_p^{-n}\abs{M_{n}[f](e^+) - M_{n}[f](e^-)}^p \\
        &\lesssim \int_{X_{e^+}}\kint_{B(x,4L_{*}^{-n})}\frac{\abs{f(x) - f(y)}^{p}}{L_{*}^{-n\pwalk}}\,\mu(dy)\,\mu(dx). 
    \end{align*}
    Taking the summation over $e \in E_{n}$, we have 
    \begin{equation*}
        \mathscr{E}_{p}(\Xi_{p,n}[f]) 
        \lesssim \int_{X}\kint_{B(x,4L_{*}^{-n})}\frac{\abs{f(x) - f(y)}^{p}}{L_{*}^{-n\pwalk}}\,\mu(dy)\,\mu(dx).
    \end{equation*}
    By letting $n \to \infty$, we get the desired inequality thanks to
    Lemma \ref{lemma: Lower-semicontinuity of energy} because
    $\Xi_{p,n}[f] \to f$ in $L^p(X,\mu)$ according to
    Lemma \ref{lemma: G_n mollifier}.
\end{proof}

We record the following version of the $(p,p)$-Poincar\'e inequality described in terms of Korevaar--Schoen $p$-energy forms. We need it in the next subsection.
\begin{proposition}\label{prop: KS-Poincare}
    There exist constants $A,C \geq 1$ so that for every Sobolev function $f \in \mathscr{F}_p$ and any ball $B := B(x,r) \subseteq X$ with radius $r \in (0,\infty)$ we have
    \begin{equation}\label{eq: KS-Poincare}
        \int_B \abs{f - f_B}^p \, d\mu \leq C r^{\pwalk} \liminf_{\varepsilon \downarrow 0}\int_{B(x,A r)}\kint_{B(y,\varepsilon)}\frac{\abs{f(y) - f(z)}^{p}}{\varepsilon^{\pwalk}}\,\mu(dz)\,\mu(dy).
    \end{equation}
\end{proposition}
\begin{proof}
    This follows from Proposition \ref{prop: EM-Poincare} and \cite[Theorem 3.12]{shimizu2024characterizationssobolevfunctionsbesovtype}.
\end{proof}

\subsection{Newton--Sobolev spaces}
The next objective is study the connection between $\mathscr{F}_p$ and the Newton--Sobolev space $N^{1,p}(X,d,\mu)$; we refer to \cite{bjorn2011nonlinear,HKST} for background. In the statement of the following theorem, we assume a slightly technical condition of quasiconvexity (Definition \ref{def: quasiconvex}) which is used to construct a non-Sobolev Lipschitz function.
It may be that this is not really needed, but for us this assumption is fairly harmless. See Remark \ref{rem: Quasiconvex} for further discussion.

\begin{definition}\label{def: quasiconvex}
    A metric space $(X,d)$ is \emph{quasiconvex} if there is a constant $C \geq 1$ so that for every pair $x,y \in X$ there is a continuous curve $\gamma : [0,1] \to X$ with $\gamma(0) = x, \gamma(1) = y$ and $\len(\gamma) \leq C \cdot d(x,y)$. Here the \emph{length} of $\gamma$ is given by
    \[
       \len(\gamma) := \sup \left\{ \sum_{i = 1}^{N - 1} d(\gamma(t_i),\gamma(t_{i+1})) : 0 = t_1 < t_2 < \dots < t_N = 1   \right\}.
    \]
\end{definition}

\begin{theorem}\label{Thm: NS-equivalence}
    If $d_{w,p} = p$, then $\mathscr{F}_p = N^{1,p}(X,d,\mu)$ with equivalent norms, and there are constants $A,C \geq 1$ so that for all Lipschitz function $f \in \Lip(X,d)$,
\begin{equation}
        \kint_{B(x,r)} \abs{f - f_{B(x,r)}} \, d\mu \leq Cr\left(\kint_{B(x,Ar)} (\Lip f)^p \, d\mu\right)^{1/p}.
\end{equation}
Here $\Lip f$ is the pointwise upper Lipschitz-constant function of $f$,
\begin{equation}\label{eq: Lipconstfcn}
    (\Lip f)(x) \coloneqq \limsup_{\varepsilon \to 0}\sup_{y \in B(x,\varepsilon)}\frac{\abs{f(x) - f(y)}}{\varepsilon}. 
\end{equation}
On the other hand, if $d_{w,p} > p$ and $(X,d)$ is quasiconvex, then there is a Lipschitz function $f \in \Lip(X,d)$ so that $f \notin \mathscr{F}_p$. In particular, $\mathscr{F}_p \neq N^{1,p}(X,d,\mu)$.
\end{theorem}

\begin{proof}
    Note that according to Lemma \ref{lemma: Properties of cond.unif.}-\eqref{item: Properties of cond.unif. (3)}, we always have $d_{w,p} \geq p$.
    Thanks to 
    \cite[Theorems 8.4.2 and 10.5.2]{HKST} 
    we know that if the $L^p$ Besov critical exponent is equal to 1, which in our case means that $d_{w,p} = p$, then the Korevaar--Schoen-Sobolev space coincides with the Newton--Sobolev space with equivalent norms, as long as the following 
    \emph{$p$-Poincar\'e inequality} holds: There are constants $A,C \geq 1$ so that for all Lipchitz functions $f : X \to \R$ and balls $B(x,r) \subseteq X$ with $r > 0$ we have
    \begin{equation}\label{eq: UG-Poincare}
        \kint_{B(x,r)} \abs{f - f_{B(x,r)}} \, d\mu \leq Cr\left(\kint_{B(x,Ar)} (\Lip f)^p \, d\mu\right)^{1/p}.
    \end{equation} 
    Note that $\Lip f < \infty$ since $f$ is Lipschitz.
    
    Let us first assume $\pwalk = p$. We obtain the desired $p$-Poincar\'e inequality \eqref{eq: UG-Poincare} by combining \eqref{eq: KS-Poincare} with the reverse Fatou lemma and H\"{o}lder's inequality. Alternatively, we could have also used the observation in Remark \ref{rem:CK}.
    
    Next assume $\pwalk > p$ and the quasiconvexity of $(X,d)$. Note that $X$ is then biLipschitz equivalent to a geodesic metric space, and that the set of functions $N^{1,p}(X,d,\mu)$ is invariant under biLipschitz changes of metrics. 
    Hence, without the loss of generality, we assume that $X$ itself is geodesic.
    
    To show $N^{1,p}(X,d,\mu) \neq \mathscr{F}_p$, it is sufficient to verify that there is a Lipchitz function $f \notin \mathscr{F}_p$. Fix any point $x_0 \in X$ and let $f(x) := d(x_0,x)$. Clearly $f$ is Lipschitz.
    Let $r > 0$ be sufficiently small and let $x \in X$ such that $d(x,x_0) \geq r$.
    By using the fact that the length of a curve is additive, there is $z \in X$ such that $d(x,z) = r/8$ and $d(x_0,z) = d(x,x_0) - r/8$. Then we estimate
    \begin{align*}
        \kint_{B(x,r)} \frac{\abs{f(x) - f(y)}^p}{r^{d_{w,p}}}\, d\mu(y) & \geq \frac{1}{\mu(B(x,r))} \int_{B(z,r/16)} \frac{\abs{f(x) - f(y)}^p}{r^{d_{w,p}}}\, d\mu(y)\\
        & \geq \frac{\mu(B(z,r/16))}{\mu(B(x,r))}
        16^{-p} r^{p - d_{w,p}} = C r^{p - d_{w,p}}.
    \end{align*}
    The second inequality follows from the triangle inequality and the choices of $x,z$.
    By integrating over all $x \in X \setminus B(x_0,r)$, we get
    \begin{align*}
        & \quad \, \int_X \kint_{B(x,r)} \frac{\abs{f(x) - f(y)}^p}{r^{d_{w,p}}} \, d\mu(y)d\mu(x)\\
        \geq & \quad \int_{X \setminus B(x,r)} \kint_{B(x_0,r)} \frac{\abs{f(x) - f(y)}^p}{r^{d_{w,p}}} \, d\mu(y)d\mu(x)\\
        \geq & \quad C(1 - \mu(B(x_0,r))) r^{p-d_{w,p}}.
    \end{align*}
    Since $d_{w,p} > p$, the last row in the previous inequalities tends to $\infty$ as $r\to 0^+$.
    We now get $f \notin \mathscr{F}_p$ from \eqref{eq:KorevaarSchoen.equivalence}. 
\end{proof}

\begin{remark}\label{rem: Quasiconvex}
    Quasiconvexity of $(X,d)$ is implied by any of the three conditions in Lemma \ref{lemma: Properties of cond.unif.} that were sufficient conditions for the conductively uniform property.
    Let us sketch the main ideas. First, we show that any of these conditions imply
    \[
        \dist\left(I_+^{(n)},I_-^{(n)},d_{G_n}\right) \leq L_*^n.
    \]
    We use this and self-similarity to construct suitable discrete paths on $(X,d)$.
    Then we take a limit using Arzela-Ascoli theorem and obtain a rectifiable curve $\gamma : [0,1] \to X$, which is of suitable length thanks to the visual metric property \eqref{eq: Visual metric}.
\end{remark}

\begin{remark}
    The proof in the case of $d_{w,p} = p$ relies only on \eqref{eq: KS-Poincare} and the fact that $(X,d,\mu)$ is a complete doubling metric measure space. The latter condition is needed to use \cite[Theorem 8.4.2]{HKST}. 
    In particular, by \cite[Proposition 5.28]{kajino2024korevaarschoenpenergyformsassociated}, we can conclude that the Sobolev space $\mathcal{W}^{p}$ defined in \cite{kigami} coincides with the Newton--Sobolev space $N^{1,p}$ when $d_{w,p} = p$. This gives a complete answer for the last question in \cite[Problem 4 in p.~113]{kigami} ($\beta_{p}$ is used in \cite{kigami} instead of $d_{w,p}$). 
\end{remark}

\begin{remark}
    It is possible that $\mathscr{F}_p$ contains Lipschitz functions even in the case $d_{w,p} > p$. In the case of Example \ref{example: Counterexample} we have $\abs{\nabla U_p} \leq L_*^{-1}$.
    By Proposition \ref{prop: Smooth function on fibers}-\eqref{item: Smooth function on fibers (1)} the optimal potential $\mathscr{U}_{p,+}$ is Lipshictz in this case.
\end{remark}

\section{Singularity of Sobolev spaces}\label{sec:SingularitySobolev}
The final part of the paper studies the more detailed properties of the Sobolev spaces.
Let us begin by stating the main theorem, which also includes Theorem \ref{intro:Thm} regarding the Laakso diamond space according to the computations in Figures \ref{intro: potential/flow} and \ref{intro: flow}.

\begin{theorem}[Singularity of Sobolev spaces] 
\label{thm:SingularitySobolev}
Let $p_1,p_2 \in (1,\infty)$ be distinct, and assume that the following three conditions hold.
\begin{enumerate}[label=(\Alph*)]
    \vspace{3pt}
    \item \label{item:A}
    For $p \in \{p_1,p_2\}$ the optimal potential function satisfies
    \begin{equation*}\label{eq: Optimal potential does not degenerate}
        \abs{\nabla U_p(e)} \neq 0 \text{ for all } e \in E_1,
    \end{equation*}
    \item \label{item:B}
    There is $e_1 \in E_1$ so that
    \begin{equation*}
        \abs{\nabla U_{p_1}(e)} \neq \abs{\nabla U_{p_2}(e)}.
    \end{equation*}
    \item \label{item:C}
    The optimal flows satisfy
    \begin{equation*}
        \mathcal{J}_{p_1}(v,w) = \mathcal{J}_{p_2}(v,w) \text{ for all } \{v,w\} \in E_1.
    \end{equation*}
\end{enumerate}
    Then for all distinct $p_1,p_2 \in (1,\infty)$,
    \[
        \mathscr{F}_{p_1} \cap \mathscr{F}_{p_2} = \{ f \in L^1(X,\mu) : \text{$f$ is constant $\mu$-almost everywhere} \}.
    \]
\end{theorem}

\begin{example}\label{ex: Laakso 2}
    The IGS in Figure \ref{fig: Capacity Laakso 2} satisfies the three conditions in Theorem \ref{thm:SingularitySobolev}.
    First, \ref{item:A} clearly holds because $\alpha_p \notin \{0,1/2,1\}$.
    \begin{figure}[!ht]
    \includegraphics[width=200pt]{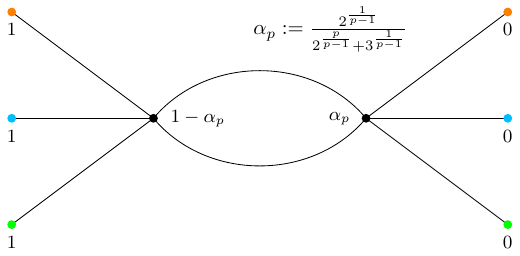}
    \caption{Another example where the singularity of Sobolev spaces hold.}
    \label{fig: Capacity Laakso 2}
    \end{figure}
    The condition \ref{item:B} holds because $\alpha_{p_1} \neq \alpha_{p_2}$. The optimal flows are easily seen to be independent of $p$ since $\abs{\mathcal{J}_p(e)} = 1/2$ for the edges in the middle and $\abs{\mathcal{J}_p(e)} = 1/3$ for the rest. Thus, \ref{item:C} also holds.
\end{example}

\subsection{Weak derivative}
One of the key tools in the remaining proofs is the \emph{weak derivative}, which allows us to perform computations that resemble BV-theory on the real line.
We shall for the rest of the work denote the discretization operator, that was introduced in Section \ref{Sec: Discretizations and mollifiers} and extended in Subsection \ref{subsec: Extending operators}, by $V_{p,n}[\cdot] : \mathscr{F}_p \to \R^{V_n}$. Namely, we emphasize that it depends on $p$.
Similarly, the averaging operators from Section \ref{sec: DPT} are denoted $V_{p,n,m}[\cdot]$.

Before proceeding, we recall the definition of $\Gamma_p\Span{\mathscr{U}_p}$ from Subsection \ref{subsec:EM}.

\begin{definition}
    Assume that the condition \ref{item:A} holds for a given $p \in (1,\infty)$.
    First we define the linear operators $D_{p,n}\Span{\cdot} : \mathscr{F}_p \to L^p(X,\Gamma_p\Span{\mathscr{U}_p})$ for all $n \in \N \cup \{0\}$,
    \begin{equation}\label{eq: Approx gradient}
        D_{p,n} \Span{f} := \sum_{\substack{e \in E_n \\ e =\{ v,w \} }} \frac{\nabla V_{p,n}[f](v,w)}{\nabla U_{p,-,n}(v,w)}  1_{X_e}.
    \end{equation}
    Recall that $\nabla g(v,w) = g(w) - g(v)$. Also we note that the above expression is well-defined thanks to the condition \ref{item:A} and Corollary \ref{cor: Gradient product}.
    Now, we define the \emph{weak derivative} $D_p\Span{\cdot} : \mathscr{F}_p \to L^p(X,\Gamma_p\Span{\mathscr{U}_p})$ as the limit
    \begin{equation}\label{eq:gradient}
        D_p\Span{f} := \lim_{n \to \infty} D_{p,n}\Span{f} \text{ in } L^p(X,\Gamma_p\Span{\mathscr{U}_p}).
    \end{equation}
\end{definition}

We first verify the existence of the limit in \eqref{eq:gradient} and some properties of $D_p\Span{\cdot}$.

\begin{proposition}\label{prop:WD}
    Assume that the condition \ref{item:A} holds for a given $p \in (1,\infty)$. Then $D_p\Span{\cdot}$ is a well-defined bounded linear operator. Furthermore, we have
    \begin{equation}\label{eq: RadonNikodym energy measure}
    \frac{d \Gamma_p\Span{f}}{d\Gamma_p\Span{\mathscr{U}_p}} = \abs{D_{p} \Span{f}}^p \text{ for all } f \in \mathscr{F}_p.
\end{equation}
In particular, $\Gamma_p\Span{f} \ll \Gamma_p\Span{\mathscr{U}_p}$ for all $f \in \mathscr{F}_p$.
\end{proposition}

\begin{proof}
    For this proof, we orient all edges $e \in E_{\#}$ so that $e = (e^+,e^-)$, and also define the \emph{signed gradient} of $g \in \R^{V_n}$ by $\nabla g(e) := g(e^-) - g(e^+)$. Note that $(e_n \cdot e_m)^{\pm} = e_n \cdot e_m^{\pm}$ by \eqref{eq:phieq}.
    
    Theorem \ref{thm: Expression of optimal potential} and the equality $U_{p,-,k} = 1 - U_{p,+,k}$ imply
    \[
        \nabla U_{p,-,n+m}(e_n \cdot e_m) = \nabla U_{p,-,n}(e_n) \nabla U_{p,-,m}(e_m).
    \]
    Similarly using \eqref{eq: Definition of U_p},
    \[
        \nabla V_{p,n+m}[\Psi_{p,n}[f]](e_n \cdot e_m) = \nabla V_{p,n}[f](e_n) \nabla U_{p,-,m}(e_m). 
    \]
    By using Theorem \ref{thm: Energy measure optimal potential} and the condition \ref{item:A},
    \begin{equation}\label{eq: M_p in nice case}
        \cM_p^{-n} = \abs{\nabla U_{p,n}(e) }^{-p} \Gamma_p\Span{\mathscr{U}_p}(X_e) \text{ for all } e \in E_n.
    \end{equation}
    
    Now, we combine these three equalities to compute
    \begin{align*}
    & \mathcal{E}_p^{(n+m)}(\Psi_{p,n+m}[f] - \Psi_{p,n}[f])\\
    = \sum_{\substack{e_n \in E_n \\ e_m \in E_m}} & \abs{\nabla V_{p,n+m}[f](e_n \cdot e_m) - \nabla V_{p,n}[f](e_n) \nabla U_{p,-,m}(e_m) }^p \cM_p^{-n-m} \\
    = \sum_{\substack{e_n \in E_n \\ e_m \in E_m}} & \left| \frac{\nabla V_{p,n+m}[f](e_n \cdot e_m)}{\nabla U_{p,-,n+m}(e_n \cdot e_m) } - \frac{\nabla V_{p,n}[f](e_n)}{\nabla U_{p,-,n}(e_n) } \right|^p \Gamma_p\Span{\mathscr{U}_p}(X_{e_n \cdot e_m})\\
    = \sum_{\substack{e_n \in E_n \\ e_m \in E_m}} & \int_{X_{e_n \cdot e_m}} \abs{D_{p,n+m}\Span{f} - D_{p,n}\Span{f}}^p \, d\Gamma_p \Span{\mathscr{U}_p}\\
    = & \int_X \abs{D_{p,n+m}\Span{f} - D_{p,n}\Span{f}}^p \, d\Gamma_p \Span{\mathscr{U}_p}.
\end{align*}
Since $\Psi_{p,n}[f]$ is a Cauchy sequence in $\mathscr{F}_p$ by Theorem \ref{thm: Mollifiers converge}, it follows from the computation above that
\begin{align*}
    \int_X \abs{D_{p,n+m}\Span{f} - D_{p,n}\Span{f}}^p \, d\Gamma_p \Span{\mathscr{U}_p} & = \mathcal{E}_p^{(n+m)}(\Psi_{p,n+m}[f] - \Psi_{p,n}[f])\\
    & \leq \mathscr{E}_p(\Psi_{p,n+m}[f] - \Psi_{p,n}[f]) \xrightarrow[]{n,m \to \infty} 0.
\end{align*}
By the completeness of $L^p$ spaces, $D_{p,n}\Span{f}$ converges to some $D_p\Span{f}$ in $L^p(X,\Gamma_p\Span{\mathscr{U}_p})$.

We move on to proving \eqref{eq: RadonNikodym energy measure}.
Fix $n \in \N$.
It follows from \eqref{eq: energy measure optimal (1)} and \eqref{eq: energy measure optimal (2)} that for all $m \in \N \cup \{ 0 \}$ and $e_n \in E_n$, $e_m \in E_m$,
\begin{align*}
    \Gamma_p\Span{\Psi_{p,n}[f]}(X_{e_n \cdot e_m}) & = \cM_p^{-n-m} \abs{\nabla V_{p,n}[f](e_n)}^p \mathscr{E}_p(\mathscr{U}_p \circ F_{e_m}) \\
    & = \frac{\abs{\nabla V_{p,n}[f](e_n)}^p}{\abs{\nabla U_{p,n}(e_n)}^p} \Gamma_p\Span{\mathscr{U}_p}(X_{e_n \cdot e_m})\\
    & = \int_{X_{e_n \cdot e_m}}\abs{D_{p,n} \Span{f}}^p \, d\Gamma_p\Span{\mathscr{U}_p}.
\end{align*}
Thus we have
\begin{equation}\label{eq: For WD}
    \Gamma_p\Span{\Psi_{p,n}[f]}(A) = \int_{A}\abs{D_{p,n} \Span{f}}^p \, d\Gamma_p\Span{\mathscr{U}_p}
\end{equation}
where $A = X_e$ for any $e \in E_m$ and $m \geq n$. 
Next assume $e \in E_k$ and $k \leq n$. We prove \eqref{eq: For WD} for $A = X_e$. 
This follows from \eqref{eq: For WD} for $A = X_{e'}$ and $e' \in E_n$, the covering $X_e = \bigcup_{e' \in \pi_{n,k}^{-1}(e)} X_{e'}$ and the fact that $d\Gamma_p\Span{\mathscr{U}_p}$ does not charge the fibers according to Proposition \ref{prop: Energy measure of fibers}:
\begin{align*}
    \Gamma_p\Span{\Psi_{p,n}[f]}(X_e) & = \sum_{e' \in \pi_{n,k}^{-1}(e)}\Gamma_p\Span{\Psi_{p,n}[f]}(X_{e'})\\
    & = \sum_{e' \in \pi_{n,k}^{-1}(e)} \int_{X_{e'}}\abs{D_{p,n} \Span{f}}^p \, d\Gamma_p\Span{\mathscr{U}_p} = \int_{X_{e}}\abs{D_{p,n} \Span{f}}^p \, d\Gamma_p\Span{\mathscr{U}_p}
\end{align*}
So we have \eqref{eq: For WD} for all $A = X_e$ and $e \in E_\#$. But then, \eqref{eq: For WD} extends to all Borel subsets $A \subseteq X$. This follows from Dynkin's $\pi$-$\lambda$ theorem because $\{X_e\}_{e\in E_\#} \cup \{\Fib(v)\}_{v \in V_\#}$ is a $\pi$-system and this family of sets generate the topology of $X$ by Proposition \ref{prop: Geom of LS}.
Finally, by using Lemma \ref{lemma: Abs.cont. of energy measures}, along with the facts that $\Psi_{p,n}[f] \to f$ in $\mathscr{F}_p$ and $D_{p,n} \Span{f} \to D_{p}\Span{f}$ in $L^p(X,\Gamma\Span{\mathscr{U}_p})$,
\begin{align*}
    \Gamma_p\Span{f}(A) & = \lim_{n \to \infty} \Gamma_p\Span{\Psi_{p,n}[f]}(A) = \lim_{n \to \infty} \int_A \abs{D_{p,n} \Span{f}}^p \, d\Gamma_p\Span{\mathscr{U}_p}\\
    & = \int_A  \abs{D_{p} \Span{f}}^p \, d\Gamma_p\Span{\mathscr{U}_p} \text{ for any Borel set } A \subseteq X.
\end{align*}

We shall conclude the proof by showing that $D_p\Span{\cdot} : \mathscr{F}_p \to L^p(X,\Gamma_p\Span{\mathscr{U}_p})$ is a bounded linear operator.
Note that $D_{p,n} \Span{\cdot} : \mathscr{F}_p \to L^p(X,\Gamma_p\Span{\mathscr{U}_p})$ is clearly linear by \eqref{eq: Approx gradient}. From Proposition \ref{prop: Properties of mollifier}-\eqref{item: Properties of mollifier (1)} and  \eqref{eq: For WD} we see that
\[
    \int_X \abs{D_{p,n}\Span{f}}^p \, d\Gamma_p\Span{\mathscr{U}_p} = \mathscr{E}_p(\Psi_{p,n}[f]) = \mathcal{E}_p^{(n)}(f) \leq \mathscr{E}_p(f).
\]
Hence the limit $D_p \Span{\cdot}$ is also bounded and linear.
\end{proof}

The key objects we base on the proof of Theorem \ref{thm:SingularitySobolev} are the measures
\begin{equation}\label{eq:BV}
    \mu_{p,n}\Span{f}(A) := \int_A \abs{D_{p,n} \Span{f}} \,d\Gamma_p\Span{\mathscr{U}_p} \text{ and } \mu_{p}\Span{f}(A) := \int_A \abs{D_{p} \Span{f}} \,d\Gamma_p\Span{\mathscr{U}_p}
\end{equation}
for $f \in \mathscr{F}_p$. 
The measure $\mu_p\Span{f}$ can be understood as a variant of the variation on the real line, $A \mapsto \int_A \abs{f'(x)}\, dx$ where $A \subseteq \R$. This is suggested by \eqref{eq:LooksBV}.

\begin{proposition}\label{prop:EqualLittlemuu}
    Assume that \ref{item:A} and \ref{item:C} hold for $p_1,p_2 \in (1,\infty)$, and let $f \in \mathscr{F}_{p_1} \cap \mathscr{F}_{p_2}$. Then the measures $\mu_{p_1}\Span{f}$ and $\mu_{p_2}\Span{f}$ are equal.
\end{proposition}

\begin{proof}
    The first step is to show that
    \begin{equation}\label{eq: Discretization ind. of p}
        V_{p_1,m}[f] = V_{p_2,m}[f] \text{ for all } m \in \N \cup \{0\}.
    \end{equation}
    Roughly speaking this follows from \ref{item:C} because the discretization operators defined in Section \ref{Sec: Discretizations and mollifiers} depend only on the flow. However, the actual proof is much more delicate because here $V_{p_i,m}[\cdot]$ are the extended operators, and the extensions for different $p$ are defined with respect to different Sobolev norms.
    To this end, we use a similar argument as in the proof of Corollary \ref{cor: Self-similarity of domain}.

    First, by the condition \ref{item:C} it follows that the averaging operators $V_{p_i,n,m}[\cdot]$ given in \eqref{eq: Discrete average} are the same for $i = 1,2$. We then use the tower rule \eqref{eq: Tower rule},
    \begin{equation}\label{eq: V_m agree}
        V_{p_1,m}[\Xi_{p_1,n}[f]] = V_{p_1,n,m}[M_n[f]] =  V_{p_2,n,m}[M_n[f]] = V_{p_2,m}[\Xi_{p_2,n}[f]]
    \end{equation}
    for all $n \geq m$. Here $\Xi_{p,n}[\cdot]$ are the rough mollifiers from Definition \ref{def: G_n operators}.
    Now, it follows from Lemma \ref{lemma: G_n mollifier} that $\{ \Xi_{p_i,n}[f] \}_{n = m}^{\infty}$ is a bounded sequence in $\mathscr{F}_{p_i}$, and that $\Xi_{p_i,n}[f] \to f$ in $L^{p_i}(X,\mu)$.
    Then Proposition \ref{prop: Equivalent change} implies that the limits
    \[
        \lim_{n \to \infty} V_{p_i,m}[\Xi_{p_i,n}[f]]
    \]
    exist in $\R^{V_m}$ and, thanks to \eqref{eq: V_m agree}, coincide for $i = 1,2$. Finally, to conclude that these limits are equal to $V_{p_i,m}[f]$, we use a similar argument based on the reflexivity of $\mathscr{F}_{p_i}$ as in the proof of Corollary \ref{cor: Self-similarity of domain}.
    Namely, we take a convergent sequence of convex combinations and apply the continuity of $V_{p_i,m}[\cdot]$.
    Thus \eqref{eq: Discretization ind. of p} follows.
    
    We are now ready to show the equality of measures $\mu_{p_1}\Span{f} = \mu_{p_2}\Span{f}$. Let us first fix $n \in \N$ and $e \in E_n$. Then by definitions \eqref{eq: Approx gradient} and \eqref{eq:BV},
    \begin{equation}\label{eq:LooksBV}
        \mu_{p_i,n}\Span{f}(X_e) = \frac{\abs{\nabla V_{p_i,n}[f](e) }}{\abs{\nabla U_{p_i,n}(e)}}\Gamma_{p_i}\Span{\mathscr{U}_{p_i}}(X_e) = \abs{\nabla V_{p_i,n}[f](e) }\abs{\mathcal{J}_{p_i,n}(e)}.
    \end{equation}
    Thus by \eqref{eq: Discretization ind. of p} and \ref{item:C}, $\mu_{p_1,n}\Span{f}(X_e) = \mu_{p_2,n}\Span{f}(X_e)$ for all $e \in E_n$.
    This equality extends to all $e \in E_m$ and $m \leq n$. To see this, we use the covering $X_e = \bigcup_{e' \in \pi_{n,m}^{-1}(e)} X_{e'}$ and the fact that $\mu_{p_i}\Span{f}$ do not charge $\Fib(X)$ according to Propositions \ref{prop: Energy measure of fibers} and \ref{prop:WD}. These facts along with Proposition \ref{prop: X_e cap X_f} give
    \[
        \mu_{p_1,n}\Span{f}(X_e) = \sum_{e' \in \pi_{n,m}^{-1}(e)}\mu_{p_1,n}\Span{f}(X_{e'}) = \sum_{e' \in \pi_{n,m}^{-1}(e)}\mu_{p_2,n}\Span{f}(X_{e'}) = \mu_{p_2,n}\Span{f}(X_e).
    \]
    By letting $n \to \infty$, we get
    \[
        \mu_{p_1}\Span{f}(A) = \mu_{p_2}\Span{f}(A)
    \]
    when $A = X_e$ for $e \in E_\#$. This equality extends to all Borel sets $A \subseteq X$ by the argument in Proposition \ref{prop:WD}.
\end{proof}

We are almost ready to prove Theorem \ref{thm:SingularitySobolev}.
Before this, we discuss an certain dichotomy arising from the condition \ref{item:A}. This is unrelated to Theorem \ref{thm:SingularitySobolev} but definitely worth recording. We discuss its relation to the attainment problem of Ahlfors regular conformal dimension in the remark below.

\begin{theorem}\label{thm:EMdichotomy}
    We have the following dichotomy.
    \begin{enumerate}
        \vspace{2pt}
        \item \label{item:Attain}
        If the condition \ref{item:A} holds for a given $p \in (1,\infty)$ then $\Gamma_p\Span{\mathscr{U}_p}$ is a doubling measure satisfying $\Gamma_p\Span{f} \ll \Gamma_p\Span{\mathscr{U}_p}$ for all $f \in \mathscr{F}_p$.
        \vspace{2pt}
        \item \label{item:NotAttain}
        If the condition \ref{item:A} fails for a given $p \in (1,\infty)$ then there is a fixed $\sigma$-porous subset $\widetilde{A} \subseteq X$ so that $\Gamma_p\Span{f}(X) = \Gamma_p\Span{f}(\widetilde{A})$ for all $f \in \mathscr{F}_p$. In particular, there is no doubling measure $\nu$ satisfying $\Gamma_p\Span{f} \ll \nu$ for all $f \in \mathscr{F}_p$. 
    \end{enumerate}
\end{theorem}

\begin{proof}
    \eqref{item:Attain}: According to Proposition \ref{prop:WD}, we only need to check that $\Gamma_p\Span{\mathscr{U}_p}$ is doubling.
    Our first step is to prove the following claim by induction on $n \in \N$: There is a constant $C \geq 1$ so that for all $n \in \N$ and $e,e' \in E_n$ so that $e \cap e' \neq \emptyset$,
    \begin{equation}\label{eq: Cellular doubling}
        \Gamma_p\Span{\mathscr{U}_p}(X_e) \leq C \Gamma_p\Span{\mathscr{U}_p}(X_{e'}).
    \end{equation}
    To this end, let us fix the constants
    \[
        a := \max_{e \in E_1} \Gamma_p\Span{\mathscr{U}_p}(X_e) \text{ and } b := \min_{e \in E_1} \Gamma_p\Span{\mathscr{U}_p}(X_e),
    \]
    and set $C := a/b$. Note that $C$ is well-defined thanks to \ref{item:A} and \eqref{eq: energy measure optimal (1)}.
    The base case $n = 1$ is trivial, so we assume the claim holds for $n \in \N$. Let $e,e' \in E_{n+1}$ so that $e \cap e' \neq \emptyset$. The case $\pi_{n+1}(e) = \pi_{n+1}(e')$ is clear since the pre-energy measure $\mathfrak{m}_p\Span{\mathscr{U}_p}$ is a Bernoulli measure. Thus, assume $\pi_{n+1}(e) = e_n \neq e_n'= \pi_{n+1}(e')$, and write $e = e_n \cdot e_1$ and $e' = e_n' \cdot e_1'$. By \ref{SM2} and the gluing rules, we have $e_n \cap e_n' \neq \emptyset$ and $\{ e_1,e_1' \} = \{ \fe(\phi_-(a)), \fe(\phi_+(a)) \}$ for some $a \in I$. By \eqref{eq:duality (potentials and flows)}, the conductive uniform property and \eqref{eq: energy measure optimal (1)}, we have $\Gamma_p\Span{\mathscr{U}_p}(X_{e_1}) = \Gamma_p\Span{\mathscr{U}_p}(X_{e_1'})$. We thus compute
    \begin{align*}
        \Gamma_p\Span{\mathscr{U}_p}(X_e) & = \Gamma_p\Span{\mathscr{U}_p}(X_{e_n}) \Gamma_p\Span{\mathscr{U}_p}(X_{e_1}) \\
        & \stackrel{\text{IH}}{\leq} C \Gamma_p\Span{\mathscr{U}_p}(X_{e_n'})\Gamma_p\Span{\mathscr{U}_p}(X_{e_1})  \\
        & = C \Gamma_p\Span{\mathscr{U}_p}(X_{e_n'})\Gamma_p\Span{\mathscr{U}_p}(X_{e_1'}) \\
        & = C\Gamma_p\Span{\mathscr{U}_p}(X_{e'}).
    \end{align*}
    This concludes the proof of \eqref{eq: Cellular doubling}.

    The final step in proving the doubling property is to use \eqref{eq: Cellular doubling} to derive estimates for measures of balls.
    This can be done using a similar argument as in \cite[Lemma 3.31]{anttila2024constructions}, which regards the Ahlfors regularity of the limit space.
    The argument is also sketched in the end of the proof of Proposition \ref{prop: Geom of LS}.
    See also \cite[Theorem 3.3.4]{Kigamiweighted}.

    \eqref{item:NotAttain}: Next we assume that \ref{item:A} fails. 
    Fix any edge $\tilde{e} \in E_1$ so that  $\abs{\nabla U_p(\tilde{e})} = 0$, and consider the subset
    \[
        A := X \setminus \left(\bigcup_{e \in E_\#} F_e(X_{\tilde{e}})\right).
    \]
    In \cite[Lemma 6.4]{anttila2024constructions} it was shown that $A$ as given above is
    porous subset of the limit space in Example \ref{example: Counterexample}. Here we use similar idea to show that
    \[
        A_n := \bigcup_{k = 0}^n \bigcup_{e \in E_k} F_e(A) \text{ for } n \in \N
    \]
    are porous subsets of $X$. To see this, note that for any $e \in E_n$ it follows from the definition of $A$ that $A_n \cap F_e(X_{\tilde{e}}) = \emptyset$. To see that this implies porosity, we use Proposition \ref{prop: Geom of LS}-\eqref{item: Stars = balls} to show that there is an open ball $B(x,r) \subseteq F_e(X_{\tilde{e}})$ where $r$ is comparable to the diameter of $X_{\tilde{e}}$. So $A_n$ are porous for all $n \in \N$, and consequently
    \[
        \widetilde{A} := \bigcup_{n = 0}^\infty A_n
    \]
    is a $\sigma$-porous subset.
    
    It remains to show that $\Gamma_p\Span{f}(X) = \Gamma_p\Span{f}(\widetilde{A})$ for all $f \in \mathscr{F}_p$. By Lemma \ref{lemma: Abs.cont. of energy measures} it is sufficient to verify this for a dense set of Sobolev functions. Thus, by Theorem \ref{thm: Mollifiers converge}, we may assume that $f = \mathscr{U}_{p,n}[g]$.

    Since we chose $\tilde{e}$ so that $\abs{\nabla U_p(\tilde{e})} = 0$, we have by \eqref{eq: energy measure optimal (1)} in Theorem \ref{thm: Energy measure optimal potential} that $\Gamma_p\Span{\mathscr{U}_p}(X) = \Gamma_p\Span{\mathscr{U}_p}(A)$. Thus $\Gamma_p\Span{f}(X) = \Gamma_p\Span{f}(A_n)$ according to \eqref{eq: energy measure optimal (2)} and the definition of $A_n$.
\end{proof}

\begin{remark}\label{rem:Confdim}
    In the current context of Laakso-type fractal spaces, the condition \ref{item:A} is closely related to a deep problem in quasiconformal geometry, the \emph{attainment problem of Ahlfors regular conformal dimension}; an interested reader may survey \cite[Introduction]{anttila2024constructions} for background and references. It was verified in a recent work by first two authors and Rainio \cite{RainioAttainemnt} that, within the current framework and the symmetry condition \ref{item: Sufficient conditions (Symmetry)}, positive answer to this attainment problem is equivalent to the condition \ref{item:A} when $p \in (1,\infty)$ is equal to the Ahlfors regular conformal dimension.
    In what follows, \ref{item:A} characterizes both the attainment problem and the problem of finding a (harmonic function) $h \in \mathscr{F}_p$ so that $\Gamma_p\Span{h}$ is a doubling measure with $\Gamma_p\Span{f} \ll \Gamma_p\Span{h}$ for all $f \in \mathscr{F}_p$. The latter characterization follows from Theorem \ref{thm:EMdichotomy} because doubling measures do not charge porous subsets. This positively resolves \cite[Conjecture 10.8]{murugan2023first} in our framework of Laakso-type fractals.
    While the discussion here only regards the analysis where $p$ is equal to the Ahlfors regular conformal dimension, these observation nevertheless seems to suggest that the Sobolev spaces $\mathscr{F}_p$, for all $p \in (1,\infty)$, capture deeper geometric features of the underlying fractals.
    For instance, when $p=2$, we expect the dichotomy of Theorem \ref{thm:EMdichotomy} to play a role for a similar question, the \emph{attainment of conformal walk dimension}, introduced and studied by Kajino and Murugan \cite{MN}.
\end{remark}

\begin{example}
It is worth to note that whether the condition \ref{item:A} holds or not can depend on $p$. For the IGSs in Figure \ref{fig: degenerate at p=2}, \ref{item:A} holds if and only if $p \neq 2$. 
To see this, the functions in the figure are $2$-harmonic on their respective graphs, and therefore equal to the 2-energy minimizer $U_2$. Clearly $\abs{\nabla U_2}$ vanishes on every vertical edge.
On the other hand, for every other $p$, it is easy to check using $p$-harmonicity that \ref{item:A} holds.
\begin{figure}[!ht]
    \includegraphics[width=315pt]{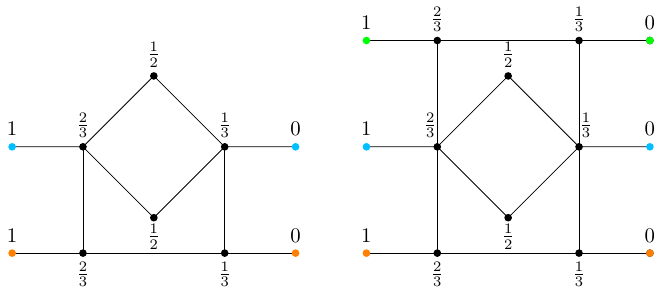}
    \caption{Figures of IGSs for which the condition \ref{item:A} holds if and only if $p \neq 2$. The functions in the figure are the 2-energy minimizers $U_{2,+}$. }
    \label{fig: degenerate at p=2}
\end{figure}
The Ahlfors regular conformal dimension of the limit space produced by the IGS in the right is equal to 2 because $\cM_2 = 1$; see Remark \ref{rem:Confdim} and \cite{RainioAttainemnt} for details. In particular, this space does not attain its Ahlfors regular conformal dimension. The space produced by the IGS in the left does attain because its Ahlfors regular conformal dimension is strictly less than 2 for which $\abs{\nabla U_p}$ does not vanish.
\end{example}

\subsection{Singularity of energy measures}
A key step in the proof of Theorem \ref{thm:SingularitySobolev} is to pass the problem of singularity of Sobolev spaces to energy measures.
This approach requires the following result, which we appropriately term \emph{singularity of energy measures}. This kind of behavior was first discovered to hold on the Sierpi\'nski gasket and some other strongly symmetric post-critically finite fractals in a recent work of Kajino and third author \cite{KajinoSGsurvey,KajinoSingularityPCF}.

\begin{proposition}[Singularity of energy measures]\label{prop:SingularityEM}
    Let $p_1,p_2 \in (1,\infty)$ and assume that the measures $\Gamma_{p_1}\Span{\mathscr{U}_{p_1}}$ and $\Gamma_{p_2}\Span{\mathscr{U}_{p_2}}$ are not equal. Then for any pair $f_1 \in \mathscr{F}_{p_1}$ and $f_2 \in \mathscr{F}_{p_2}$ the measures $\Gamma_{p_1}\Span{f_1}$ and $\Gamma_{p_2}\Span{f_2}$ are mutually singular.
\end{proposition}

\begin{proof}
    Recall that, according to Proposition \ref{prop: pre-energy measure optimal potential}, the pre-energy measure $\mathfrak{m}_p\Span{\mathscr{U}_p}$ is a Bernoulli measure given by the weights $\abs{\nabla U_p(e)} \abs{\mathcal{J}_p(e)}$. Because the measures $\Gamma_{p_1}\Span{\mathscr{U}_{p_1}}$ and $\Gamma_{p_2}\Span{\mathscr{U}_{p_2}}$ are not equal, there has to be an edge where the weights differ. 
    Now, the rest of the proof is an adaptation of the classical Law of large numbers argument which  shows that Bernoulli measures with different weights are mutually singular.
    This method was used in the note of Schioppa \cite{AndreaSchioppa2015} in a somewhat similar context.
    We nevertheless go through the details.

    During this argument, it is convenient to identify $\Sigma$ with the space of sequences $E_1^{\N}$ as it was discussed in Remark \ref{rem: Alternative symbol space}.
    We may further reduce it to the subset $\Sigma_0 := \Sigma \setminus \chi^{-1}(\Fib(X)) \subseteq \Sigma$ because the energy measures do not charge the fibers according to Proposition \ref{prop: Energy measure of fibers}.
    Let us now define the following Borel subsets for both $p_1$ and $p_2$,
    \[
        A_i := \left\{ (e_j)_{j = 1}^{\infty} \in \Sigma_0 : \lim_{n \to \infty} \frac{\abs{ \{ 1 \leq j \leq n : e_j = e \}}}{n} = \mathfrak{m}_{p_i}\Span{\mathscr{U}_{p_i}}(\Sigma_e) \text{ for all } e \in E_1 \right\}.
    \]
    Because the weights of the two Bernoulli measures are not equal, $A_1 \cap A_2 = \emptyset$. Then by noting that $\chi|_{\Sigma_0}$ is injective according to Proposition \ref{prop: X_e cap X_f}, the two Borel subsets $\widetilde{A}_i := \chi(A_i) \subseteq X$ also have an empty intersection.
    In what follows, it is sufficient to verify that $\Gamma_{p_i}\Span{f_i}(X) = \Gamma_{p_i}\Span{f_i}(\widetilde{A}_i)$ for all $f_i \in \mathscr{F}_{p_i}$.
    Yet again, it is by Lemma \ref{lemma: Abs.cont. of energy measures} sufficient to verify this equality for a dense subset of Sobolev functions, so we assume $f_{i} = \mathscr{U}_{p_i,n}[g]$.

    We first use the Strong law of large numbers to get
    \[
    \Gamma_{p_i}\Span{\mathscr{U}_{p_i}}(X) = \Gamma_{p_i}\Span{\mathscr{U}_{p_i}}(\widetilde{A}_i).
    \]
    We refer to \cite{ash2000probability} for details.
    Second, we note that $F_e(\widetilde{A}_{i}) \subseteq \widetilde{A}_{i}$ for all $e \in E_{\#}$. To see this, it clearly follows from the definition of $A_{i}$ that $\sigma_{e}(A_{i}) \subseteq A_{i}$. So we have
    \[
        F_e(\widetilde{A}_{i}) = F_e(\chi(A_{i})) = \chi(\sigma_e(A_{i})) \subseteq \chi(A_{i}) = \widetilde{A}_{i}.
    \]
    By combining the previous two observations with the formula \eqref{eq: energy measure optimal (2)} in Theorem \ref{thm: Energy measure optimal potential} and $f_{p_i} = \mathscr{U}_{p_i,n}[g]$, we have arrived to the desired conclusion.
\end{proof}

On an unrelated note to Theorem \ref{thm:SingularitySobolev}, we also have the following singularity result regarding the measure $\mu$. See \cite{kajino2020singularity,hino2005singularity} for similar results on Dirichlet forms.

\begin{proposition}\label{prop: Energy measure vs mu}
    For all $p \in (1,\infty)$ the following hold.
    \begin{enumerate}
    \vspace{2pt}
        \item \label{item: Energy measure vs mu (1)}
        If $d_{w,p} = p$ then $\Gamma_p\Span{\mathscr{U}_p} = \mu$ and $\Gamma_p\Span{f} \ll \mu$ for all $f \in \mathscr{F}_p$.
        \vspace{2pt}
        \item \label{item: Energy measure vs mu (2)}
        If $d_{w,p} > p$ then $\Gamma_{p}\Span{f}$ and $\mu$ are mutually singular for all $f \in \mathscr{F}_p$.
    \end{enumerate}
\end{proposition}

\begin{proof}
    We recall from Lemma \ref{lemma: Properties of cond.unif.}-\eqref{item: Properties of cond.unif. (3)} that the condition $d_{w,p} = p$ is equivalent to $\abs{\nabla U_p} \equiv L_*^{-1}$ and $\cM_p = \abs{E_1}/L_*^p$. Thus, the measures $\Gamma_p\Span{\mathscr{U}_p}$ and $\mu$ are equal if and only if $d_{w,p} = p$.

    \eqref{item: Energy measure vs mu (1)}: This now follows from the previous observation and Proposition \ref{prop:WD}.

    \eqref{item: Energy measure vs mu (2)}: We apply the Law of large numbers as in the proof of Proposition \ref{prop:SingularityEM}.
\end{proof}

\subsection{Singularity of Sobolev spaces}\label{subsec:singSob}
We finally have all the ingredients to prove our theorem regarding singularity of Sobolev spaces.

\begin{proof}[Proof of Theorem \ref{thm:SingularitySobolev}]
    Let $p_1,p_2 \in (1,\infty)$ be distinct. Note that constant functions are contained in the Sobolev spaces by construction.

    Now, let $f \in \mathscr{F}_{p_1} \cap \mathscr{F}_{p_2}$. To prove that $f$ is a constant function, we show that $\mathscr{E}_{p_i}(f) = \Gamma_{p_i}\Span{f}(X) =  0$. This is sufficient according to Theorem \ref{thm: Energy analytic properties}-\eqref{item: Poincare}. As it follows from Proposition \ref{prop:WD}, we may as well show that $\mu_{p_i}\Span{f}(X) = 0$.

    First, we note that the optimal flow $\mathcal{J}_p$ does not vanish on any edge. This follows from \ref{item:A} and the duality, Lemma \ref{lemma:duality}. Now, by combining all the conditions \ref{item:A}, \ref{item:B} and \ref{item:C}, there in particular has to be an edge $\tilde{e} \in E_1$ so that $\abs{\nabla U_{p_1}(\tilde{e})} \abs{\mathcal{J}_{p_1}(\tilde{e})} \neq \abs{\nabla U_{p_2}(\tilde{e})} \abs{\mathcal{J}_{p_2}(\tilde{e})}$. In what follows, according to Theorem \ref{thm: Energy measure optimal potential}, $\Gamma_{p_1}\Span{\mathscr{U}_{p_1}} \neq \Gamma_{p_2}\Span{\mathscr{U}_{p_2}}$.
    Then Proposition \ref{prop:SingularityEM} implies that $\Gamma_{p_1}\Span{f}$ and $\Gamma_{p_2}\Span{f}$ are mutually singular. Since $\mu_{p_i}\Span{f}$ and $\Gamma_{p_i}\Span{f}$ are mutually absolutely continuous by Proposition \ref{prop:WD}, it follows that $\mu_{p_1}\Span{f}$ and $\mu_{p_2}\Span{f}$ are mutually singular as well.

    On the other hand, the measures $\mu_{p_1}\Span{f}$ and $\mu_{p_2}\Span{f}$ were shown to be equal in Proposition \ref{prop:EqualLittlemuu}. So the measures $\mu_{p_1}\Span{f}$ and $\mu_{p_2}\Span{f}$ are simultaneously equal and mutually singular. The only possibility that allows this is $\mu_{p_i}\Span{f}(X) = 0$.
\end{proof}

The remainder of the work is dedicated to discussing singularity of Sobolev spaces in a broader context. Because this does not happen in $W^{1,p}(\R^n)$ or in the standard settings studied in analysis on metric spaces \cite{HKST}, we focus on fractal spaces. 
But before this, we take a quick side-tour to Lipschitz functions.

\subsection{Lipschitz functions}
Using the same method as in the proof of Theorem \ref{thm:SingularitySobolev}, we prove that some of these Sobolev spaces do not contain any non-constant Lipschitz functions. This generalizes Theorem \ref{intro:NoLips}.

\begin{theorem}\label{thm:NoLips}
    Assume that for $p \in (1,\infty)$ the conditions \ref{item:A}, $d_{w,p} > p$ and 
    \begin{equation}\label{eq:ForLip}
        \sum_{e \in E_1} \abs{\mathcal{J}_p(e)} = L_*
    \end{equation}
    all hold. Then $\mathscr{F}_p$ does not contain any non-constant Lipschitz functions, i.e.,
    \[
        \mathscr{F}_p \cap \Lip(X,d) = \{ f \in L^1(X,\mu) : \text{$f$ is constant $\mu$-almost everywhere} \}.
    \]
\end{theorem}

\begin{proof}
    Fix a $K$-Lipschitz function $f \in \mathscr{F}_p \cap \Lip(X,d)$. We will prove that $f$ is constant by showing that $\mu_p\Span{f}(X) = 0$. This is sufficient by Theorem \ref{thm: Energy analytic properties}-\eqref{item: Poincare}.

    We need some preparations first.
    Let $\mathfrak{m}$ be the Bernoulli measure on $\Sigma$ given by the weights $\alpha(e) := L_*^{-1}\abs{\mathcal{J}_p(e)}$, and $\nu$ be the push-forward measure $\nu := \chi_*(\mathfrak{m})$ on $X$.
    Note that the sum of these weights is 1 due to \eqref{eq:ForLip}.
    By the same argument as in the proof of Proposition \ref{prop: Geom of LS}-\eqref{item: diams/measures of cells}, we see that $\nu$ does not charge the fibers. In particular,
    \begin{equation}\label{eq:computenu}
        \nu(X_e) = L_*^{-n}\abs{\mathcal{J}_{p,n}(e)} \text{ for all } e \in E_\#.
    \end{equation}

    The rest of the argument is similar to the proof of Theorem \ref{thm:SingularitySobolev}. Namely, we show that $\mu_p\Span{f}$ is simultaneously absolutely continuous and mutually singular with $\nu$.
    
    First, the singularity follows from the Law of large numbers argument in the proof of Proposition \ref{prop:SingularityEM} once we know that $\Gamma_p\Span{\mathscr{U}_p} \neq \nu$.
    We use $d_{w,p} > p$ to prove this. According to Lemma \ref{lemma: Properties of cond.unif.}-\eqref{item: Properties of cond.unif. (3)} there is $\tilde{e} \in E_1$ with $\abs{\nabla U_p(\tilde{e})} \neq L_*^{-1}$. Then it follows from Proposition \ref{prop: pre-energy measure optimal potential} and $\abs{\mathcal{J}_p(\tilde{e})} \neq 0$ that the weights determining $\nu$ and $\Gamma_p\Span{\mathscr{U}_p}$ differ at $\tilde{e}$. Hence, $\nu \neq \Gamma_p\Span{\mathscr{U}_p}$, and therefore the singularity follows.

    Next, we use the fact that $f$ is Lipschitz to show that $\mu_p\Span{f} \ll \nu$.
    First we analyze the measures $\mu_{p,n}\Span{f}$. Because $f$ is continuous, its discretizations $V_{p,n}[f]$ are given by the integrals in Definition \ref{def: Discretization operator}. Also note that $\diam(X_e) = L_*^{-n}$ for any $e \in E_n$ according to Proposition \ref{prop: Geom of LS}-\eqref{item: diams/measures of cells}. We then compute
    \[
        \abs{\nabla V_{p,n}[f](e)} \leq \int_{\Fib(e^+)} \int_{\Fib(e^-)}\abs{f(x) - f(y)} \, d\mathcal{J}_{p,e^-}(x) d\mathcal{J}_{p,e^+}(y) \leq K L_*^{-n}.
    \]
    In the last inequality we used that $f$ is $K$-Lipschitz and that $\mathcal{J}_{p,v}$ are probability measures.
    By the computation in \eqref{eq:LooksBV} and \eqref{eq:computenu},
    \[
        \mu_{p,n}\Span{f}(X_e) \leq K L_*^{-n}\abs{\mathcal{J}_{p,n}(e)} = K\nu(X_e) \text{ for all } e \in E_n.
    \]
    Since neither of these measures charge the fibers, the inequality
    \begin{equation}\label{eq:nuIneq}
        \mu_{p}\Span{f}(A) \leq K\nu(A)
    \end{equation}
    holds for all $A=X_e$ where $e\in E_\#$.
    This follows from the argument in the proof of Proposition \ref{prop:EqualLittlemuu}, which is used to prove the equality $\mu_{p_1}\Span{f} = \mu_{p_2}\Span{f}$.
    By again using the fact that these measures do not charge the fibers, we see that \eqref{eq:nuIneq} holds when $A$ is a countable union of sets in $\{ X_e\}_{e \in E_{\#}}$.
    In particular, \eqref{eq:nuIneq} holds for all open sets.
    By using the fact that $\mu_p\Span{f}$ and $\nu$ are Radon measures, we get \eqref{eq:nuIneq} for all Borel sets.
    Hence $\mu_p\Span{f} \ll \nu$, and this proves the claim.

    %The final step is to verify \eqref{eq:nuIneq} for a suitable family of Borel sets $\mathcal{A}$ and apply Caratheodory's extension theorem. We set
    %\begin{align*} \mathcal{A}_0 := \{ X_e \}_{e \in E_\#} & \cup \{ \Fib(v) \}_{v \in V_\#} \cup \{ X_e \setminus \Fib(v) \}_{v \in e \in E_\#}\\& \cup \{ X_e \setminus (\Fib(e^+) \cup \Fib(e^-))\}_{e \in E_\#}\end{align*}
    %and choose $\mathcal{A}$ to be the family of Borel sets that can be expressed as finite unions of sets in $\mathcal{A}$. Using one last time that fiber's are not charged, we see that \eqref{eq:nuIneq} holds for sets in $\mathcal{A}$. We also see from Proposition \ref{prop: X_e cap X_f} that $\mathcal{A}$ is an algebra of sets, so we can use Caratheodory's extension theorem. By applying it to $K\nu - \mu_p\Span{f}$ it follows that $\mu_p\Span{f} \leq K\nu$, and in particular $\mu_p\Span{f} \ll \nu$.
\end{proof}

\begin{remark}
    Consider the Laakso diamond space $(X,d,\mu)$ from Example \ref{intro: Construction example} and let $\nu$ be as in the proof of Theorem \ref{thm:NoLips}. Then $(X,d,\nu)$ supports $(1,1)$-Poincar\'e inequality of Heinonen--Koskela \cite{HK}. This follows from Cheeger--Kleiner \cite{CK_PI}.
\end{remark}

\begin{remark}
    Theorem \ref{thm:NoLips} fails in general if any of the three conditions are absent. This follows from Examples \ref{example: Laakso space}, \ref{example: Counterexample} and \ref{ex: Singularities not equivalent} because Proposition \ref{prop: Smooth function on fibers} implies that the optimal potential functions in these examples are Lipschitz. Moreover, out of these three conditions we only know $d_{w,p} > p$ to be necessary for the conclusion of Theorem \ref{thm:NoLips} due to Theorem \ref{Thm: NS-equivalence}.
\end{remark}

\subsection{Non-trivial intersections of Sobolev spaces}
Let us get back to singularity of Sobolev spaces.
We shall discuss two examples where the singularity of Sobolev spaces does not happen. We begin with the Laakso space.

\begin{example}\label{ex:LaaksoNoSing}
    Consider the Laakso space in Example \ref{example: Laakso space} and see the computations in Figure \ref{fig:LaaksoSpaceNoSing}.
    \begin{figure}[!ht]
    \begin{tikzpicture}[scale=0.5]
    \input{plots/capacity_LaaksoSpaceNoSing}
    \end{tikzpicture} 
    \caption{Sobolev spaces on the Laakso space have non-trivial intersection.}
    \label{fig:LaaksoSpaceNoSing}
    \end{figure}
    Observe that the conditions \ref{item:A} and \ref{item:C} hold but $\ref{item:B}$ fails. 
    Indeed, the optimal potential $U_{p}$ is the same function for all $p \in (1,\infty)$.
    In what follows, the continuous optimal potential $\mathscr{U}_{\bullet}$ is contained in each $\mathscr{F}_p$, so we have a non-trivial intersection of Sobolev spaces. This also follows from the equality $\mathscr{F}_p = N^{1,p}(X,d,\mu)$ which by Hölder's inequality imply $\mathscr{F}_{p_1} \subseteq \mathscr{F}_{p_2}$ for $p_2 \in (1,p_1)$; see Theorem \ref{Thm: NS-equivalence}. Or we also could say $\Lip(X,d) \subseteq N^{1,p}(X,d,\mu) = \mathscr{F}_p$.
    
\end{example}

In the following example, we see that singularity of Sobolev spaces can fail in a far more interesting way compared to the previous example.

\begin{example}\label{ex: Singularities not equivalent}
Consider the IGS in Figure \ref{fig: singularities not equiv}.
We give a short argument why the energy measures are mutually singular in the sense of Proposition \ref{prop:SingularityEM} but the Sobolev spaces have a non-trivial intersection.
\begin{figure}[!ht]
    \includegraphics[width=150pt]{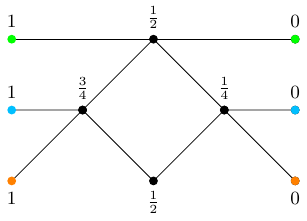}
    \caption{An example where the singularity of energy measures holds but singularity of Sobolev spaces fail.}
    \label{fig: singularities not equiv}
\end{figure}
First, the latter is easily seen to hold because the optimal potential $U_{p}$, which is computed in Figure \ref{fig: singularities not equiv} and the validity follows from $p$-harmonicity, is the same for all $p \in (1,\infty)$.
Hence $\mathscr{U}_{\bullet}$ is contained in each $\mathscr{F}_p$. On the other hand, the energy measures $\Gamma_{p_1}\Span{\mathscr{U}_{\bullet}}$ and $\Gamma_{p_2}\Span{\mathscr{U}_{\bullet}}$ are nevertheless different for every distinct $p_1,p_2 \in (1,\infty)$.
To see this, we simply apply the duality, Lemma \ref{lemma:duality} to show that $\mathcal{J}_{p_1} \neq \mathcal{J}_{p_2}$ for all distinct $p_1,p_2 \in (1,\infty)$. In conclusion, the singularity of energy measures holds while the singularity of Sobolev spaces fails. The reason why the proof of Theorem \ref{thm:SingularitySobolev} fails in this example is the lack of Proposition \ref{prop:EqualLittlemuu}.
\end{example}

In what follows from Example \ref{ex: Singularities not equivalent}, the two phenomena of singularity of energy measures and singularity of Sobolev spaces, are NOT equivalent.
This observation seems to indicate that the question of singularity of Sobolev spaces on the Sierpi\'nski gasket is not deducible from the singularity of energy measures which is known by \cite{KajinoSingularityPCF,KajinoSGsurvey}. There does not seem to be any reason to expect that a ``$p$-independent optimal flow'', meaning an analogue to condition \ref{item:C}, exists on the Sierpi\'nski gasket.

\subsection{Open problems}

Let us finish the paper by discussing some open problems.
We begin by reviewing two examples where we do not know whether the singularity of Sobolev space holds or not.

\begin{example}
    Consider the IGS in Figure \ref{fig:LaaksoDdeg}. Because $\abs{\nabla U_p(e)} = 0$ for the edge in the middle, the condition \ref{item:A} fails.
    Thus, the proof of Theorem \ref{thm:SingularitySobolev} does not work for this example.
    \begin{figure}[!ht]
    \begin{tikzpicture}
    \input{plots/LaaksoDDeg}
    \end{tikzpicture} 
    \caption{An example where the condition \ref{item:A} fails but \ref{item:B} and \ref{item:C} hold.}
    \label{fig:LaaksoDdeg}
\end{figure}
\end{example}

In the following example, the analysis whether singularity of Sobolev spaces holds or not seems to require a completely different approach.

\begin{example}
    For the IGS in Figure \ref{fig:p_independence} both the optimal potential $U_p$ and optimal flow $\mathcal{J}_p$ depend on $p$.
    \begin{figure}[!ht]
    \begin{tikzpicture}
    \input{plots/p_independence}
    \end{tikzpicture} 
    \caption{An example where both the optimal potentials and the optimal flows depend on $p$. The conditions \ref{item:A} and \ref{item:B} hold that \ref{item:C} fails.}
    \label{fig:p_independence}
\end{figure}
It is easy to check that $\Gamma_p\Span{\mathscr{U}_{p_1}}$ and $\Gamma_p\Span{\mathscr{U}_{p_2}}$ are not equal for distinct $p_1,p_2$. So according to Proposition \ref{prop:SingularityEM} singularity of energy measures holds. But since the flow depends on $p$, as we concluded in Example \ref{ex: Singularities not equivalent}, the approach used in the proof of Theorem \ref{thm:SingularitySobolev} is not available. We do not know whether singularity of Sobolev spaces holds for this example.
\end{example}

Lastly, it would be very interesting to clarify whether the two questions, singularity of energy measures and singularity of Sobolev spaces, are related. By the Laakso space in Example \ref{example: Laakso space}, both singularities may be false. By the Laakso diamond, both singularities may be true. According to Example \ref{ex:LaaksoNoSing}, it is possible that energy measures are singular while Sobolev spaces are not. However, we do not know about the remaining possibility.
To conclude the paper, we pose in a general context whether it can be realized. The precise interpretation of the question below is left to the reader.

\begin{question}\label{problem:Sing}
    Is there a natural family of $(1,p)$-Sobolev spaces and energy measures for $p \in (1,\infty)$ where singularity of Sobolev spaces holds while singularity of energy measures fails?
\end{question}

If the answer to Question \ref{problem:Sing} is positive, this would mean that the two problems of singularities of Sobolev spaces and are \emph{independent of each other}.

\bibliographystyle{acm}
\bibliography{clp}

\end{document}

%% file: plots/capacity_laakso.tex
\node at (-1,0) [circle,fill,inner sep=1.2pt]{};
\node at (3,0) [circle,fill,inner sep=1.2pt]{};
\node at (6,3) [circle,fill,inner sep=1.2pt]{};
\node at (6,-3) [circle,fill,inner sep=1.2pt]{};
\node at (9,0) [circle,fill,inner sep=1.2pt]{};
\node at (13,0) [circle,fill,inner sep=1.2pt]{};

\draw (-1,0)--(3,0)--(6,3)--(9,0)--(13,0);
\draw (3,0)--(6,-3)--(9,0);

% \node at (0,.5) {\small $x_0$}; \node at (3,.5) {\small $x_2$}; \node at (6,3.5) {\small $x_3$}; \node at (6,-3.5) {\small $x_4$}; \node at (9,.5) {\small $x_5$}; \node at (12,.5) {\small $x_1$};

\node at (4.3,0) {\small $1-\alpha_p$};
\node at (6,2.2) { $\frac{1}{2}$};
\node at (6,-2.2) { $\frac{1}{2}$};
\node at (8.2,0) {\small $\alpha_p$};
\node at (-1,.5) {\small $1$};
\node at (13,.5) {\small $0$};

\node at (-1,-1) {$v_+$};
\node at (13,-1) {$v_-$};

\node at (11,3) {$\alpha_p := \frac{2^{\frac{2 - p}{p-1}}}{2^{\frac{1}{p-1}} + 1}$};

%% file: plots/LaaksoDflow.tex
\node at (-1,0) [circle,fill,inner sep=1.5pt]{};
\node at (3,0) [circle,fill,inner sep=1.5pt]{};
\node at (6,3) [circle,fill,inner sep=1.5pt]{};
\node at (6,-3) [circle,fill,inner sep=1.5pt]{};
\node at (9,0) [circle,fill,inner sep=1.5pt]{};
\node at (13,0) [circle,fill,inner sep=1.5pt]{};

\draw (-1,0)--(3,0)--(6,3)--(9,0)--(13,0);
\draw (3,0)--(6,-3)--(9,0);

% \node at (0,.5) {\small $x_0$}; \node at (3,.5) {\small $x_2$}; \node at (6,3.5) {\small $x_3$}; \node at (6,-3.5) {\small $x_4$}; \node at (9,.5) {\small $x_5$}; \node at (12,.5) {\small $x_1$};

\node at (1,0.5) {\small $\overrightarrow{1}$};

\node at (4.25,2.25) {\small $\overrightarrow{\frac{1}{2}}$};

\node at (4.25,-2.25) {\small $\overrightarrow{\frac{1}{2}}$};

\node at (7.75,2.25) {\small $\overrightarrow{\frac{1}{2}}$};

\node at (7.75,-2.25) {\small $\overrightarrow{\frac{1}{2}}$};

\node at (11,0.5) {\small $\overrightarrow{1}$};

\node at (-1,-1) {$v_+$};
\node at (13,-1) {$v_-$};

%% file: plots/capacity_LaaksoSpaceNoSing.tex
\node at (0,3) [circle,fill,inner sep=1.5pt]{};
\node at (0,-3) [circle,fill,inner sep=1.5pt]{};
\node at (3,0) [circle,fill,inner sep=1.5pt]{};
\node at (6,3) [circle,fill,inner sep=1.5pt]{};
\node at (6,-3) [circle,fill,inner sep=1.5pt]{};
\node at (9,0) [circle,fill,inner sep=1.5pt]{};
\node at (12,3) [circle,fill,inner sep=1.5pt]{};
\node at (12,-3) [circle,fill,inner sep=1.5pt]{};

\draw (0,-3)--(3,0)--(6,3)--(9,0)--(12,3);
\draw (0,3)--(3,0)--(6,-3)--(9,0)--(12,-3);

\node at (0,-3.75) {\small $1$};
\node at (0,3.75) {\small $1$};
\node at (3,0.85) { $\frac{3}{4} $};
\node at (6,3.85) { $\frac{1}{2}$};
\node at (6,-3.85) { $\frac{1}{2}$};
\node at (9,0.85) { $\frac{1}{4}$};
\node at (12,-3.75) {\small $0$};
\node at (12,3.75) { \small $0$};

\node at (1.75,2) {\small $\overrightarrow{\frac{1}{2}}$};

\node at (1.75,-2) {\small $\overrightarrow{\frac{1}{2}}$};

\node at (3.9,2) {\small $\overrightarrow{\frac{1}{2}}$};

\node at (3.9,-2) {\small $\overrightarrow{\frac{1}{2}}$};

\node at (8,2) {\small $\overrightarrow{\frac{1}{2}}$};

\node at (8,-2) {\small $\overrightarrow{\frac{1}{2}}$};

\node at (10.25,2) {\small $\overrightarrow{\frac{1}{2}}$};

\node at (10.25,-2) { \small $\overrightarrow{\frac{1}{2}}$};

%% file: plots/LaaksoDDeg.tex
\draw (-.5,0)--(1,0)--(2,1)--(3,0)--(4.5,0);
\draw (1,0)--(2,-1)--(3,0);
\draw (2,1)--(2,-1);

\node at (-.5,0) [circle,fill,color = mycolor,inner sep=1.5pt]{};
\node at (1,0) [circle,fill,inner sep=1.2pt]{};
\node at (2,1) [circle,fill,inner sep=1.2pt]{};
\node at (2,-1) [circle,fill,inner sep=1.2pt]{};
\node at (3,0) [circle,fill,inner sep=1.2pt]{};
\node at (4.5,0) [circle,fill,color = mycolor,inner sep=1.5pt]{};

%% file: plots/p_independence.tex
\draw (0,0)--(3,0)--(5,-1)--(7,0)--(10,0);
\draw [-] (3,0) to [bend left = 5em] (7,0);

\node at (0,0) [circle,fill,color = mycolor,inner sep=1.5pt]{};
\node at (3,0) [circle,fill,inner sep=1.2pt]{};
\node at (5,-1) [circle,fill,inner sep=1.2pt]{};
\node at (7,0) [circle,fill,inner sep=1.2pt]{};
\node at (10,0) [circle,fill,color = mycolor,inner sep=1.5pt]{};

\node at (3,-.5) {\small $1-\alpha_p$};
\node at (5,-1.5) { $\frac{1}{2}$};
\node at (7,-0.5) {\small $\alpha_p$};
\node at (0,.5) {\small $1$};
\node at (10,.5) {\small $0$};

\node at (8,1.25) {\small $\alpha_p := \frac{(1 - 2^{1-p})^{\frac{1}{p-1}}}{1 + 2(1 - 2^{1-p})^{\frac{1}{p-1}}}$};